\documentclass[oneside, a4paper,reqno]{amsart}
\usepackage{pdfsync}
\usepackage{mathabx}
\usepackage{stmaryrd}
\usepackage{mathrsfs}
\usepackage{datetime}
\usepackage{chemarrow}
\usepackage{mathtools}

\usepackage{amsmath, amsthm, amscd, amssymb, latexsym, eucal, url}
\usepackage[all]{xy}

 \calclayout \makeatletter
 \calclayout \makeatletter
\def\serieslogo@{} \def\@setcopyright{} \makeatother

\makeatletter
\renewcommand*\env@matrix[1][c]{\hskip -\arraycolsep
  \let\@ifnextchar\new@ifnextchar
  \array{*\c@MaxMatrixCols #1}}
\makeatother

\numberwithin{equation}{section}
\newtheorem{thm}{Theorem}[section]
\newtheorem{cor}[thm]{Corollary}
\newtheorem{lem}[thm]{Lemma}
\newtheorem{prop}[thm]{Proposition}

\theoremstyle{definition}
\newtheorem{defn}[thm]{Definition}
\newtheorem{rem}[thm]{Remark}
\newtheorem{exam}[thm]{Example}

\newtheorem*{nota}{Notation}

\newcommand{\lxr}{\longrightarrow}
\newcommand{\rxr}{\longleftarrow}

\newcommand{\ab}{\mathfrak{Ab}}
\newcommand{\A}{\mathscr A}
\newcommand{\B}{\mathscr B}
\newcommand{\C}{\mathscr C}

\newcommand{\E}{\mathscr E}

\newcommand{\K}{\mathcal K}
\newcommand{\M}{\mathcal M}

\newcommand{\R}{\mathcal R}
\newcommand{\T}{\mathcal T}
\newcommand{\U}{\mathcal U}
\newcommand{\V}{\mathcal V}

\newcommand{\X}{\mathcal X}
\newcommand{\Y}{\mathcal Y}
\newcommand{\Z}{\mathcal Z}

\DeclareMathOperator*{\Ker}{\mathsf{Ker}}
 \DeclareMathOperator*{\Image}{\mathsf{Im}}
\DeclareMathOperator*{\Coker}{\mathsf{Coker}}

\DeclareMathOperator*{\coker}{\mathsf{coker}}

 \DeclareMathOperator{\pd}{\mathsf{pd}}
\DeclareMathOperator*{\id}{\mathsf{id}}

 \DeclareMathOperator*{\gd}{\mathsf{gl.dim}}

\DeclareMathOperator*{\Mod}{\mathsf{Mod}-\!}

\DeclareMathOperator*{\End}{\mathsf{End}}
 \DeclareMathOperator*{\smod}{\mathsf{mod}-\!}
 \DeclareMathOperator*{\lsmod}{\!-\mathsf{mod}}

\DeclareMathOperator*{\Gor}{\mathsf{G}-\!}
 \DeclareMathOperator*{\umod}{\underline{\mathsf{mod}}-\!}

\DeclareMathOperator*{\inj}{\mathsf{inj}}
\DeclareMathOperator*{\proj}{\mathsf{proj}}
\DeclareMathOperator*{\Inj}{\mathsf{Inj}}
\DeclareMathOperator*{\Proj}{\mathsf{Proj}}

\DeclareMathOperator*{\ddom}{\mathsf{dom.dim}}

\DeclareMathOperator*{\Fac}{\mathsf{Fac}}

\DeclareMathOperator*{\Add}{\mathsf{Add}}
 
\DeclareMathOperator*{\add}{\mathsf{add}}

\DeclareMathOperator{\Hom}{\mathsf{Hom}}

\DeclareMathOperator*{\Ext}{\mathsf{Ext}}

  \DeclareMathOperator*{\op}{\mathsf{op}}
  \DeclareMathOperator*{\silp}{\mathsf{silp}}
  \DeclareMathOperator*{\spli}{\mathsf{spli}}

        \DeclareMathOperator*{\uHom}{\underline{\mathsf{Hom}}}
         \DeclareMathOperator*{\uEnd}{\underline{\mathsf{End}}}
         
         \DeclareMathOperator*{\Cell}{{\mathsf{Cell}}}

\newcommand{\unA}{\underline A}
\newcommand{\unB}{\underline B}
\newcommand{\unC}{\underline C}

\newcommand{\unK}{\underline K}
\newcommand{\unX}{\underline X}

\newcommand{\unf}{\underline f}
\newcommand{\ung}{\underline g}

\newcommand{\un}{\underline}

\newsavebox{\proofbox}
\savebox{\proofbox}{\begin{picture}(7,7)%
  \put(0,0){\framebox(7,7){}}\end{picture}}

\begin{document}

\title[Rigid Localizations]{Rigid Objects, Triangulated Subfactors  and Abelian Localizations}

\author[A. Beligiannis]{Apostolos Beligiannis}
\address{Department of Mathematics, University of Ioannina, 45110
Ioannina, Greece}
\email{abeligia@cc.uoi.gr}



\dedicatory{Dedicated to Idun Reiten on the occasion of her 70th birthday}

\begin{abstract}
We show that the abelian category $\smod\X$ of coherent functors over a contravariantly finite rigid subcategory $\X$ in a triangulated category $\T$ is equivalent to the Gabriel-Zisman localization at the class of regular maps of a certain factor category of $\T$, and moreover it can be calculated by left and right fractions.  Thus we generalize recent results of Buan and Marsh.  We also extend recent results of Iyama-Yoshino concerning subfactor triangulated categories arising from mutation pairs in $\T$. In fact we give a classification of thick triangulated subcategories of a natural pretriangulated factor category of $\T$ and a classification of functorially finite rigid subcategories of $\T$ if the latter has Serre duality. In addition we characterize $2$-cluster tilting subcategories along these lines. Finally we extend basic results of Keller-Reiten concerning the Gorenstein and the Calabi-Yau property for categories arising from certain rigid, not necessarily cluster tilting, subcategories, as  well as several results of the literature concerning the connections between $2$-cluster tilting subcategories of triangulated categories  and tilting subcategories of the associated abelian category of coherent functors.        
\end{abstract}

\maketitle

\setcounter{tocdepth}{1} \tableofcontents

 \section{Introduction}
 
  Recently Keller and Reiten  in their study of cluster tilting theory in Calabi-Yau triangulated categories \cite{KR1}, generalizing previous work of Buan, Marsh and Reiten \cite{BRM}, observed that certain factor categories of triangulated categories arising from cluster tilting objects are abelian. This nice and rather rare phenomenon, which can be regarded as a tilting theorem in this context, provided the basis of many important developments in cluster tilting theory, and  was explored further by Koenig and Zhu \cite{KZ} and Iyama and Yoshino \cite{IY}.

 On the other hand Buan and Marsh, based on the work of Keller and Reiten, initiated the study of the structure and the localization theory of certain factor categories of triangulated categories with  Serre duality induced by rigid objects.     
 More precisely, in a series of two papers, see \cite{BM1, BM2}, Buan and Marsh showed that if $\T$ is a triangulated category with Serre duality over a field, and if $T$ is a {\em rigid} object of $\T$, i.e.  $\T(T,T[1]) = 0$, then the factor category $\T/T^{\bot}$ enjoys special exactness properties; here $T^{\bot}$ denotes the full subcategory of $\T$ consisting of all objects receiving no non-zero map starting from $T$.    Moreover they showed that the localization of $\T/T^{\bot}$, in the sense of Gabriel-Zisman \cite{GZ}, at the class of {\em regular maps} in $\T/T^{\bot}$, that is maps which are both monics and epics, is equivalent to the category of finitely generated modules over the endomorphism algebra of $T$. Finally they showed that this localization can be calculated by left and right fractions.  In their work Buan and Marsh used  results of Rump \cite{Rump} and others concerning the structure of semiabelian and related categories.   
 
 Our first aim in this paper is to give a simple short proof of the following more general result, see Theorem $4.6$ for more details. Recall that a full subcategory $\X$ of a triangulated category $\T$ is called {\em rigid} if $\T(\X,\X[1]) = 0$. We set $\X^{\bot} = \{A \in \T \, | \, \T(\X,A) = 0\}$ and we denote by $\smod\X$ the category of coherent functors over $\X$. 
 
 \medskip
 
 {\bf Theorem A.} {\em Let $\T$ be a triangulated category and $\X$ a contravariantly finite rigid subcategory of $\T$. Then the class $\mathcal R$ of regular morphisms in the factor category $\T/\X^{\bot}$ admits a calculus of left and right fractions, and the Gabriel-Zisman localized category $(\T/\X^{\bot})[\mathcal R^{-1}]$ is equivalent to the abelian category $\smod\X$:}
 \[
 \begin{CD}
  (\T/\X^{\bot})[\mathcal R^{-1}] \,\, @> \approx >> \,\, \smod\X 
 \end{CD}
 \]

 \medskip

Moreover the localization functor $\mathsf{P} \, \colon \, \T/\X^{\bot} \lxr (\T/\X^{\bot})[\mathcal R^{-1}]$ is faithful, preserves kernels and cokernels and  admits a fully faithful left adjoint. In addition we show that if the subcategory $\X^{\bot}$ is contravariantly finite in $\T$, then the factor category $\T/\X^{\bot}$ is integral, in the sense that $\T/\X^{\bot}$ has kernels and cokernels, and the class of epimorphisms, resp. monomorphisms, is closed under pull-backs, resp. push-outs.  As a consequence, for any morphism $f \colon A \lxr B$ in $\T/\X^{\bot}$, the canonical morphism $\widetilde{f} \colon \Coker(\mathsf{ker} f) \lxr \Ker(\coker f)$ is regular, so $\T/\X^{\bot}$ is semiabelian in the sense of \cite{Rump}, and the failure of $\T/\X^{\bot}$ to be abelian is measured by the failure of the  regular morphisms in $\T/\X^{\bot}$ to be invertible. An interesting case where $\X^{\bot}$ is contravariantly finite is when $\X$ is functorially finite and  $\T$ admits Serre duality. This is the situation studied by Buan and Marsh in \cite{BM1, BM2} in case $\X = \add T$ is the additive closure of a rigid object $T$ in $\T$, so Theorem A extends the results of Buan-Marsh to more general situations. 
However we give examples where Theorem A can not be covered by the methods of \cite{BM1, BM2}.  
 
 As an application of Theorem A, we characterize when the homological restricted Yoneda functor $\mathsf{H} \colon \T \lxr \smod\X$ given by $\mathsf{H}(A) = \T(-,A)|_{\X}$ is full. We show that this happens if and only if $\T$ consists of the direct factors of direct sums $A_{1} \oplus A_{2}$, where $A_{1}$ lies in the extension category $\X\star \X[1]$ and $A_{2}$ lies in $\X^{\bot}$ (compare with \cite[Theorem 3.3]{BK}). This is exactly the case when  the factor category $\T/\X^{\bot}$ is abelian. Using this characterization we give an example where $\mathsf{H}$ is not full and where the factor category $\T/\X^{\bot}$ is not abelian. Note that $\T/\X^{\bot}$ is abelian if $\X$ is $2$-cluster tilting in which case Keller-Reiten proved that the category $\smod\X$ is $1$-Gorenstein and stably $3$-Calabi-Yau (if $\T$ is $2$-Calabi-Yau), see \cite{KR1}. 

 In this connection we define a contravariantly finite rigid subcategory $\X$ of $\T$ to be {\em fully rigid} if the associated homological functor $\mathsf{H} \colon \T \lxr \smod\X$ is full. Then we have the following second main result of the paper which generalizes and improves the above mentioned basic results of Keller-Reiten and the results of Koenig-Zhu \cite{KZ}, see Theorems $5.2$, $5.4$. 
 
 \medskip
 
 {\bf Theorem B.} {\em Let $\X$ be a fully rigid subcategory of a Hom-finite  $k$-linear triangulated category  
 $\T$ over a field $k$.
 \begin{enumerate}
 \item If either $(\alpha)$ $\X = \add T$, or else $(\beta)$ $\X$ is covariantly finite and $\T$ has Serre duality, then the category $\smod\X$ of coherent functors over $\X$ is $1$-Gorenstein.  
\item If $\T$ is $2$-Calabi-Yau, then the stable triangulated category $\underline{\mathsf{CM}}(\X)$ of Cohen-Macaulay coherent functors over $\X$ is $3$-Calabi-Yau.
\end{enumerate}
}

 \medskip
 
Under the presence of a fully rigid subcategory $\X$ of a triangulated category $\T$ with Serre duality, we use  Theorem B to analyze the connections between $2$-cluster tilting subcategories of $\T$ and tilting subcategories of $\smod\X$.  In fact we show that $\X$ being $2$-cluster tilting is equivalent to the existence of a natural  bijective correspondence between $2$-cluster subcategories of $\T$ and special tilting subcategories of the category of coherent functors over $\X$. More precisely  we have the following third main result of the paper, see Theorem $6.6$,  which generalizes several related results of the literature, see \cite{Smith}, \cite{FuLiu}, \cite{IT}, \cite{JH}.

\medskip

 {\bf Theorem C.} {\em Let $\T$ be a Hom-finite  $k$-linear triangulated category over a field $k$ with a Serre functor $\mathbb S$, and let $\X$ be a fully rigid subcategory of $\T$ such that the category of coherent functors $\smod\X$ has finite global dimension. Then the following are equivalent.
\begin{enumerate}
\item $\X$ is a $2$-cluster tilting subcategory of $\T$.
\item $\mathbb S(\X) = \X[2]$, and the map $\Y \ \longmapsto \ \mathsf{H}(\Y)$ gives a bijective correspondence between:
\begin{enumerate}
\item[$\mathsf{(I)}$] $2$-cluster tilting subcategories $\Y$ of $\T$ with no non-zero direct summands from $\X^{\bot}$. 
\item[$\mathsf{(II)}$] Tilting subcategories $\U$ of $\smod\X$ such that:
\begin{enumerate}
\item[$(\alpha)$] $\mathsf{H}^{-1}(\U)$ is contravariantly finite in $\T$, and 
\item[$(\beta)$] $\mathbb S\mathsf{H}^{-1}(\U) = \mathsf{H}^{-1}(\U)[2]$. 
\end{enumerate}
\end{enumerate}
\end{enumerate}
}

\medskip
 
 It should be noted that the map $\mathsf{(II)} \longmapsto \mathsf{(I)}$ in part (ii) of Theorem C  exists without assuming finiteness of global dimension of $\smod\X$. As a consequence,  tilting subcategories of $\smod\X$ satisfying the conditions in $\mathsf{(II)}$ above lift to $2$-cluster tilting subcategories of $\T$. 
 
 In the same conceptual framework of cluster tilting theory, Iyama and Yoshino, working in the context of the theory of mutations, proved that certain factors of extension closed subcategories of triangulated categories are triangulated. More precisely they showed, see \cite[Theorem 4.2]{IY}, that if $\T$ is a triangulated category and $\U$ is an extension closed subcategory of $\T$ equipped with a full subcategory $\X$ such that $(\U,\U)$ forms an $\X$-mutation pair in the sense of \cite{IY}, then the subfactor category $\U/\X$ is triangulated.  This result played a critical role in the development of mutation theory of cluster tilting objects in general triangulated categories.  In this paper we call $(\U,\U)$ an $\X$-{\em mutation pair} if $(\U,\U)$ is an $\X$-mutation pair in the sense of Iyama-Yoshino satisfying the additional assumption that  $\U$ is closed under extensions in $\T$. 
 
Now recall from \cite{BR, B:3cats}, see also \cite{PJ}, that for any functorially finite subcategory $\X$ of a triangulated category $\T$, the factor category $\T/\X$, which itself is a Gabriel-Zisman localization of $\T$, is in a natural way a pretriangulated category. This means that there is an adjoint pair of endofunctors $(\Sigma^{1}_{\X},\Omega^{1}_{\X})$ on the factor category $\T/\X$, and $\T/\X$ carries a left triangulated structure with loop functor $\Omega^{1}_{\X}$ and a right triangulated structure with suspension functor $\Sigma^{1}_{\X}$ and the left and right triangulated structures are compatible in a certain way, see \cite{BR} for more details. Under reasonable assumptions, this construction carries over to extension closed subcategories of a triangulated category. First we need some definitions. We denote by  $\mathsf{Gh}_{\X}(\T)$, resp. $\mathsf{CoGh}_{\X}(\T)$, the ideal of $\T$ consisting of the $\X$-ghost, resp. $\X$-coghost, maps in $\T$, where a map $f$ is called $\X$-{\em ghost}, resp. $\X$-{\em coghost}, if the induced map $\T(\X,f)$, resp. $\T(f,\X)$, is zero. Let $\U$ be a full subcategory of $\T$ which is closed under extensions and direct summands, and let $\X$ be a functorially finite subcategory of $\U$. Then $\U$ is called $\X$-{\bf Frobenius} if: 
\begin{enumerate}
\item[$(\alpha)$] $\Omega^{1}_{\X}(\U) \subseteq \U \supseteq\Sigma^{1}_{\X}(\U)$, and
\item[$(\beta)$] $\mathsf{Gh}_{\X}(\U,\X[1]) = 0 = \mathsf{CoGh}_{\X}(\X[-1],\U)$. 
\end{enumerate}
In this context we have the following fourth main result of the paper, see Theorem $3.3$ and Corollary $3.11$,  which in particular gives a  classification of triangulated subfactors.

\medskip
 
 {\bf Theorem D.} {\em  Let $\T$ be a triangulated category. 
 \begin{enumerate} 
 \item[$(\alpha)$] If $\U$ is an $\X$-Frobenius subcategory of $\T$, then the subfactor category $\U/\X$ is triangulated.

 \item[$(\beta)$] For any functorially finite subcategory $\X \subseteq \T$,  the maps 
\[
\U \,\, \longmapsto \,\, \U/\X \,\,\,\,\,\,\,\,\, \text{and} \,\,\,\,\,\,\,\, \mathcal S \longmapsto   \pi^{-1}(\mathcal S)
\]
where $\pi \colon \T \lxr \T/\X$ is the projection functor, give mutually inverse bijections between:
\begin{enumerate}
\item[$\mathsf{(I)}$] $\X$-Frobenius subcategories of $\T$.  
\item[$\mathsf{(II)}$] Thick triangulated subcategories $\mathcal S$ of $\T/\X$, such that $\pi^{-1}(\mathcal S)$ is closed under extensions. 
\end{enumerate} 
\end{enumerate}
}

\medskip

Specializing Theorem D to rigid subcategories $\X$ we show that  $\U$ is $\X$-Frobenius if and only if $(\U,\U)$ forms an $\X$-mutation pair in the sense of Iyama-Yoshino, so  Theorem D  generalizes, and gives a simple proof to, Iyama-Yoshino's theorem.  Furthermore   if the vanishing condition $\T(\U,\U[-1]) = 0$ holds, then we show that the subfactor triangulated category $\U/\X$ is algebraic; in fact  $\U$ is an exact Frobenius category with $\X$ as the full subcategory of projective-injective objects. It should be noted that  the triangulated structure of the stable derived category of a self-injective algebra, in the sense of Wheeler \cite{Wheeler}, can be obtained from Theorem D. On the other hand specializing Theorem D to the case $\U = \T$ and assuming that $\T$ admits Serre duality, we recover a recent result of J{\o}rgensen \cite{PJ}  characterizing when the factor category $\T/\X$ is triangulated. This has some interesting applications to $n$-Calabi-Yau categories, for $ 0 \leq n\leq  2$.  

In the same vein we show that Theorem D can be used to give the following classification of functorially finite rigid subcategories $\X$ in a triangulated category $\T$ with Serre functor $\mathbb S$ satisfying $\mathbb S(\X) = \X[2]$, see Theorem $3.15$.  First note that if  $\U\subseteq \T$ satisfies $\mathbb S(\U) = \U[2]$, then ${^{\bot}}\U[1] = \U^{\bot}[-1]$.  Then we define  the {\em heart} of an extension closed subcategory $\U$ of $\T$ such that $\mathbb S(\U) = \U[2]$ to be the full subcategory $\U \cap {^{\bot}}\U[1] = \U \cap \U^{\bot}[-1]$, and we say that $\U$ is {\em maximal} if $\U$ is maximal among those with the same heart.

\medskip 

{\bf Theorem E.} {\em Let $\T$ be a triangulated category with a Serre functor $\mathbb S$ over a field $k$. Then the maps
\[
\U \ \ \longmapsto \ \ \X := \U \cap {^{\bot}}\U[1]  \ \ \ \ \  \text{and} \ \ \ \ \ \X \ \ \longmapsto \ \ \U := {^{\bot}}\X[1] 
\]
give mutually inverse bijections between:
\begin{enumerate}
\item[$\mathsf{(I)}$] Maximal functorially finite extension closed subcategories $\U$ of $\T$ such that $\mathbb S(\U) = \U[2]$.     
\item[$\mathsf{(II)}$] Functorially finite rigid subcategories $\X$ of $\T$ such that $\mathbb S(\X) = \X[2]$.   
\end{enumerate}
Under the above bijections, $\U$ is $\X$-Frobenius and the subfactor category $\U/\X$ is triangulated.}

\medskip

Note that Theorem E generalizes some results of Buan-Iyama-Reiten-Scott \cite{BIRS} on cluster structures and   can be regarded as a triangulated analog of the classification of Wakamatsu tilting modules over an Artin algebra by Mantese-Reiten \cite{MR}.

It should be mentioned that in case $\T$ has a Serre functor $\mathbb S$, in particular if $\T$ is $2$-Calabi-Yau, Theorems A, D and E combined together show that the theory of contra(co)variantly finite rigid subcategories $\X$ of $\T$ satisfying $\mathbb S(\X) = \X[2]$ is surprisingly  rich: on one hand the pretriangulated category $\T/\X$ contains $\X^{\bot}[-1]/\X$ as a maximal thick triangulated subcategory and on the other hand the factor pretriangulated category $\T/\X^{\bot}$ admits $\smod\X$ as an abelian localization. In this connection, we show further that if the pretriangulated category $\T/\X$ is abelian, then  the triangulated subfactor  category $\X^{\bot}[-1]/\X$ is semisimple abelian. Using this 
we have the sixth main result of the paper, see Theorem $7.3$, which, complementing and extending related results of Koenig-Zhu \cite{KZ}, characterizes $2$-cluster tilting subcategories in a triangulated category with Serre duality:

\medskip

{\bf Theorem F.} {\em Let $\T$ be a connected triangulated category with  a Serre functor $\mathbb S$ and let $\X$ be a contravariantly finite rigid subcategory of $\T$.  Then the following are equivalent. 
\begin{enumerate}
\item $\X$ is $2$-cluster tilting.
\item The pretriangulated category $\T/\X$ is abelian and $\mathbb S(\X) = \X[2]$.   
\end{enumerate} 
}

\medskip

It follows that for connected $2$-Calabi-Yau triangulated categories $\T$ the $2$-cluster tilting property of $\X$ is equivalent to the abelianness of $\T/\X$. In fact we show that $\X$ is a $2$-cluster tilting  if and only if the category $\smod\T/\X$ of coherent functors over $\T/\X$ is the free abelian, or Auslander, category of $\X$, in the sense of \cite{B:freyd}, so it has nice homological properties.  
     
 A general convention used in the paper is that the composition of morphisms in a given category, but not the composition of functors, is meant in the diagrammatic order. Our additive categories contain split idempotents and finite direct sums and their subcategories are assumed to be closed under isomorphisms and direct summands.

\section{Homological Finiteness, Coherent Functors and  Rigidity}

Throughout the paper $\T$ denotes a triangulated category with split idempotents. We fix a full additive subcategory $\X$ of $\T$ which is closed under isomorphisms and direct summands. 

Recall that $\X$ is called {\em contravariantly finite} in $\T$ if for any object $A$ in $\T$ there exists a map $f_{A} \colon X_{A} \lxr A$, where $X_{A}$ lies in $\X$ such that the induced map $\T(\X,f_{A}) \colon \T(\X,X_{A}) \lxr \T(\X,A)$ is surjective. In this case the map $f_{A}$ is called a {\em right} $\X$-{\em approximation} of $A$.  Dually we have the notions of {\em covariant finiteness} and {\em left approximations}.  Then $\X$ is called {\em functorially finite} if it is both contravariantly and covariantly finite. 

We denote by $\T/\X$ the stable or factor category of $\T$ with respect to $\X$; recall that the objects of $\T/\X$ are the objects of $\T$, and the morphism space $\Hom_{\T/\X}(A,B)$  is the quotient $\T(A,B)/\T_{\X}(A,B)$, where $\T_{\X}(A,B)$ is the subgroup of $\T(A,B)$ consisting of all maps $A \lxr B$ factorizing through an object from $\X$. We then have an additive functor $\pi \colon \T \lxr \T/\X$, $\pi(A) = \unA$ and $\pi(f) = \unf$. It should be noted that the functor $\pi$ represents $\T/\X$ as the localization $\T/\X = \T[\mathcal S^{-1}]$ of $\T$, in the sense of Gabriel-Zisman \cite{GZ}, at the class $\mathcal S$ of {\em stable equivalences}, i.e.  those maps $f \colon A \lxr B$ in $\T$ for which there exists a map $g \colon B \lxr A$ such that  both maps $1_{A} - g\circ f$ and $1_{B} - f\circ g$ factorize through an object from $\X$.    
For a class of objects  $\V \subseteq \T$ we denote by $\add \V$ the full subcategory of $\T$ consisting of the direct factors of finite direct sums of copies of objects from $\V$.

\subsection{Contravariant finiteness} If $\X$ is contravariantly finite in $\T$, then $\forall A \in \T$ there exists a triangle 
\begin{equation}
\begin{CD}
\Omega^{1}_{\X}(A) \,\ @> g^{0}_{A} >>  \,\, X^{0}_{A} \,\ @> f^{0}_{A} >> \,\, A \,\, @> h^{0}_{A} >>  \,\ \Omega^{1}_{\X}(A)[1]
\end{CD}
\end{equation}
where the middle map is a right $\X$-approximation of $A$. Iterating this construction we obtain triangles 
\[
\begin{CD}
\Omega^{t+1}_{\X}(A) \,\ @> g^{t}_{A} >>  \,\, X^{t}_{A} \,\ @> f^{t}_{A} >>  \,\, \Omega^{t}_{\X}(A) \,\, @> h^{t}_{A} >>  \,\ \Omega^{t+1}_{\X}(A)[1]
\end{CD}
\]
called the $\X$-{\em approximation triangles associated to} $A$, where the middle map is a right $\X$-approximation of $\Omega^{t}_{\X}(A)$, $\forall t \geq 0$, $\Omega^{0}_{\X}(A) := A$.  It is not difficult to see that the assignment $A \longmapsto \Omega^{1}_{\X}(A)$ induces a well-defined functor 
\[
\begin{CD}
\Omega^{1}_{\X} \,\, \colon \,\, \T/\X \,\, @> >>  \,\, \T/\X
\end{CD}
\]
on the stable category $\T/\X$, see \cite{B} for more details.  We denote by $\smod\X$ the category of coherent functors over $\X$. Recall that an additive functor $F \colon \X^{\op} \lxr \mathfrak{Ab}$ is called {\em coherent}, if there exists an exact sequence $\X(-,X^{1}) \lxr \X(-,X^{0}) \lxr F \lxr 0$. It is well known that $\smod\X$ is abelian if and only if $\X$ has weak kernels. An easy consequence of contravariant finiteness of $\X$ in $\T$ and the fact that $\T$, as a triangulated category, has weak kernels, is that $\smod\X$ is abelian. Moreover the restricted Yoneda functor  
\[
\begin{CD}
\mathsf{H} \, \colon \, \T \, @> >> \, \smod\X, \,\,\,\, \ \  \mathsf{H}(A) = \T(-,A)|_{\X}
\end{CD} 
\] 
is homological and induces an equivalence between $\X$ and the full subcategory of projective objects of $\smod\X$. We denote by $\X^{\bot}$ the full subcategory $\Ker \mathsf{H}$ of $\T$:
\[
\X^{\bot} = \big\{A \in \T \,\ | \,\, \T(\X,A) = 0\big\}
\]
 
 \subsection{Covariant finiteness} Dually let $\X$ be covariantly finite in $\T$. Then  $\forall A \in \T$ there exist triangles 
\begin{equation}
\begin{CD}
\Sigma^{t+1}_{\X}(A)[-1] \,\ @> h^{A}_{t} >> \,\, \Sigma^{t}_{\X}(A) \,\ @> f^{A}_{t} >>  \,\, X^{t}_{A} \,\, @> g^{A}_{t} >>  \,\ \Sigma^{t+1}_{\X}(A)
\end{CD}
\end{equation}
where the middle map is a left $\X$-approximation of $\Sigma^{t}_{\X}(A)$, $\forall t \geq 0$. As in $2.1$ the assignment $A \longmapsto \Sigma^{1}_{\X}(A)$ induces a well-defined functor $\Sigma^{1}_{\X} \, \colon \, \T/\X \, \lxr \, \T/\X$.  Since $\X$ is covariantly finite, the category $\X\lsmod := \smod\X^{\op}$ of covariant coherent functors over $\X$ is abelian and the restricted Yoneda functor  $
\mathsf{H}^{\op} \, \colon \, \T^{\op} \, \lxr \, \X\lsmod$, $\mathsf{H}^{\op}(A) = \T(A,-)|_{\X} 
$
is cohomological and induces a duality between $\X$ and the full subcategory of projective objects of $\X\lsmod$. In the sequel we shall use  the full subcategory $\Ker \mathsf{H}^{\op}$ of $\T$: \,\! ${^{\bot}}\X = \{A \in \T \,\ | \,\, \T(A,\X) = 0\}$.

\subsection{Functorial Finiteness} Now let $\X$ be  a functorially finite subcategory of $\T$. Then both functors $\Omega^{1}_{\X}, \Sigma^{1}_{\X} \colon$ $\T/\X \lxr \T/\X$ are defined and we have the following triangulated analog of \cite[Proposition 2.5]{B:gor}. For the notions of homotopy pull-backs/push-outs, in triangulated categories we refer to Neeman's book \cite[Chapter 1]{Neeman}.

 \begin{lem} There is  an adjoint pair
 \[
\big(\Sigma^{1}_{\X},\Omega^{1}_{\X}\big) \, \colon \, \xymatrix@C=2.75pc {\T/\X \,\, \ar@<0.5ex>[r]^-{{\mathsf{}}} & \,\, \ar@<0.5ex>[l]^-{{\,\, \mathsf{}}}\T/\X}
\]
where the unit $\delta : \mathsf{Id}_{\T/\X} \lxr \Omega^{1}_{\X}\Sigma^{1}_{\X}$ is epic and the counit $\varepsilon : \Sigma^{1}_{\X}\Omega^{1}_{\X} \lxr \mathsf{Id}_{\T/\X}$ is monic.
\begin{proof} Taking a left $\X$-approximation $f^{\Omega^{1}_{\X}(A)}_{0}  \colon \Omega^{1}_{\X}(A) \lxr X^{\Omega^{1}_{\X}(A)}_{0}$ of $\Omega^{1}_{\X}(A)$ and a right $\X$-approximation $f^{0}_{\Sigma^{1}_{\X}(A)} \colon$ $X^{0}_{\Sigma^{1}_{\X}(A)} \lxr \Sigma^{1}_{\X}(A)$ of $\Sigma^{1}_{\X}(A)$, and using the triangles $(2.1)$ and $(2.2)$,   we have morphisms of triangles:
\[
  \xymatrix{
    \Omega^{1}_{\X}(A) \ar[r]^{f^{\Omega^{1}_{\X}(A)}_{0}}\ar@{=}[d]&
    X^{\Omega^{1}_{\X}(A)}_{0} \ar[r]^{g^{\Omega^{1}_{\X}(A)}_{0}}\ar[d]^{\rho}&
    \Sigma^{1}_{\X}\Omega^{1}_{\X}(A) \ar[d]^{\varepsilon_{A}} \ar[r]^{-h^{\Omega^{1}_{\X}(A)}_{0}}&
    \Omega^{1}_{\X}(A)[1] \ar@{=}[d] \\
    \Omega^{1}_{\X}(A) \ar[r]^{g^{0}_{A}}& X^{0}_{A} \ar[r]^{f^{0}_{A}}& A\ar[r]^{h^{0}_{A}} & \Omega^{1}_{\X}(A)[1]}
    \]

and  
  \[\xymatrix{
    \Sigma^{1}_{\X}(A)[-1] \ar[r]^{h^{A}_{0}}\ar@{=}[d]&
    A \ar[r]^{f^{A}_{0}}\ar[d]^{\delta_{A}}&
    X^{0}_{A} \ar[d]^{\sigma} \ar[r]^{g^{A}_{0}}&
    \Sigma^{1}_{\X}(A)\ar@{=}[d] \\
    \Sigma^{1}_{\X}(A)[-1] \ar[r]^{-h^{0}_{\Sigma^{1}_{\X}(A)}}& \Omega^{1}_{\X}\Sigma^{1}_{\X}(A) \ar[r]^{g^{0}_{\Sigma^{1}_{\X}(A)}}& X^{0}_{\Sigma^{1}_{\X}(A)}\ar[r]^{f^{0}_{\Sigma^{1}_{\X}(A)}} & \Sigma^{1}_{\X}(A)}  
    \]
    
We leave to the reader as an easy exercise to check that the induced maps $\un\varepsilon_{A} : \Sigma_{\X}\Omega_{\X}(\unA) \lxr \unA$  and
$\un\delta_{A} : \unA \lxr \Omega^{1}_{\X}\Sigma^{1}_{\X}(\unA)$ are natural and define the counit and the unit of an adjoint pair $(\Sigma^{1}_{\X},\Omega^{1}_{\X})$ in $\T/\X$. If $\underline{\alpha} \colon \unC \lxr \Sigma^{1}_{\X}\Omega^{1}_{\X}(A)$ is a map in $\T/\X$ such that $\underline{\alpha} \circ \underline{\varepsilon}_{A} = 0$, then $\alpha \circ \varepsilon_{A}$ factorizes through an object from $\X$, hence it factorizes through the right $\X$-approximation $X^{0}_{A}$ of $A$.  Since the diagram on the left is homotopy pull-back, it follows that $\alpha$ factorizes through  $X^{\Omega^{1}_{\X}(A)}_{0}$, so $\underline{\alpha} = 0$ in $\T/\X$ and $\underline{\varepsilon}_{A}$ is monic. Dually $\underline{\delta}_{A}$ is epic. 
\end{proof}
\end{lem}

Note that the stable category $\T/\X$ is a pretriangulated category, see \cite{BR}, in the sense that $\T/\X$ is left triangulated with loop functor $\Omega^{1}_{\X}$ and right triangulated with suspension functor $\Sigma^{1}_{\X}$ and moreover the left and right triangulated structures are compatible in a certain way, we refer to \cite{BR} for more details. Clearly then $\T/\X$ is triangulated if and only if $\Sigma^{1}_{\X}$, or equivalently $\Omega^{1}_{\X}$, is an equivalence.  Here we only point out how the left and right triangles in $\T/\X$ are defined:

\begin{enumerate}
\item[\textsf{(LT)}] A left triangle in $\T/\X$ is by definition a diagram   which  is isomorphic in $\T/\X$ to a diagram 
\[
\begin{CD}
\Omega^{1}_{\X}(\unC) \, @> >> \, \unA \, @> >> \, \unB \, @> >> \, \unC 
\end{CD}
\eqno (l)
\]  
arising by forming the homotopy pull-back in $\T$ 
\[
  \xymatrix@C=1.5cm{
    \Omega^{1}_{\X}(C) \ar[r]^{}\ar@{=}[d]&
   A \ar[r]^{}\ar[d]^{}&
    B \ar[d]^{} \ar[r]^{}&
    \Omega^{1}_{\X}(C)[1] \ar@{=}[d] \\
    \Omega^{1}_{\X}(C) \ar[r]^{g^{0}_{C}}& X^{0}_{C} \ar[r]^{f^{0}_{C}}& C \ar[r]^{h^{0}_{C}\ \ \ \ } & \Omega^{1}_{\X}(C)[1]} 
\]
of a morphism $B \lxr C$ along the $\X$-approximation triangle $\Omega^{1}_{\X}(C) \lxr X^{0}_{C} \lxr C \lxr \Omega^{1}_{\X}(C)[1]$.
\item[\textsf{(RT)}] A right triangle in $\T/\X$ is by definition a diagram  which is isomorphic in $\T/\X$ to a diagram 
\[\begin{CD}
\unA \, @> >> \, \unB \, @> >> \, \unC \, @> >> \, \Sigma^{1}_{\X}(\unA) 
\end{CD}
\eqno (r)
\]  
arising by forming  the homotopy push-out in $\T$ 
\[
  \xymatrix@C=1.5cm{
    \Sigma^{1}_{\X}(A)[-1] \ar[r]^{ \ \ \ \ \ h^{A}_{0}}\ar@{=}[d]&
   A \ar[r]^{f^{A}_{0}}\ar[d]^{}&
    X^{A}_{0} \ar[d]^{} \ar[r]^{g^{A}_{0} \ \ }&
    \Sigma^{1}_{\X}(A) \ar@{=}[d] \\
    \Sigma^{1}_{\X}(A)[-1] \ar[r]^{}& B \ar[r]^{}& C \ar[r]^{} & \Sigma^{1}_{\X}(A)}
\]  
of a morphism $A \lxr B$ along the $\X$-approximation triangle $\Sigma^{1}_{\X}(A)[-1] \lxr A \lxr X^{A}_{0} \lxr \Sigma^{1}_{\X}(A)$. 
\item[$\bullet$] Equivalently it is easy to see  that the left, resp. right, triangles in $\T/\X$ are the diagrams which are isomorphic in $\T/\X$ to diagrams $(l)$, resp. $(r)$, arising from triangles $A \lxr B \lxr C \lxr A[1]$ in $\T$, where the map $\T(\X,B) \lxr \T(\X,C)$, resp. $\T(B,\X) \lxr \T(A,\X)$,  is epic. 
\end{enumerate} 

We always consider $\T/\X$ as a pretriangulated category with the above pretriangulated structure.

\subsection{Rigidity} A full subcategory $\X$ of $\T$ is called {\em rigid}, resp. {\em corigid}, if $\T(\X,\X[1]) = 0$, resp. $\T(\X,\X[-1]) = 0$. Clearly $\X$ is rigid if and only if $\X[1] \subseteq \X^{\bot}$ if and only if $\X[-1] \subseteq {^{\bot}}\X$.  In the sequel we shall need the following observation which gives an interesting interplay between rigidity and contra(co)variant finiteness.    

\begin{lem} Let $\X$ be a rigid subcategory of $\T$. 
\begin{enumerate}
\item If $\X$ is contravariantly finite, then $\X^{\bot}$ is covariantly finite and \ $\Omega^{t}_{\X}(A)[1] \in \X^{\bot}$, $\forall A\in \T$,  $\forall t \geq 1$.  
\item If $\X$ is covariantly finite, then ${^{\bot}}\X$ is contravariantly finite and  \ $\Sigma^{t}_{\X}(A)[-1] \in {^{\bot}}\X$,  $\forall A\in \T$,  $\forall t \geq 1$.   
\end{enumerate}
\begin{proof} (i) For any object $A$ in $\T$ and any integer $t \geq 0$, consider the approximation triangles 
\[
\begin{CD}
\Omega^{t+1}_{\X}(A) \,\ @> g^{t}_{A} >>  \,\, X^{t}_{A} \,\ @> f^{t}_{A} >>  \,\, \Omega^{t}_{\X}(A) \,\, @> h^{t}_{A} >>  \,\ \Omega^{t+1}_{\X}(A)[1]
\end{CD}
\]
where $\Omega^{0}_{\X}(A) := A$. Since the map $f^{t}_{A}$ is a right $\X$-approximation, the induced map $\T(\X,f^{t}_{A})$ is surjective, and since  $\X$ is rigid, we have $\T(\X,X^{t}_{A}[1]) = 0$. Then applying the homological functor $\T(\X,-)$ to the above triangles we infer that  $\T(\X,\Omega^{t}_{\X}(A)[1]) = 0$.  Hence $\Omega^{t}_{\X}(A)[1] \in \X^{\bot}$, $\forall A\in \T$,  $\forall t \geq 1$.  Now if $\alpha \colon A \lxr B$ is a map, where $B$ lies in $\X^{\bot}$, then the composition $f^{0}_{A} \circ \alpha$ is zero, hence $\alpha$ factorizes through $\Omega^{1}_{\X}(A)[1]$, i.e. the map $A \lxr \Omega^{1}_{\X}(A)[1]$ is a left $\X^{\bot}$-approximation of $A$.  We infer that $\X^{\bot}$ is covariantly finite.  Part (ii) is dual.  
\end{proof}
\end{lem}

\section{Triangulated Subfactors} 

Let as before $\T$ be a triangulated category with split idempotents. Our aim in this section is to  prove that certain subfactor categories of $\T$ carry, in a natural way, a triangulated structure which permits us to give classifications of functorially finite rigid or extension closed subcategories.  

\subsection{Triangulated Structures} Let $\X$ be a full subcategory of $\T$. Recall that a map $f \colon A \lxr B$ is called $\X$-{\bf ghost} if $\T(\X,f) = 0$, i.e. $\mathsf{H}(f) = 0$. Dually $f$ is called $\X$-{\bf coghost} if $\T(f,\X) = 0$, i.e. $\mathsf{H}^{\op}(f) = 0$. We denote by $\mathsf{Gh}_{\X}(A,B)$, resp. $\mathsf{CoGh}_{\X}(A,B)$, the subset of $\T(A,B)$ consisting of all $\X$-ghost, resp. $\X$-coghost, maps. Clearly $\mathsf{Gh}_{\X}(A,B)$ is a subgroup
of $\T(A,B)$ and it is easy to see that in this way we obtain an
ideal $\mathsf{Gh}_{\X}(\T)$ of $\T$. Dually we have the ideal $\mathsf{CoGh}_{\X}(\T)$ of $\X$-coghost maps.  

Assume that $\X$ is functorially finite in $\T$.   Let $A$ be an object of $\T$ and consider the approximation triangles
\begin{subequations}
\begin{equation}
  \Omega^{1}_{\X}(A) \,\ \stackrel{g^{0}_{A}}{\lxr} \,\, X^{0}_{A} \,\ \stackrel{f^{0}_{A}}{\lxr} \,\, A \,\, \stackrel{h^{0}_{A}}{\lxr} \,\ \Omega^{1}_{\X}(A)[1]
\end{equation}    
\begin{equation}
 \Sigma^{1}_{\X}(A)[-1] \,\ \stackrel{h^{A}_{0}}{\lxr} \,\, A \,\ \stackrel{f^{A}_{0}}{\lxr} \,\, X^{A}_{0} \,\, \stackrel{g^{A}_{0}}{\lxr} \,\ \Sigma^{1}_{\X}(A)
 \end{equation}
\end{subequations}

\,

Clearly the map $h^{0}_{A} \colon A \lxr \Omega^{1}_{\X}(A)[1]$ is $\X$-ghost and any $\X$-ghost map $A \lxr B$ factorizes through $h^{0}_{A}$, i.e. $h^{0}_{A}$ is a {\em universal $\X$-ghost map} out of $A$; dually $h^{A}_{0} \colon \Sigma^{1}_{\X}(A)[-1] \lxr A$ is a {\em universal $\X$-coghost map} into $A$.   

\begin{prop} Let $\X$ be a functorially finite subcategory of $\T$ and $A$ an object of $\T$. 
\begin{enumerate}
\item $\mathsf{Gh}_{\X}(A,\X[1]) = 0$ if and only if  the counit $\Sigma^{1}_{\X}\Omega^{1}_{\X}(A) \lxr A$ is invertible in $\T/\X$.
\item  $\mathsf{CoGh}_{\X}(\X[-1],A) = 0$ if and only if the unit $A \lxr \Omega^{1}_{\X}\Sigma^{1}_{\X}(A)$ is invertible in $\T/\X$.
\item The adjoint pair $(\Sigma^{1}_{\X}, \Omega^{1}_{\X})$ induces an equivalence: 
\[
\begin{CD}
\big\{A  \in \T \ | \ \mathsf{Gh}_{\X}(A,\X[1]) = 0\big\}/\X \ \ @> \approx >> \ \ \big\{A  \in \T \ | \ \mathsf{CoGh}_{\X}(\X[-1],A) = 0\big\}/\X 
\end{CD}
\]
\end{enumerate}
In particular the pretriangulated category $\T/\X$ is triangulated if and only if  for any object $A \in \T$:
\[\mathsf{Gh}_{\X}(A,\X[1]) = 0 = \mathsf{CoGh}_{\X}(\X[-1], A)\] 
\begin{proof} (i)  ``$\Longrightarrow$'' Assume that $\mathsf{Gh}_{\X}(A,\X[1]) = 0$. Since the map $h^{0}_{A} \colon A \lxr \Omega^{1}_{\X}(A)[1]$ is $\X$-ghost,  for any map $\Omega^{1}_{\X}(A) \lxr X$, where $X \in \X$, the composition $A \lxr \Omega^{1}_{\X}(A)[1] \lxr X[1]$ is $\X$-ghost and therefore it is zero. Hence the map $\Omega^{1}_{\X}(A)[1] \lxr X[1]$ factorizes through the cone $X^{0}_{A}[1]$ of $h^{0}_{A}$. Then $\Omega^{1}_{\X}(A) \lxr X$ factorizes through  $X^{0}_{A}$, and  this means that the map $g^{0}_{A} \colon \Omega^{1}_{\X}(A) \lxr X^{0}_{A}$ is a left $\X$-approximation of $\Omega^{1}_{\X}(A)$. By construction, see the proof of Lemma $2.1$,  this clearly implies that the counit $\varepsilon_{A} \colon \Sigma^{1}_{\X}\Omega^{1}_{\X}(A) \lxr A$ is invertible in $\T/\X$. 

``$\Longleftarrow$'' Assuming that the counit  $\Sigma^{1}_{\X}\Omega^{1}_{\X}(A) \lxr A$ is invertible in $\T/\X$, we first show that the map $g^{0}_{A} \colon \Omega^{1}_{\X}(A) \lxr X^{0}_{A}$ is a left $\X$-approximation of $\Omega^{1}_{\X}(A)$. Since $\varepsilon_{A} \colon \Sigma^{1}_{\X}\Omega^{1}_{\X}(A) \lxr A$ is invertible in $\T/\X$, there is a map $\mu_{A} \colon A \lxr \Sigma^{1}_{\X}\Omega^{1}_{\X}(A)$ such that the map $1_{A} - \mu_{A} \circ \varepsilon_{A}$ factorizes through an object of $\X$ and therefore it factorizes through the right $\X$-approximation $f^{0}_{A} \colon X^{0}_{A} \lxr A$ of $A$. Hence we have a factorization $1_{A} - \mu_{A} \circ \varepsilon_{A} = \kappa \circ f^{0}_{A}$ for some map $\kappa \colon A \lxr X^{0}_{A}$. Then $(1_{A} - \mu_{A} \circ \varepsilon_{A})\circ h^{0}_{A} = \kappa \circ f^{0}_{A} \circ h^{0}_{A} = 0$ and therefore $h^{0}_{A} = \mu_{A} \circ \varepsilon_{A} \circ h^{0}_{A}$. As in the proof of Lemma $2.1$ we have $\varepsilon_{A} \circ h^{0}_{A} = - h^{\Omega^{1}_{\X}(A)}_{0}$, hence $h^{0}_{A} = \mu_{A} \circ (- h^{\Omega^{1}_{\X}(A)}_{0})$.  It follows that there exists a map $\sigma \colon X^{0}_{A} \lxr X^{\Omega^{1}_{\X}(A)}_{0}$ making the following diagram a morphism of triangles:
\[
\xymatrix@C=1.5cm{
  \Omega^{1}_{\X}(A) \ar[r]^{g^{0}_{A}} \ar@{=}[d] & X_A^0 \ar@{-->}[d]^{\sigma} \ar[r]^{f^{0}_{A}} & A \ar[d]^{\mu_{A}} \ar[r]^{h^{0}_{A}} & \Omega^{1}_{\X}(A)[1] \ar@{=}[d]  \\
  \Omega^{1}_{\X}(A) \ar[r]^{f^{\Omega^{1}_{\X}(A)}_{0}} & X^{\Omega^{1}_{\X}(A)}_{0} \ar[r]^{g^{\Omega^{1}_{\X}(A)}_{0}} & \Sigma^{1}_{\X}\Omega^{1}_{\X}(A) \ar[r]^{-h^{\Omega^{1}_{\X}(A)}_{0}} & \Omega^{1}_{\X}(A)[1]   }
\]  
If $\rho \colon \Omega^{1}_{\X}(A) \lxr X$ is a map, where $X \in \X$, then $\rho$ factorizes through the left $\X$-approximation $f^{\Omega^{1}_{\X}(A)}_{0}$ of $\Omega^{1}_{\X}(A)$. Since from the above commutative diagram we have $f^{\Omega^{1}_{\X}(A)}_{0} = g^{0}_{A} \circ \sigma$, we infer that $\rho$ factorizes through $g^{0}_{A}$. This means that  the map $g^{0}_{A} \colon \Omega^{1}_{\X}(A) \lxr X^{0}_{A}$ is a left $\X$-approximation of $\Omega^{1}_{\X}(A)$. 

Finally let $\alpha \colon A \lxr X[1]$ be an $\X$-ghost map, where $X \in \X$. Then clearly $\alpha = h^{0}_{A} \circ \beta$ for some map $\beta \colon \Omega^{1}_{\X}(A)[1] \lxr X[1]$. Since $\beta[-1] \colon \Omega^{1}_{\X}(A) \lxr X$ factorizes through the left $\X$-approximation $g^{0}_{A}$ of $\Omega^{1}_{\X}(A)$,  we have $\alpha[-1] = h^{0}_{A}[-1] \circ g^{0}_{A} \circ \tau$ for some map $\tau \colon X^{0}_{A} \lxr X$. Since $h^{0}_{A}[-1] \circ g^{0}_{A} = 0$, it follows that  $\alpha[-1]$, hence $\alpha$, is zero. We infer that $\mathsf{Gh}_{\X}(A,X[1]) = 0$.    

Part (ii) is dual, and part (iii) and the last assertion follow directly from (i), (ii).  
\end{proof} 
\end{prop}

The above result naturally suggests the following definition. 

\begin{defn} Let $\X$ be a full subcategory of $\T$. A full subcategory $\U$ of $\T$ is called $\X$-{\bf Frobenius} if:
\begin{enumerate}
\item $\U$ closed under extensions and direct summands in $\T$.  
\item $\X \subseteq \U$ is functorially finite in $\U$. 
\item \begin{enumerate}
\item $\mathsf{Gh}_{\X}(\U,\X[1]) = 0 = \mathsf{CoGh}_{\X}(\X[-1], \U)$.
\item $\Omega^{1}_{\X}(\U) \subseteq \U \supseteq\Sigma^{1}_{\X}(\U)$.
\end{enumerate} 
\end{enumerate}
\end{defn}

Clearly if $\X = \{0\}$, then $\T$ is $\X$-Frobenius, and if $\X = \T$, then $\T$ is $\X$-Frobenius if and only if $\T = \{0\}$. In the first case we have the triangulated factor $\T/\X = \T$ and in the second case the triangulated factor $\T/\X = \{0\}$. 
The following result, which is Theorem D from the Introduction,  characterizes along these lines when the stable category of an extension closed subcategory of $\T$ carries a natural triangulated structure. Recall first that a full triangulated subcategory of a (pre-)triangulated category is called {\em thick} if it is closed under direct factors.

\begin{thm} \begin{enumerate}
\item[$\mathrm{(i)}$] The stable category $\U/\X$ of an $\X$-Frobenius subcategory $\U$ of $\T$ is triangulated.  
\item[$\mathrm{(ii)}$] If $\X$ is a functorially finite subcategory of $\T$, then the maps 
\[
\U \,\, \longmapsto \,\, \U/\X \,\,\,\,\,\,\,\,\, \text{and} \,\,\,\,\,\,\,\, \mathcal S \,\, \longmapsto \,\,    \pi^{-1}(\mathcal S)
\]
give mutually inverse bijections between:
\begin{enumerate}
\item[$\mathsf{(I)}$] $\X$-Frobenius subcategories $\U$ in $\T$.
\item[$\mathsf{(II)}$] Thick triangulated subcategories $\mathcal S$ of $\T/\X$ such that $\pi^{-1}(\mathcal S)$ is closed under extensions in $\T$.   
\end{enumerate} 
\end{enumerate}
\begin{proof} (i)  Since $\U$ is $\X$-Frobenius, the functors $\Sigma^{1}_{\X}, \Omega^{1}_{\X} \colon \U/\X \lxr \U/\X$ are defined and then as noted in $2.3$ we have an adjoint pair $
\big(\Sigma^{1}_{\X},\Omega^{1}_{\X}\big) \, \colon  \xymatrix@C=1.5pc {\U/\X \ \ar@<0.4ex>[r]^-{{\mathsf{}}} &  \ar@<0.4ex>[l]^-{{\mathsf{}}} \ \U/\X}$. Recall that the left and right triangles in $\T/\X$ are defined using homotopy pull-backs and homotopy push-outs respectively. Since $\U$ is closed under extensions, it follows easily that the proof of the fact that $\T/\X$ is pretriangulated when $\X$ is functorially finite in $\T$ works if we replace $\T$ with any extension closed subcategory $\U$ of $\T$ containing $\X$ as a functorially finite subcategory satisfying condition (ii)(b) of Definition $3.2$, see \cite{B:3cats, B, PJ}.  Hence $\U/\X$ is a pretriangulated category and it remains to show that $\Omega^{1}_{\X}$, or equivalently $\Sigma^{1}_{\X}$, is an equivalence.  This follows from Proposition $3.1$. 

(ii) If  $\U$ is $\X$-Frobenius, then by (i) the subfactor category $\U/\X$ is triangulated, and in fact it is a triangulated subcategory of $\T/\X$.  If $\underline{U} = \unA \oplus \unB$ is a direct sum decomposition in $\T/\X$, where $U \in\U$, then there exist objects $X_{1}, X_{2}$ in $\X$ and a direct sum decomposition $U \oplus X_{1} \cong A \oplus B \oplus X_{2}$ in $\T$. Since $\U$ is closed under direct summands and contains $\X$, the objects $A$ and $B$ lie in $\U$, i.e. $\unA, \unB \in \U/\X$. We infer that $\mathcal S := \U/\X$ is a thick triangulated subcategory of $\T/\X$. Conversely, let $\mathcal S$ be a thick triangulated subcategory of $\T/\X$ such that $\U := \pi^{-1}(\mathcal S)$ is closed under extensions in $\T$. Then $\U$ is a full subcategory of $\T$ containing $\X$ and clearly $\U$ is closed under direct summands in $\T$.  Since $\X$ is functorially finite in $\T$ and $\mathcal S$ is a triangulated subcategory of $\T/\X$, it follows that $\Omega^{1}_{\X}(\U) \subseteq \U \supseteq\Sigma^{1}_{\X}(\U)$. Proposition $3.1$ shows that the remaining condition (iii)(a) of Definition $3.2$ holds and therefore  $\U$ is $\X$-Frobenius.  Clearly the maps $\mathsf{(I)} \mapsto \mathsf{(II)}$ and $\mathsf{(II)} \mapsto \mathsf{(I)}$ are mutually inverse.  
\end{proof}
\end{thm}  

\begin{exam} Let $\T$ be a triangulated category with enough $\E$-projectives and enough $\E$-injectives for a proper class of triangles $\E$ in $\T$ in the sense of \cite{B:3cats}. Let  $\mathcal P(\E)$, resp. $\mathcal I(\E)$, be the full subcategory of $\E$-projective, resp. $\E$-injective, objects.  We say that $\T$ is $\E$-{\em Frobenius} if $\mathcal P(\E) = \mathcal I(\E)$.  Then by Proposition $3.1$ we have the following: $\T$ is $\E$-Frobenius if and only if the factor category $\T/\mathcal P(\E)$, equivalently the factor category $\T/\mathcal I(\E)$, is triangulated if and only if $\T$ is $\mathcal P(\E)$-Frobenius if and only if $\T$ is $\mathcal I(\E)$-Frobenius in the sense of  Definition $3.2$, see also \cite[Theorem 7.2]{B:3cats}. 

For instance  the bounded derived category ${\bf D}^{b}(\smod\Lambda)$ of a self-injective algebra $\Lambda$ is $\X$-Frobenius, where $\X$ is the full subcategory of ${\bf D}^{b}(\smod\Lambda)$ consisting of the complexes with projective components and zero differential. Hence the stable category $\underline{{\bf D}}^{b}(\smod\Lambda) := {\bf D}^{b}(\smod\Lambda)/\X$, known as the  stable derived category of $\Lambda$ (in the sense of Wheeler \cite{Wheeler}), is triangulated; we refer to \cite{Wheeler}, \cite{B:3cats} for more details.      
\end{exam}

If $\X$ is rigid, then clearly $\mathsf{Gh}_{\X}(A,\X[1]) = \T(A,\X[1])$ and $\mathsf{CoGh}_{\X}(\X[-1],A) = \T(\X[-1],A) = \T(\X,A[1])$. This implies that there are no non-trivial triangulated factors of $\T$ by rigid functorially finite subcategories:

\begin{cor} Let $\X$ be a functorially finite rigid subcategory of $\T$.
\begin{enumerate}
\item $\T(A,\X[1]) = 0$ if and only if the counit $\Sigma^{1}_{\X}\Omega^{1}_{\X}(A) \lxr A$ is invertible in $\T/\X$. 
\item $\T(\X,A[1]) = 0$ if and only if  the unit $A \lxr \Omega^{1}_{\X}\Sigma^{1}_{\X}(A)$ is invertible in $\T/\X$.
\item ${^{\bot}}\X[1]/\X = \big\{\unA \in \T/\X \ | \ \text{the counit}  \ \Sigma^{1}_{\X}\Omega^{1}_{\X}(\unA) \lxr \unA \ \text{is invertible}\big\}$.
\item $\X^{\bot}[-1]/\X = \big\{\unA \in \T/\X \ | \ \text{the unit}  \ \unA \lxr \Omega^{1}_{\X}\Sigma^{1}_{\X}(\unA) \ \text{is invertible}\big\}$.
\item The adjoint pair $(\Sigma^{1}_{\X}, \Omega^{1}_{\X})$ induces an equivalence: ${^{\bot}}\X[1]/\X \approx \X^{\bot}[-1]/\X$. 
\end{enumerate}
In particular $\T/\X$ is triangulated, i.e. $\T$ is $\X$-Frobenius,  if and only if $\X = 0$. 
\end{cor}

A $k$-linear triangulated category $\T$ over a field $k$ is called {\em Hom-finite} if the $k$-vector space $\T(A,B)$ is finite-dimensional, $\forall A,B \in \T$. A Hom-finite category $\T$ admits {\em Serre duality} if there is a triangulated equivalence $\mathbb S \colon \T \lxr \T$, the {\em Serre functor}, equipped with natural bifunctorial  isomorphisms: 
\[
\begin{CD}
\mathsf{D}\T(A,B)  \,\, @> \approx >> \,\, \T\big(B,\mathbb S(A)\big), \,\,\,\,\,\,\ \forall A, B \in \T
\end{CD}
\]
where $\mathsf{D} = \Hom_{k}(-,k)$ denotes duality with respect to the base field. If $\mathbb S$ is a Serre functor of $\T$, then $\T$ is called $n$-{\em Calabi-Yau}, for some $n \geq 1$, if $\mathbb S(?) \cong (?)[n]$ as triangulated functors.   

The next result, part (iii) of which  recovers and gives a short proof to a result of J\o rgensen, see \cite[Theorem 2.3]{PJ}, shows that if $\X$ is not rigid but $\T$ admits Serre duality, then the vanishing conditions in Proposition $3.1$ take a nice form. It follows that  $\T$ is $\X$-Frobenius if and only if the Serre functor $\mathbb S$ satisfies $\mathbb S(\X) = \X[1]$.

\begin{cor} Let $\T$ be a Hom-finite $k$-linear triangulated category over a field $k$ which admits a Serre functor $\mathbb S$, and let $\X$ be  a  functorially finite subcategory of $\T$. 
\begin{enumerate}
\item $\mathsf{Gh}_{\X}(-,\X[1]) = 0$ \,  if and only if \, $\X[1] \subseteq \mathbb S(\X)$. 
\item $\mathsf{CoGh}_{\X}(\X[-1],-) = 0$ \, if and only if \, $\mathbb S(\X) \subseteq \X[1]$. 
\item $\T/\X$ is triangulated \, if and only if \,  $\mathbb S(\X)  = \X[1]$. 
\end{enumerate}
\begin{proof} (i) ``$\Longleftarrow$'' Let $f \in \mathsf{Gh}_{\X}(A,X[1])$ be $\X$-ghost, where $X \in \X$.  Since $\X[1] \subseteq \mathbb S(\X)$, $X[1] = \mathbb S(X^{\prime})$ for some object $X^{\prime} \in \X$ and we have the $\X$-ghost map $f \colon A \lxr \mathbb S(X^{\prime})$, i.e. $\T(\X,f) \colon \T(\X,A)$ $\lxr \T(\X,\mathbb S(X^{\prime}))$ is zero. Using Serre duality this easily implies that  the map $\T(X^{\prime},\X) \lxr \T(\mathbb S^{-1}(A),\X)$ is zero. The image of $1_{X^{\prime}}$ under this map is the map  $\mathbb S^{-1}(f) \colon \mathbb S^{-1}(A) \lxr X^{\prime}$. Hence $\mathbb S^{-1}(f) = 0$ and then $f = 0$. Therefore $\mathsf{Gh}_{\X}(A,\X[1]) = 0$. 

``$\Longrightarrow$'' Let $0 \neq X \in \X$ be an indecomposable object; then the identity map of $X$, via  Serre duality $\mathsf{D}\T(X,X) \stackrel{\cong}{\lxr} \T(X,\mathbb S(X)) \stackrel{\cong}{\lxr} \T(\mathbb S^{-1}(X)[1],X[1])$, gives us a map $\mathbb S^{-1}(X)[1] \lxr X[1]$ which induces an Auslander-Reiten triangle $(*) \colon X \lxr A \lxr \mathbb S^{-1}(X)[1] \lxr X[1]$ in $\T$ in the sense of Happel \cite{Happel}, see also \cite{VdBR}.  If $\mathbb S^{-1}(X)[1]$ does not lies in $\X$, then any map $X^{\prime} \lxr \mathbb S^{-1}(X)[1]$, where $X^{\prime} \in \X$, is not split mono and therefore it factorizes through $A$; equivalently the composition $X^{\prime} \lxr \mathbb S^{-1}(X)[1] \lxr X[1]$ is zero. This means that the map $\mathbb S^{-1}(X)[1] \lxr X[1]$ is $\X$-ghost and therefore it is zero. Then $X \lxr \mathbb S(X)$ is zero and this is impossible since $X\neq 0$. Hence $\mathbb S^{-1}(X)[1] = \mathbb S^{-1}(X[1])$ lies in $\X$ and therefore $X[1] \in \mathbb S(\X)$ for any indecomposable object $X$ of $\X$. Since $\T$ is a Krull-Schmidt category we infer that $\X[1] \subseteq \mathbb S(\X)$. 

Part (ii) is dual and (iii) follows from (i), (ii) and Proposition $3.1$. 
\end{proof}
\end{cor}  

If $\X$ is a functorially finite subcategory of $\T$ and the factor category $\T/\X$ is triangulated, then the following observation gives an alternative description of $\T/\X$ if, in addition, $\X$ is a thick subcategory of $\T$. To avoid confusion between the two factors we denote by $\T\! \sslash \! \X$ the Verdier quotient of $\T$ by a thick subcategory $\X$ of $\T$.

\begin{lem} Let $\T$ be a triangulated category and $\X$ a functorially finite subcategory of $\T$.   If the factor category $\T/\X$ is triangulated and $\X$ is a thick subcategory of $\T$, then there is a triangle equivalence: 
\[
\begin{CD}
\T/\X  \,\, @> \approx >> \,\, \T\!\sslash \! \X
\end{CD}
\]
\begin{proof} If $\X$ is thick, then we have two triangulated categories: the factor category $\T/\X$ and the Verdier quotient $\T\!\sslash \! \X$. Recall that $\T/\X$  is the localization $\T[\mathcal S^{-1}]$ of $\T$ at the class  $\mathcal S$ of stable equivalences, and $\T\!\sslash\!\X$ is the localization $\T[\mathcal R^{-1}]$ of $\T$ at the class $\mathcal R$ of those maps  in $\T$ whose cone lies in $\X$. Since the objects of $\X$ become zero in $\T\!\sslash\!\X$, trivially any stable equivalence lies in $\mathcal R$, hence $\mathcal S \subseteq \mathcal R$. Now let $\sigma \colon A \lxr B$ be in $\R$ and let $X \lxr A \lxr B \lxr X[1]$ be a triangle in $\T$, where $X$, hence $X[1]$, lies in $\X$. Taking the homotopy pull-back of $\sigma$ along a right $\X$-approximation $X^{0}_{B}$ of $B$, we have, as in subsection 2.3, a triangle $\Omega^{1}_{\X}(\unB) \lxr \unC \lxr \underline{A} \lxr \unB$ in $\T/\X$ and a triangle $X \lxr C \lxr X^{0}_{B} \lxr X[1]$ in $\T$. Since $\X$ is closed under extensions, it follows that $C$ lies in $\X$. Since $\T/\X$ is triangulated, this implies that $\unC = 0$ and therefore $\underline{\sigma}$ is invertible in $\T/\X$. Then $\sigma$ lies in $\mathcal S$ and therefore $\R = \mathcal S$. 
This implies that $\T/\X = \T[\mathcal S^{-1}] = \T[\mathcal R^{-1}] = \T\!\sslash \X$. 
\end{proof}
\end{lem} 

As a direct consequence we have the following. 

\begin{cor} Let $\T$ be a  $k$-linear triangulated category and $\X$ a functorially finite subcategory of $\T$.  
\begin{enumerate} 
\item If $\T$ is $1$-Calabi-Yau, then the factor category $\T/\X$ is triangulated.
\item If $\T$ is $n$-Calabi-Yau, where $n = 0, 2$, then the factor category $\T/\X$ is triangulated  if and only if $\X = \X[1]$.
\item If $\T$ is $n$-Calabi-Yau, where $0 \leq n \leq 2$, and $\X$ is a thick subcategory of $\T$, then there is a triangle equivalence $\T/\X \approx \T\!\sslash \! \X$, where $ \T\!\sslash \! \X$ is the Verdier quotient. 
\end{enumerate}
\begin{proof} Parts (i) and (ii) follow from Corollary $3.6$ and part (iii) follows from (i), (ii) and Lemma  $3.7$.  
\end{proof} 
\end{cor}

\begin{exam} (i) Let $\T$ be either $(\alpha)$ the category $\proj \Lambda$ of finitely generated projective modules over a deformed preprojective
algebra $\Lambda$ of generalized Dynkin type, or $(\beta)$ the stable module category $\umod k[t]/(t^{n})$, for $n \geq 3$, of the self-injective algebra $k[t]/(t^{n})$.  Then $\T$ is a $1$-Calabi-Yau triangulated category of finite representation type, see \cite{Amiot} for $(\alpha)$ and  \cite{ES} for $(\beta)$. Hence for any subcategory $\X$ of $\T$, the stable category $\T/\X$ is triangulated. 

(ii)  The bounded derived category ${\bf D}^{b}(\mathsf{Coh}\,\mathbb E)$ of coherent sheaves over an elliptic curve $\mathbb E$ is $1$-Calabi-Yau. Hence for any functorially finite subcategory $\X$ of ${\bf D}^{b}(\mathsf{Coh}\,\mathbb E)$, the factor category ${\bf D}^{b}(\mathsf{Coh}\,\mathbb E)/\X$  is triangulated. 

(iii) It follows from results of Rickard, see \cite{Rickard}, that the category $\K^{b}(\proj\Lambda)$ of perfect complexes over a symmetric finite dimensional $k$-algebra $\Lambda$ over a field $k$ is $0$-Calabi-Yau. Hence, by Corollary $3.8$, for any functorially finite subcategory $\X \subseteq \K^{b}(\proj\Lambda)$, the factor category $\K^{b}(\proj\Lambda)/\X$ is triangulated if and only if $\X  = \X[1]$.    
\end{exam}

\begin{rem} Assume that $\T$ admits a Serre functor $\mathbb S$ and $\X$ is a full, not necessarily rigid, subcategory of $\T$. If $\X^{\bot}$ is contravariantly finite and ${^{\bot}}\X$ is covariantly finite, then  since, as easily seen, $\mathbb S({^{\bot}}\X) = \X^{\bot}$, it follows that both $\X^{\bot}$ and ${^{\bot}}\X$ are functorially finite and we have an equivalence $\T/{^{\bot}}\X \approx \T/\X^{\bot}$. Then by Corollary $3.6$, the factor category $\T/\X^{\bot} \approx \T/{^{\bot}}\X$ is triangulated if and only if $\mathbb S({^{\bot}}\X) = {^{\bot}}\X[1]$; this is equivalent to: $\X^{\bot}[-1] = {^{\bot}}\X$.     
\end{rem}

\subsection{Mutation Pairs} Recently Iyama and Yoshino in their study of mutations of cluster tilting objects in triangulated categories constructed in \cite[Theorem 4.2]{IY} subfactor triangulated categories out of  $\X$-mutation pairs in the following sense: 

\begin{defn} (See \cite[Sections 2 and 4]{IY}) Let $\X$ be a full subcategory of $\T$. For a full subcategory $\U$ of $\T$,  we say that $(\U,\U)$ forms an $\X$-{\em mutation pair} if the following conditions hold:
\begin{enumerate}
\item $\U$ is closed under extensions and direct summands in $\T$, and contains $\X$. 
\item  For any object $U \in \U$, there exist objects $X^{\pm} \in \X$ and $U^{\pm} \in \U$ and  triangles in $\T$:
 \[ U^{+} \ \lxr \  X^{+} \ \lxr \ U \ \lxr \ U^{+}[1] \,\,\,\,\,\,\,\,\,\, \text{and}  \,\,\,\,\,\,\,\,\,\, U^{-}[-1] \ \lxr \ U \ \lxr \ X^{-} \  \lxr \ U^{-}  
 \]
\item $\T(\U,\X[1]) = 0 = \T(\X,\U[1])$.  
\end{enumerate} 
\end{defn} 

It is easy to see that $(\U,\U)$ is an $\X$-mutation pair in the sense of $3.11$ iff $(\U,\U)$ is an $\X$-mutation pair as defined by Iyama-Yoshino in \cite[Section 2]{IY} and $\U$ is closed under extensions. Now if $(\U,\U)$ is an $\X$-mutation pair in $\T$, then clearly the map $X^{+} \lxr U$ is a right $\X$-approximation of $U$ and $\Omega^{1}_{\X}(U) = U^{+} \in \U$, and the map $U \lxr X^{-}$ is a left $\X$-approximation of $U$ and $\Sigma^{1}_{\X}(U) = U^{-} \in \U$. Since $\X$ is rigid, we infer that $\U$ is an $\X$-Frobenius subcategory of $\T$, and then Theorem $3.3$ reduces to the following result which generalizes \cite[Theorem 4.2]{IY} giving at the same time  a classification of all $\X$-mutation pairs in $\T$  in terms of special thick triangulated subcategories of $\T/\X$.

\begin{cor} \begin{enumerate}
\item If $\X$ and  $\U$ are full subcategories of $\T$, then the following conditions are equivalent:
\begin{enumerate}
\item $(\U,\U)$ is an $\X$-mutation pair in $\T$.
\item $\X$ is rigid and  $\U$ is $\X$-Frobenius.
\end{enumerate}
 In this case the subfactor category $\U/\X$ is triangulated.  
\item For any functorially finite rigid subcategory $\X$ of $\T$,  the maps 
\[
(\U,\U) \,\, \longmapsto \,\, \U/\X \,\,\,\,\,\,\,\,\, \text{and} \,\,\,\,\,\,\,\, \mathcal S \longmapsto   \big(\pi^{-1}(\mathcal S),\pi^{-1}(\mathcal S)\big)
\]
give mutually inverse bijections between:
\begin{enumerate}
\item[$\mathsf{(I)}$] $\X$-mutation pairs $(\U,\U)$ in $\T$.
\item[$\mathsf{(II)}$] Thick triangulated subcategories $\mathcal S$ of $\T/\X$ such that $\pi^{-1}(\mathcal S)$ is closed under extensions in $\T$.   
\end{enumerate} 
\item If $\T$ is $2$-Calabi-Yau and $\X$ is a rigid subcategory of $\T$, then for any  $\X$-Frobenius subcategory $\U$ of $\T$, the subfactor category $\U/\X$ is $2$-Calabi-Yau. 
\end{enumerate}
\begin{proof} Parts (i) and  (ii) follow from Theorem $3.3$ and the above discussion.

(iii) Fix bifunctorial isomorphisms $\omega_{A,B} \colon \T(B[-1],A) \stackrel{\cong}{\lxr} \mathsf{D}\T(A,B[1])$, $\forall A,B\in \T$, and define maps:
\[
\begin{CD}
\widetilde{\omega}_{U,V} \,\, \colon \,\, \Hom_{\U/\X}\big(\Omega^{1}_{\X}(V), U\big)  \,\, @>  >> \,\,   \mathsf{D}\Hom_{\U/\X}\big(U,\Sigma^{1}_{\X}(V) \big)
\end{CD}
\]
$\forall U,V \in \U$, as the composition of the following maps:
\[
\small{\begin{CD}
\Hom_{\U/\X}\big(\Omega^{1}_{\X}(V), U\big)  \,\, @> \mu_{U,V} >> \,\, \T(V[-1],U) \,\, @> \omega_{U,V} >> \,\, 
\mathsf{D}\T(U,V[1])  \,\, @> \mathsf{D}\nu_{U,V} >> \,\,   \mathsf{D}\Hom_{\U/\X}(U,\Sigma^{1}_{\X}(V)) 
\end{CD}}
\]
where $\mu_{U,V}(\underline{\alpha}) = h^{0}_{V}[-1] \circ \alpha$, $\forall \underline{\alpha} \in \Hom_{\U/\X}\big(\Omega^{1}_{\X}(V), U\big)$, and $\nu_{U,V}(\underline{\beta}) = \beta\circ h^{V}_{0}[1]$, $\forall \underline{\beta} \in \Hom_{\U/\X}\big(U,\Sigma^{1}_{\X}(V)\big)$. 
Note that since $\U$ is $\X$-Frobebius, it follows that $\mu$ and $\nu$ are well-defined natural maps which are monomorphisms. Using that $\X$ is rigid and condition (iii)(a) in Definition $3.2$, it is easy to see that $\mu$ and $\nu$, hence also $\mathsf{D}\nu$, are invertible. Hence $ \widetilde{\omega}_{U,V}$ is a natural isomorphism and therefore  $\U/\X$ is $2$-Calabi-Yau. 
\end{proof}
\end{cor}

\subsection{Rigidity and Extension closure} We show that if the triangulated category $\T$ admits a Serre functor $\mathbb S$, then  the functorially finite rigid subcategories $\X$ of $\T$ such that $\mathbb S(\X) = \X[2]$ are in bijective correspondence with the, maximal in a certain sense, functorially finite extension-closed subcategories of $\T$ such that $\mathbb S(\U) = \U[2]$. 

To this end we need the following result which generalizes, and is inspired by, results of \cite{BIRS} and which in particular shows that for any functorially finite rigid subcategory  $\X$ of $\T$, there is a maximal thick triangulated subcategory of $\T/\X$ whose preimage under $\pi \colon \T \lxr \T/\X$ is closed under extensions.

\begin{prop} Let $\T$ be a triangulated category  with a Serre functor $\mathbb S$.
\begin{enumerate} 
\item Let $\U$ be  a  functorially finite extension-closed subcategory of $\T$ such that $\mathbb S(\U) = \U[2]$.  Then $\X := \U \cap {^{\bot}}\U[1]$ is a functorially finite rigid subcategory of $\T$, $\mathbb S(\X) = \X[2]$ and $\U$ is $\X$-Frobenius. 
\item Let $\X$ be a functorially finite rigid subcategory of $\T$ such that $\mathbb S(\X) = \X[2]$.  Then $\U := {^{\bot}}\X[1]$ is $\X$-Frobenius, $\mathbb S(\U) = \U[2]$,  and the subfactor $\U/\X$ is a maximal thick triangulated subcategory of $\T/\X$. 
\end{enumerate} 
\begin{proof} (i) The condition $\mathbb S(\U) = \U[2]$ implies easily that $\U^{\bot} = {^{\bot}}\U[2]$ or equivalently $\U^{\bot}[-1] = {^{\bot}}\U[1]$. Then $\T(\X,\X[1]) = \T(\U\cap {^{\bot}}\U[1], \U[1]\cap {^{\bot}}\U[2]) = \T(\U\cap {^{\bot}}\U[1], \U[1]\cap \U^{\bot}) = 0$, hence $\X$ is rigid.
Since $\U[1]$ is covariantly finite in $\T$, for any object $U \in \U$, there exists a triangle $X_{U} \lxr U \lxr V[1] \lxr X_{U}[1]$, where the middle map is a minimal left $\U[1]$-approximation of $\U$. By Wakamatsu's Lemma, see \cite{BR}, \cite{IY}, the minimality condition implies that $X_{U} \in {^{\bot}}\U[1]$ and then $X_{U} \in \X$ since $\U$ is closed under extensions. Clearly then the map $X_{U} \lxr U$ is a right $\X$-approximation of $U$ and by construction $\Omega^{1}_{\X}(U) = V \in \U$. Hence $\X$ is contravariantly finite in $\U$ and $\Omega^{1}_{\X}(\U) \subseteq \U$.  By duality $\X$ is covariantly finite in $\U$ and $\Sigma^{1}_{\X}(\U) \subseteq \U$. Hence $\U$ is $\X$-Frobenius.  Finally we show that $\mathbb S(\X) = \X[2]$. Since $\mathbb S(\X) = \mathbb S(\U \cap {^{\bot}}\U[1]) = \mathbb S(\U) \cap \mathbb S({^{\bot}}\U[1]) = \U[2] \cap \mathbb S({^{\bot}}\U[1])$, it suffices to show that $\mathbb S({^{\bot}}\U[1]) = {^{\bot}}\U[3]$ or equivalently $\mathbb S({^{\bot}}\U) = {^{\bot}}\U[2]$. Since ${^{\bot}}\U[2] = \U^{\bot}$, this follows directly by Serre duality.

(ii) As in (i), the condition $\mathbb S(\X) = \X[2]$ implies that $\X^{\bot} = {^{\bot}}\X[2]$ and clearly then  $\U := \X^{\bot}[-1] = {^{\bot}}\X[1]$ is closed under extensions and direct summands in $\T$. For any object $A \in \T$, we have a triangle $\Omega^{1}_{\X}(A) \lxr X^{0}_{A} \lxr A \lxr \Omega^{1}_{\X}(A)[1]$ and $\Sigma^{1}_{\X}(A)[-1] \lxr A \lxr X^{A}_{0} \lxr \Sigma^{1}_{\X}(A)$, see (3.1). Using that $\X$ is rigid,  Lemma $2.1$ ensures that  
$\Omega^{1}_{\X}(A)[1]$ lies in $\X^{\bot}$ and $\Sigma^{1}_{\X}(A)[-1]$ lies in ${^{\bot}}\X$, hence we have:  $\Omega^{1}_{\X}(A)$ lies in $\X^{\bot}[-1] = \U$ and $\Sigma^{1}_{\X}(A)$ lies in ${^{\bot}}\X[1] = \U$. Finally we have $\T(\X,\U[1]) =
 \T(\X,\X^{\bot}[-1][1]) = \T(\X,\X^{\bot}) = 0$ and  $\T(\U,\X[1]) = \T({^{\bot}}\X[1],\X[1]) =\T({^{\bot}}\X,\X) = 0$. Therefore $\U$ is $\X$-Frobenius, hence, by Theorem $3.3$, $\U/\X$ is a thick triangulated subcategory of $\T/\X$.   If $\mathcal S$ be a thick triangulated subcategory of $\T/\X$ such that $\V := \pi^{-1}(\mathcal S)$ is closed under extensions, then by Corollary $3.12$ we have an $\X$-mutation pair  $(\V,\V)$ in $\T$. In particular $\T(\V,\X[1]) = \T(\V[-1],\X) = 0$, hence $\V[-1] \subseteq {^{\bot}}\X$ and then $\V \subseteq {^{\bot}}\X[1] = \U$. Therefore $\mathcal S  = \V/\X \subseteq \U/\X$. Finally as in (i) we have $\mathbb S(\U) = \U[2]$. 
\end{proof}
\end{prop}   

The above result suggests to introduce the following:

\begin{defn} Let $\T$ be a triangulated category with Serre duality and $\U$ a full subcategory of $\T$ closed under extensions.  The {\bf heart} of $\U$ is defined to be the full rigid subcategory $\U \cap {^{\bot}}\U[1] = \U \cap \U^{\bot}[-1]$. The subcategory $\U$ is called {\bf maximal} if it is maximal among those with the same heart: for any extension closed subcategory $\V$ of $\T$, the inclusion $\U\subseteq \V$ is an equality provided that $\U$ and $\V$ have the same heart. \end{defn}

The following result, which is Theorem E from the Introduction, generalizes some results of \cite{BIRS} concerning cluster structures on $2$-Calabi-Yau categories and  gives an interesting connection between rigid subcategories and extension closed subcategories (satisfying some additional conditions).  It should be noted that part (ii) can be regarded as a triangulated analog of the bijective correspondence between Wakamatsu-tilting modules and special coresolving subcategories over an Artin algebra observed by Mantese-Reiten in \cite[Theorems $2.10$, $3.4$]{MR}.   

\begin{thm} Let $\T$ be a triangulated category  with a Serre functor $\mathbb S$. 
\begin{enumerate}
\item The assignments   
\[
 \X \ \ \longmapsto \ \ \Psi(\X) := {^{\bot}}\X[1] \ \ \ \ \ \ \ \text{and} \ \ \ \ \ \ \ \ \ 
  \U \ \ \longmapsto \ \ \Phi(\U) := \U \cap {^{\bot}}\U[1]
\]
induce maps $\Psi \colon \mathsf{Rigid}(\T) \lxr \mathsf{Ext}(\T)$ and $\Phi \colon \mathsf{Ext}(\T) \lxr \mathsf{Rigid}(\T)$, where:
\begin{enumerate}
\item $\mathsf{Rigid}(\T)$  is the class of functorially finite rigid subcategories $\X$ of $\T$ such that $\mathbb S(\X) = \X[2]$.
\item $\mathsf{Ext}(\T)$ is the class of functorially finite extension closed subcategories $\U$ of $\T$ such that $\mathbb S(\U) = \U[2]$. 
\end{enumerate}
\item We have  
\[
\Phi\Psi(\X) = \X, \ \ \forall \X \in \mathsf{Rigid}(\T) \ \ \ \ \   \text{and}  \ \ \ \ \ \U \subseteq \Psi\Phi(\U),  \ \ \forall \U \in \mathsf{Ext}(\T)
\] Moreover  the maps $\Phi$, $\Psi$ induce  mutually inverse bijections 
\[
\Psi \,  \colon \,   \mathsf{Rigid}(\T) \  \lxr \  \mathsf{MaxExt}(\T) \ \ \ \ \ \ \ \text{and} \ \ \ \ \ \ \ \ \ 
\Phi \, \colon \,  \mathsf{MaxExt}(\T) \  \lxr \  \mathsf{Rigid}(\T) 
\]
between $\mathsf{Rigid}(\T)$ and the subclass $\mathsf{MaxExt}(\T)$  of $\mathsf{Ext}(\T)$ consisting of all maximal subcategories. 
\item For any $\U \in  \mathsf{Ext}(\T)$, the subfactor category $\U/\Phi(\U)$ is a triangulated subcategory of $\T/\Phi(\U)$, and for any $\X \in \mathsf{Rigid}(\T)$, the subfactor category $\Psi(\X)/\X$ is a maximal triangulated subcategory of $\T/\X$. 
\end{enumerate}  
\begin{proof} (i) If $\U$ lies in $\mathsf{Ext}(\T)$, then by Proposition $3.13$(i) it follows that $\Phi(\U) = \U \cap {^{\bot}}\U[1]$ is a rigid subcategory of $\T$ which is functorially finite in $\U$. Since $\U$ is functorially finite in $\T$, it follows directly that $\Phi(\U)$ is functorially finite in $\T$.  
Hence the assignment $\U \longmapsto \Phi(\U) = \U \cap {^{\bot}}\U[1]$ gives us a map $\Phi \colon \mathsf{Ext}(\T) \lxr \mathsf{Rigid}(\T)$.  

On the other hand if $\X$ lies in $\mathsf{Rigid}(\T)$, then clearly $\Psi(\X) = {^{\bot}}\X[1] = \X^{\bot}[-1]$ is an extension closed subcategory of $\T$. Since $\X$ is rigid, by Lemma $2.2$, the subcategory $\X^{\bot}$ is covariantly finite in $\T$ and the subcategory ${^{\bot}}\X$ is contravariantly finite in $\T$.  This implies that   $\X^{\bot}[-1]$ is covariantly finite in $\T$ and the subcategory ${^{\bot}}\X[1]$ is contravariantly finite in $\T$. Since $\X^{\bot}[-1] = {^{\bot}}\X[1]$, we infer that $\Psi(\X)$ is a functorially finite  extension closed subcategory of $\T$ and by Proposition $3.13$(ii) we have $\mathbb S\Psi(\X) = \Psi(\X)[2]$.  
Hence the assignment $\X \longmapsto \Psi(\X) = {^{\bot}}\X[1]$ gives us a map $\Psi \colon \mathsf{Rigid}(\T) \lxr \mathsf{Ext}(\T)$. 

(ii) We show that $\Phi\Psi(\X) = \X$, $\forall \X \in \mathsf{Rigid}(\T)$, or equivalently $\X = {^{\bot}}\X[1] \cap {^{\bot}}({^{\bot}}\X[1])[1]$. First since $\X$ is rigid, we have $\X \subseteq {^{\bot}}\X[1]$. On the other hand we have 
\[
\T(\X[-1], {^{\bot}}\X[1]) = \T(\X[-1],\X^{\bot}[-1]) = \T(\X,\X^{\bot}) = 0
\]
 hence $\X[-1] \subseteq  {^{\bot}}({^{\bot}}\X[1])$ or equivalently $\X \subseteq {^{\bot}}({^{\bot}}\X[1])[1]$. We infer that $\X \subseteq {^{\bot}}\X[1] \cap {^{\bot}}({^{\bot}}\X[1])[1]$. Conversely let $A\in {^{\bot}}\X[1] \cap {^{\bot}}({^{\bot}}\X[1])[1]$. Then $A[-1] \in {^{\bot}}({^{\bot}}\X[1]) = {^{\bot}}(\X^{\bot}[-1])$ and therefore $\T(A[-1],\X^{\bot}[-1]) = 0$ or equivalently $\T(A,\X^{\bot}) = 0$. This implies that the left $\X^{\bot}$-approximation $h^{0}_{A} \colon A \lxr \Omega^{1}_{\X}(A)[1]$ of $A$ is zero and therefore $A$ lies in $\X$ as a direct summand of $X^{0}_{A}$. We infer that $\X = {^{\bot}}\X[1] \cap {^{\bot}}({^{\bot}}\X[1])[1]$.

 Let $\U \in \mathsf{Ext}(\T)$; then $\Psi\Phi(\U) = {^{\bot}}(\U\cap {^{\bot}}\U[1])[1]$. For any $U \in \U$,  using that $\U^{\bot} = {^{\bot}}\U[2]$, we have 
 \[
 \T(U[-1], \U\cap {^{\bot}}\U[1]) = \T(U,\U[1]\cap  {^{\bot}}\U[2]) = \T(U,\U[1]\cap \U^{\bot}) = 0
 \]
  Hence $U[-1] \in {^{\bot}}( \U\cap {^{\bot}}\U[1])$ and therefore $U \in  {^{\bot}}(\U\cap {^{\bot}}\U[1])[1] = \Psi\Phi(\U)$. We infer that $\U \subseteq \Psi\Phi(\U)$. On the other hand,   
 by (i) the extension closed functorially finite subcategories $\U$ and $\Psi\Phi(\U)$ of $\T$ have the same heart $\Phi(\U)$. Hence if $\U$ lies in $\mathsf{MaxExt}(\T)$, we have $\U = \Psi\Phi(\U)$. Conversely it follows by Proposition $3.13$(ii) that $\Psi(\X)$ lies in $\mathsf{MaxExt}(\T)$, for any $\X \in \mathsf{Rigid}(\T)$. We infer that the map $\Psi \colon \mathsf{Rigid}(\T) \lxr \mathsf{MaxExt}(\T)$ is a bijection with inverse the map $\Phi \colon \mathsf{MaxExt}(\T) \lxr \mathsf{Rigid}(\T)$.   
 
 (iii) Follows from (i), (ii) and Proposition $3.13$. 
\end{proof}
\end{thm} 

\begin{rem} Using the correspondence described in part (i) of Theorem $3.15$, it is easy to see that the intersection  $\mathsf{Rigid}(\T) \bigcap \mathsf{MaxExt}(\T)$ is the class of all functorially finite (extension closed) subcategories $\X$ of $\T$ such that $\X^{\bot} = \X[1]$. In other words,  $\mathsf{Rigid}(\T) \bigcap \mathsf{MaxExt}(\T)$ is the class of $2$-cluster tilting subcategories of $\T$.  
\end{rem} 

\begin{cor} If $\T$ is a $2$-Calabi-Yau triangulated category, then there is a bijective correspondence between: 
\begin{enumerate}
\item  Functorially finite rigid subcategories of $\T$.
\item Maximal functorially finite extension closed subcategories of $\T$. 
\end{enumerate}
In particular this gives a bijection between basic rigid objects of $\T$ and maximal functorially finite extension closed subcategories of $\T$ whose heart has an additive generator. 
\end{cor}

\subsection{Algebraic Subfactors} Finally we observe that if the extension closed subcategory $\U$ satisfies an additional vanishing condition, namely that $\U$ is corigid, that is $\T(\U,\U[-1]) = 0$, then the subfactor triangulated category $\U/\X$ of Corollary $3.12$ is algebraic. Recall that a triangulated category is called {\em algebraic} if it is triangle equivalent to the stable category of an exact Frobenius category, see \cite{Happel} for more details.  

\begin{cor}  Let $\U$ be an $\X$-Frobenius subcategory of $\T$, where $\X$ is rigid, i.e. $(\U,\U)$ is an $\X$-mutation pair in $\T$. If $\T(\U,\U[-1]) = 0$, then $\U$ is an exact Frobenius category with $\X$ as the full subcategory of projective-injective objects and the triangulated subfactor category $\U/\X$ is algebraic. 
\begin{proof} By a result of Dyer \cite{Dyer} the class of diagrams $0 \lxr A \lxr B \lxr C \lxr 0$ in $\U$, such that there exists a triangle  $A \lxr B \lxr C \lxr A[1]$ in $\T$, is the class of admissible short exact sequences for an exact structure (in the sense of Quillen) in $\U$. Alternatively it is easy to see that in our setting the functor $\mathsf{H} \colon \U \lxr \smod\X$, $\mathsf{H}(U) = (-,U)|_{\X}$ is a full embedding and $\Image\mathsf{H}$ is closed under extensions in $\smod\X$, so $\U$ is an exact category with exact structure as indicated above. If $X \in \X$, then since $\T(\X,\U[1]) = 0$, any map $X \lxr C$ factorizes through $B$, and since $\T(\U[-1],\X) = 0$, any map $A \lxr X$ factorizes through $B$. It follows that $\X$ consists of projective-injective objects of $\U$. Since $\X$ is functorially finite in $\U$, condition (ii) in Definition $3.2$ implies that $\U$ has enough projective and injective objects. Hence $\U$ is a Frobenius category, so   $\U/\X$ is algebraic.  
\end{proof}  
\end{cor}

\section{Abelian Localizations and Categories of Fractions} 

In this section, in contrast to section 3,  we are interested in abelian categories arising, via Gabriel-Zisman localization, from certain factor categories of triangulated categories equipped with a contra(co)variantly  finite rigid subcategory. Our main result is Theorem $4.6$ which is Theorem A from the Introduction.

\subsection{Rigidity and categories of extensions} Let $\T$ be a triangulated category and $\X$ a full subcategory of $\T$ closed under direct summands and isomorphisms. In this subsection we study consequences of rigidity on $\X$ to certain factor categories of $\T$ associated to $\X$ and also to the category of coherent functors over $\X$. 
As in section 2, we denote by $\mathsf{H} \colon \T \lxr \smod\X$ the Yoneda restriction homological functor, $\mathsf{H}(A) = \T(-,A)|_{\X}$.

\begin{lem} Let $\X$ be a contravariantly finite rigid subcategory of $\T$ and consider a triangle 
\[
\begin{CD}
A \,\ @> \alpha >> \,\, B \,\ @> \beta >> \,\, C \,\ @> \gamma >> \,\, A[1]
\end{CD}
\] 
\begin{enumerate}
\item $\underline{\beta} = 0$ in $\T/\X^{\bot}$ if and only if $\mathsf{H}(\beta) = 0$, i.e. $\beta$ is $\X$-ghost, if and only if $\beta$ factorizes through $\Omega^{1}_{\X}(B)[1]$.  
\item $\underline{\beta}$ is an epimorphism in $\T/\X^{\bot}$ if and only if $\mathsf{H}(\beta)$ is an epimorphism in $\smod\X$ if and only if $\underline{\gamma} = 0$ in $\T/\X^{\bot}$ if and only if $\gamma$ factorizes through $\Omega^{1}_{\X}(C)[1]$.
\item $\underline{\beta}$ is a monomorphism in $\T/\X^{\bot}$ if and only if $\mathsf{H}(\beta)$ is a monomorphism in $\smod\X$ if and only if $\underline{\alpha} = 0$ in $\T/\X^{\bot}$ if and only if $\alpha$ factorizes through $\Omega^{1}_{\X}(A)[1]$.
\end{enumerate}
\begin{proof} (i) Since $\Omega^{1}_{\X}(B)[1]$ lies in $\X^{\bot}$, we have $\underline{\beta} = 0$ in $\T/\X^{\bot}$ if and only if $\beta$ factorizes through an object from $\X^{\bot}$ if and only if $\beta$ factorizes through the left $\X^{\bot}$-approximation $\Omega^{1}_{\X}(B)[1]$ of $B$. If this holds, then clearly $\mathsf{H}(\beta) = 0$ since $\X^{\bot} = \Ker\mathsf{H}$.  Conversely if $\mathsf{H}(\beta) = 0$, then the composition $f^{0}_{B} \circ \beta$ is zero, so $\beta$ factorizes through the cone $\Omega^{1}_{\X}(B)[1]$ of $f^{0}_{B}$. Since $\mathsf{H}$ is homological, parts (ii), (iii) follow directly from (i).
\end{proof} 
\end{lem}

An object $R$ in an additive category $\A$ is {\em projective} if for any epic $B \lxr C$ in $\A$, the induced map $\A(R,B) \lxr \A(R,C)$ is surjective. $\A$ has {\em enough projective objects}, if for any object $A$ there is an epic $P \lxr A$, where $P$ is projective. Dually one defines {\em injective objects} and when $\A$ has {\em enough injective objects}.

\begin{lem} Let $\X$ be a contravariantly finite rigid subcategory of $\T$. Then the stable category $\T/\X^{\bot}$ has cokernels, enough projectives and    $\Proj(\T/\X^{\bot}) = \X/\X^{\bot} \approx \X$.
 \begin{proof}  Let $\underline{\alpha} \colon \unA \lxr \unB$ be a map in $\T/\X^{\bot}$ and consider the homotopy push-out diagram in $\T$  
 \begin{equation}
  \xymatrix@C=1.5cm{
    X^{0}_{A}  \ar[r]^{f^{0}_{A}}  \ar@{=}[d]  &  
   A \ar[r]^{h^{0}_{A} \ \ \ }\ar[d]^{\alpha}&
    \Omega^{1}_{\X}(A)[1] \ar[d]^{\delta} \ar[r]^{\ \ -g^{0}_{A} \ \ }&
    X^{0}_{A}[1] \ar@{=}[d] \\
    X^{0}_{A} \ar[r]^{e}& B \ar[r]^{c}& C \ar[r]^{d} & X^{0}_{A}[1]}
\end{equation}  
of the rotated triangle $(2.1)$ along the map $\alpha$. Since $\Omega^{1}_{\X}(A)[1]$ lies in $\X^{\bot}$, the composition $\underline{\alpha} \circ \underline{c}$ is zero in $\T/\X^{\bot}$. Since $X^{0}_{A}[1]$ lies in $\X^{\bot}$, $\underline{c}$ is an epimorphism by Lemma $4.1$.  Now let $\underline{g} \colon \underline{B} \lxr \underline{D}$ be a map in $\T/\X^{\bot}$ such that $\underline{\alpha} \circ \underline{g} = 0$. Then the composition $\alpha \circ g$ factorizes through an object from $\X^{\bot}$, hence  it factorizes through the left $\X^{\bot}$-approximation $h^{0}_{A} \colon A \lxr \Omega^{1}_{\X}(A)[1]$ of $A$, that is $\alpha \circ g = h^{0}_{A} \circ \rho$, for some map $\rho \colon \Omega^{1}_{\X}(A)[1] \lxr D$. Since $(4.1)$ is a homotopy push-out diagram, there is a map $\sigma \colon C \lxr D$ such that $c\circ \sigma = g$, hence $\underline{c} \circ \underline{\sigma} = \underline{g}$.  The map $\sigma$ is unique in $\T/\X^{\bot}$ since $\underline{c}$ is an epic. 
We infer that $\underline{c}$ is the cokernel of $\underline{\alpha}$ in $\T/\X^{\bot}$ and therefore $\T/\X^{\bot}$ has cokernels. Now assume that $\underline{\alpha}$ is epic and let $\underline{\rho} \colon \unX \lxr \unB$ be a map in $\T/\X^{\bot}$, where $X \in \X$.  Since $\underline{\alpha}$ is an epic, it follows that its cokernel $\unC$ is zero, i.e. $C$ lies in $\X^{\bot}$, hence the composition $\rho \circ c \colon X \lxr C$ is zero. Therefore $\rho$ factorizes through $e = f^{0}_{A} \circ \alpha$ and then $\underline{\rho}$ factorizes through $\underline{\alpha}$. We infer that $\unX$ is projective in $\T/\X^{\bot}$. Finally for any object $A \in \T$, the diagram $(4.1)$ shows that the cokernel of $\unf^{0}_{A} \colon \unX^{0}_{A} \lxr \unA$ in $\T/\X^{\bot}$ is the map $\underline{h}^{0}_{A} \colon \unA \lxr \Omega^{1}_{\X}(A)[1]$. It follows that $\unf^{0}_{A} $ is an epimorphism in $\T/\X^{\bot}$ since  $\Omega^{1}_{\X}(A)[1]\in \X^{\bot}$.  We conclude that $\T/\X^{\bot}$ has enough projective objects. If $\unA$ is projective, then $\unf^{0}_{A}$ splits in $\T/\X^{\bot}$ and therefore $\unA$ lies in $\X/\X^{\bot}$. Hence $\Proj(\T/\X^{\bot}) = \X/\X^{\bot}$.  Clearly the projection functor $\X \lxr \X/\X^{\bot}$, $X \longmapsto \underline{X}$ is an equivalence.     
\end{proof}
\end{lem}

Next we show that the extension category $\X\star\X[1]$ is contravariantly finite in $\T$ and the associated factor category $(\X\star \X[1])/\X^{\bot}$ is coreflective in $\T/\X^{\bot}$. 
In this respect we need the following construction. Let $\U$ and $\V$ be full subcategories of $\T$. Recall that the {\em extension category} $\U\star\V$ consists of all direct summands of objects $A$ for which there exists a triangle $U\lxr A \lxr V \lxr U[1]$, where $U$ lies in $\U$ and $V$ lies in $\V$.  
Now let $A$ be an object of $\T$ and form the homotopy push-out of the triangle 
$(2.1)$ along the map $h^{1}_{A} \colon \Omega^{1}_{\X}(A) \lxr \Omega^{2}_{\X}(A)[1]$:

\begin{equation}
\xymatrix@C=1.5cm{
                                 X^{1}_{A} \ar@{=}[r] \ar[d]  & X^{1}_{A}  \ar[d] &                \\
\Omega^{1}_{\mathcal X}(A) \ar[d]_{h^{1}_{A}} \ar[r]^{g^{0}_{A}} & X^{0}_{A} \ar[r]^{f^{0}_{A}} \ar[d]^{\alpha_{A}}           & A \ar@{=}[d]   \ar[r]^{h^{0}_{A} \ \ \ \ \ }   & \Omega^{1}_{\mathcal X}(A)[1] \ar[d]^{h^{1}_{A}[1]} \\
\Omega^{2}_{\mathcal X}(A)[1]   \ \ar[r]^{\beta_{A}} \ar[d]^{}            &   \ \mathsf{Cell}_{1}(A) \ar[r]^{\ \ \gamma_{A}} \ar[d]          & A  \ar[r]^{\omega_{A}\ \ \ \ \ }         & \Omega^{2}_{\mathcal X}(A)[2]              \\
                                 X^{1}_{A}[1] \ar@{=}[r]                      & X^{1}_{A}[1] &
}
\end{equation}

\, 
 
\begin{rem} \begin{enumerate}
\item[(i)] It follows from $(4.2)$ that the object $\mathsf{Cell}_{1}(A)$ lies in $\X \star \X[1]$. The morphism $\gamma_{A} \colon \Cell_{1}(A) \lxr A$ is monic and epic in $\T/\X^{\bot}$, and $\mathsf{H}(\gamma_{A})$ is invertible in $\smod\X$.  Indeed this follows from Lemma $4.1$ using that $\Omega^{2}_{\X}(A)[1]$ lies in $\X^{\bot}$ and  the map $\omega_{A}$ factorizes through $\Omega^{1}_{\X}(A)[1]$. 

\item[(ii)] The map $\gamma_{A} \colon \mathsf{Cell}_{1}(A) \lxr A$, which by Lemma $4.4$ below is a right $(\X\star\X[1])$-approximation of $A$, is called the {\em first cellular approximation of $A$ with respect to $\X$}.  For a justification of the terminology and the notation, which are of topological origin, we refer to \cite{B}. 

\item[(iii)] Note that we have an equality $(\X\star \X[1])/\X^{\bot} = (\X\star \X[1])/\X[1]$.  Indeed the inclusion $(\X\star \X[1])/\X[1] \subseteq (\X\star \X[1])/\X^{\bot}$ is clear since $\X$ is rigid, so $\X[1] \subseteq \X^{\bot}$. The other inclusion follows from the fact that a map $\alpha \colon A \lxr B$, where $A \in \X\star \X[1]$ factors through an object of $\X^{\bot}$ if and only if $\alpha$ factors through $X^{1}[1] \in \X[1]$, where $X^{1} \lxr X^{0} \lxr A \lxr X^{1}[1]$ is a triangle in $\T$ and the $X^{i} \in \X$. 
\end{enumerate}
\end{rem} 

\,

Consider now the inclusions $\X[-1]\star \X \, \subseteq \, \T \, \supseteq \X \star \X[1]$ which induce fully faithful functors
\[
\begin{CD}
\mathsf{L} \, \colon \,  (\X\star \X[1])/\X^{\bot} \,\ @>  >> \,\  \T/\X^{\bot} \,\,\,\ \ \  \text{and} \,\,\,\,\,\ \ \  \mathsf{R}^{\op} \, \colon \,   (\X[-1]\star \X)/{^{\bot}}\X \,\  @>  >> \,\  \T/{^{\bot}}\X
\end{CD}
\]

\,

\begin{lem} Let $\X$ be a contravariantly finite rigid subcategory of $\T$. 
\begin{enumerate}
\item $\X\star\X[1]$ is contravariantly finite in $\T$ and the factor category $(\X\star\X[1])/\X^{\bot}$ is abelian.  More precisely the homological functor $\mathsf{H} \colon \T \lxr \smod\X$ induces an equivalence
\[
\begin{CD}
\underline{\mathsf{H}} \,\ \colon \,\, (\X\star\X[1])/\X^{\bot} \,\,\, @> \approx >>  \,\,\, \smod\X
\end{CD}
\] 
\item  $(\X\star\X[1])/\X^{\bot}$ is coreflective in $\T/\X^{\bot}$; there is an adjoint pair $(\mathsf{L},\mathsf{R})$ where $\mathsf{L}$ is fully faithful:
\[
(\mathsf{L}, \mathsf{R}) \, \colon \, \xymatrix@C=3.5pc {(\X\star\X[1])\big/\X^{\bot} \ \ \ar@<0.5ex>[r]^-{} & \ar@<0.5ex>[l]^-{{}} \ \ \T/\X^{\bot}}
\]
\end{enumerate}
\begin{proof} Let $A$ be an object of $\T$ and consider the triangle constructed in the diagram (4.2):
\[
\begin{CD}
\Omega^{2}_{\X}(A)[1] \,\ @>  \beta_{A} >>  \,\, \mathsf{Cell}_{1}(A) \,\ @> \gamma_{A} >>  \,\, A \,\, @> \omega_{A} >> \,\ \Omega^{2}_{\X}(A)[2] 
\end{CD} \eqno (\dag)
\]

(i) Note that $\omega_{A} = h^{0}_{A} \circ h^{1}_{A}[1]$ factorizes through $\Omega^{1}_{\X}(A)[1] \in \X^{\bot}$  and we know that $\mathsf{Cell}_{1}(A)$ lies in $\X\star \X[1]$.  Let $\delta \colon M \lxr A$ be a map in $\T$, where $M \in \X\star\X[1]$.  We claim that the composition $\delta \circ \omega_{A} \colon M\lxr A \lxr  \Omega^{2}_{\X}(A)[2]$ is zero. Indeed let $X^{0} \stackrel{\kappa}{\lxr} M \stackrel{\lambda}{\lxr} X^{1}[1] \lxr X^{0}[1]$ be a triangle in $\T$, where the $X^{i}$ lie in $\X$.  Then the composition $X^{0} \lxr M \lxr A \lxr \Omega^{2}_{\X}(A)[2]$ is zero since the last map factorizes through $\Omega^{1}_{\X}(A)[1]$. Hence the map $M \lxr A \lxr \Omega^{2}_{\X}(A)[2]$ factorizes through $X^{1}[1]$, say via a map $X^{1}[1] \lxr \Omega^{2}_{\X}(A)[2]$. However this last map is zero since $\Omega^{2}_{\X}(A)[1] \in \X^{\bot}$. We infer that the composition  $M \lxr A \lxr \Omega^{2}_{\X}(A)[2]$ is zero and therefore the map $M \lxr A$ factorizes through $\mathsf{Cell}_{1}(A)$. This means that the map $\gamma_{A} \colon \mathsf{Cell}_{1}(A) \lxr A$ is a right $(\X\star\X[1])$-approximation of $A$. We infer that $\X\star \X[1]$ is contravariantly finite in $\T$.  The last assertion follows from Remark $4.3$ and a result of Keller-Reiten, see \cite[Lemma 5.1]{KR1}.

(ii) By (i)  the map $\gamma_{A} \colon \mathsf{Cell}_{1}(A) \lxr A$ is a right $(\X\star\X[1])$-approximation of $A$. This clearly implies that the map $\underline{\gamma}_{A} \colon  \underline{\mathsf{Cell}_{1}(A)} \lxr \underline{A}$ is a right $(\X\star\X[1])/\X^{\bot}$-approximation of $\unA$. We claim that $\underline{\gamma}_{A}$ is the coreflection of $\unA$ in $(\X\star\X[1])/\X^{\bot}$. To this end it suffices to show that if $\underline{\alpha}, \underline{\beta} \colon \underline{M} \lxr \underline{\mathsf{Cell}_{1}(A)}$ are maps, where $M \in \X\star \X[1]$ is as above,  such that $\underline{\alpha} \circ \underline{\gamma}_{A} =  \underline{\beta} \circ \underline{\gamma}_{A}$, then $\underline{\alpha} = \underline{\beta}$ in $\T/\X^{\bot}$. Indeed then the map $(\alpha - \beta)\circ \gamma_{A} \colon M \lxr A$ factorizes through an object from $\X^{\bot}$ and therefore it factorizes through the left $\X^{\bot}$-approximation $\lambda \colon M \lxr X^{1}[1]$ of $M$, say as $(\alpha - \beta)\circ \gamma_{A} = \lambda \circ \sigma$ for some map $\sigma \colon X^{1}[1] \lxr A$. Then $\kappa \circ (\alpha - \beta)\circ \gamma_{A} = \kappa \circ \lambda \circ \sigma = 0$, hence from $(\dag)$ there is a map $\tau \colon X^{0} \lxr \Omega^{2}_{\X}(A)[1]$ such that $\kappa \circ (\alpha - \beta) = \tau \circ \beta_{A}$. However $\tau = 0$ since $X^{0}\in \X$ and $\Omega^{2}_{\X}(A)[1] \in \X^{\bot}$ and therefore $\kappa \circ (\alpha - \beta) = 0$. It follows that $\alpha - \beta$ factorizes through the cone $X^{1}[1] \in \X^{\bot}$ of $\kappa$ and consequently $\underline{\alpha} = \underline{\beta}$ in $\T/\X^{\bot}$. We infer that any object $\unA$ in $\T/\X^{\bot}$ admits a coreflection in $(\X\star\X[1])/\X^{\bot}$. Then we obtain the right adjoint $\mathsf{R}$ of the inclusion $\mathsf{L}$ by  setting $\mathsf{R} \, \colon \, \T/\X^{\bot} \, \lxr \,(\X\star\X[1])/\X^{\bot}$, $\mathsf{R} (\unA) = \underline{\mathsf{Cell}_{1}(A)}$. 
\end{proof}
\end{lem}

We now state without proof the dual versions of Lemmas $4.2$ and $4.4$.
 
\begin{lem} Let $\X$ be a covariantly finite rigid subcategory of $\T$. 
\begin{enumerate}
\item The factor category $\T/{^{\bot}}\X$ has kernels, enough injectives and  $\Inj(\T/{^{\bot}}\X) = \X/{^{\bot}}\X \approx \X$.
\item $\X[-1]\star \X$ is covariantly finite in $\T$ and the factor category $(\X[-1]\star\X)/{^{\bot}}\X$ is abelian. More precisely the cohomological functor $\mathsf{H}^{\op} \colon \T \lxr \X\lsmod$ induces a duality
\[
\begin{CD}
\underline{\mathsf{H}}^{\op} \,\, \colon \,\, (\X[-1]\star\X)/{^{\bot}}\X \,\, @> \approx >> \,\, (\X\lsmod)^{\op}
\end{CD}
\]
\item $(\X[-1]\star\X)/{^{\bot}}\X$ is reflective in $\T/{^{\bot}}\X$; there is an adjoint pair $(\mathsf{L}^{\op},\mathsf{R}^{\op})$ where $\mathsf{R}^{\op}$ is fully faithful: 
\[
(\mathsf{R}^{\op}, \mathsf{L}^{\op}) \, \colon \, \xymatrix@C=3.5pc {(\X[-1]\star\X)\big/{^{\bot}}\X \ \ \ar@<0.5ex>[r]^-{{}} & \ar@<0.5ex>[l]^-{{}}\ \ \T/{^{\bot}}\X}
\]
\end{enumerate}
\end{lem}

\subsection{Categories of Fractions}
Recall that the Gabriel-Zisman localization of an additive category $\C$ at a class of morphisms $\Sigma$ of $\C$ is an additive category $\C[\Sigma^{-1}]$ equipped with an additive functor $\mathsf{P} \colon \C \lxr \C[\Sigma^{-1}]$ which makes the morphisms in $\Sigma$ invertible, i.e. $\mathsf{P}(s)$ is an isomorphism, for any $s \in \Sigma$, and is universal with this property. That is, for any functor $F \colon \C \lxr \M$ making the morphisms in $\Sigma$ invertible in $\M$, there exists a unique, up to a natural isomorphism,  functor $F^{*} \colon \C[\Sigma^{-1}] \lxr \M$ such that $F^{*}\circ \mathsf{P} = F$.

Also recall that a class of morphisms $\Sigma$ in $\C$ admits a {\em calculus of right fractions} if:
\begin{enumerate}
\item $\Sigma$ is closed under composition and contains the identities of $\C$,
\item If $f,g \colon A \lxr B$ are morphisms in $\C$ and $s \colon B \lxr D$ is a map in $\Sigma$ such that $f \circ s = g \circ s$, then there exists a map $t \colon C \lxr A$ in $\Sigma$ such that $t \circ f = t\circ g$. 
\item Any diagram $A \stackrel{f}{\lxr} C \stackrel{s}{\rxr} B$ in 
$\C$, where $s \in \Sigma$, can be completed to a commutative square with $t \in \Sigma$:
\[
\xymatrix@C=1.5cm{
  D \ar@{-->}[d]_{t \,\, } \ar@{-->}[r]^{g}  & B \ar[d]^{\, s}  \\
  A \ar[r]^{f} &  C  } 
\] 
\end{enumerate}

The notion of {\em calculus of left fractions} is defined dually;  we refer to Gabriel-Zisman's book \cite{GZ} for more details on localizations and calculus of fractions. 

Now we are ready to prove Theorem A from the Introduction. Recall that a map $f$ in a given category is called {\em regular} if $f$ is monic and epic.

\begin{thm} Let $\T$ be a triangulated category and $\X$ a contravariantly finite rigid subcategory of $\T$. Let $\mathcal R$ be the class of regular maps in $\T/\X^{\bot}$. Then the class $\mathcal R$ admits a calculus of  left and right fractions, the localization category $(\T/\X^{\bot})[\mathcal R^{-1}]$ exists, the canonical functor $\mathsf{P} \colon \T/\X^{\bot} \lxr (\T/\X^{\bot})[\mathcal R^{-1}]$ is faithful, preserves kernels and cokernels, and admits a fully faithful left adjoint, and  there is an equivalence
\[
\begin{CD}
 (\T/\X^{\bot})[\mathcal R^{-1}]\,\,\  @> \approx >> \,\,\ \smod\X
\end{CD}
\] 
\begin{proof} Consider the adjoint pair 
$\big(\mathsf{L}, \mathsf{R}\big)$ of Lemma $4.4$ and let $\mathcal R$ be the class of maps $\underline{\alpha} \colon \unA \lxr \unB$ in $\T/\X^{\bot}$ such that the map $\mathsf{R}(\unf)$ is invertible. By a classical result of Gabriel-Zisman \cite[Chapter 1]{GZ}, the localization $\mathsf{P} \colon \T/\X^{\bot} \lxr (\T/\X^{\bot})[\mathcal R^{-1}]$ of $\T/\X^{\bot}$ at  $\mathcal R$ exists, and the class $\mathcal R$ admits a calculus of right fractions and is saturated.  Recall that saturation means that $\underline{\alpha}$ lies in $\mathcal R$ if and only if $\mathsf{P}(\underline{\alpha})$ is invertible (in fact $\mathcal R$ is the saturation of the class $\{\varepsilon_{\unA} \colon \underline{\Cell_{1}(A)} \lxr \unA \,\, | \,\, \unA \in \T/\X^{\bot}\}$, where $\varepsilon$ is the counit of the adjoint pair $(\mathsf{L}, \mathsf{R}))$. Moreover there is an equivalence 
\[
\begin{CD}
\mathsf{R}^{*} \,\ \colon \,\,   (\T/\X^{\bot})[\mathcal R^{-1}]\,\,\  @> \approx >> \,\,\ (\X\star\X[1])/\X^{\bot}
\end{CD}
\] 
 such that $\mathsf{R}^{*} \circ \mathsf{P} = \mathsf{R}$, i.e. the following diagram commutes:
\begin{equation}
\begin{CD}
\T/\X^{\bot} @> \mathsf{P} >> \big(\T/\X^{\bot}\big)[\mathcal R^{-1}] \\
\Big\| & \ &  @V{\approx}V{\mathsf{R}^{*}}V  \\ 
\T/\X^{\bot} @> \mathsf{R} >> \big(\X\star\X[1]\big)/\X^{\bot}
\end{CD}
\end{equation} 
Composing $\mathsf{R}^{*}$ with the equivalence $\underline{\mathsf{H}} \, \colon \,  (\X\star\X[1])/\X^{\bot} \, \stackrel{\approx}{\lxr} \,\smod\X$
of Lemma $4.4$, we obtain an equivalence:
\[
\begin{CD}
\underline{\mathsf{H}} \circ \mathsf{R}^{*} \,\ \colon \,\,   (\T/\X^{\bot})[\mathcal R^{-1}]\,\,\  @> \approx >> \,\,\ \smod\X
\end{CD}
\] 
Clearly the functor $\mathsf{H} \colon \T \lxr \smod\X$ admits a factorization $ \mathsf{H} = \underline{\mathsf{H}} \circ \mathsf{R}^{*} \circ \mathsf{P} \circ \pi \ \colon \ \T \lxr \smod\X$:
\[
\begin{CD}
  \T   @> \pi >>  \T/\X^{\bot}   @> \mathsf{P} >>  (\T/\X^{\bot})[\mathcal R^{-1}]   @> \mathsf{R}^{*} >>  (\X\star\X[1])/\X^{\bot}   @> \underline{\mathsf{H}} >>  \smod\X 
\end{CD} \eqno (\dag\dag)
\] 
where $\pi$ is  the projection functor and the last two functors are equivalences. It follows directly from this that a map $\underline{\alpha}$ in $\T/\X^{\bot}$ lies in $\mathcal R$ if and only if $\mathsf{P}(\underline{\alpha})$ is invertible if and only if $\mathsf{H}(\alpha)$ is invertible. Moreover if $\mathsf{P}(\underline{\alpha}) = 0$, then $\underline{\mathsf{H}} \circ \mathsf{R}^{*} \circ \mathsf{P}(\underline{\alpha})
 = 0$ so  $\mathsf{H}(\alpha)  = 0$. Then $\underline{\alpha} = 0$ in $\T/\X^{\bot}$ by Lemma $4.1$ and $\mathsf{P}$ is faithful.   It remains to show that $\mathcal R$ admits a calculus of left fractions and consists of the regular maps in $\T/\X^{\bot}$. 

First let $\underline{\alpha}$ be a regular map in $\T/\X^{\bot}$, i.e. $\underline{\alpha}$ is monic and epic. Since the class $\R$ admits a calculus of right fractions, it follows by the dual of \cite[Proposition I.3.1, page 16]{GZ} that the functor $\mathsf{P}$ preserves kernels. Hence $\mathsf{P}(\underline{\alpha})$ is monic. On the other hand, by Lemma $4.1$, the map $\mathsf{H}(\alpha)$ is epic. Since $\mathsf{H}(\alpha) = \underline{\mathsf{H}} \mathsf{R}^{*}  \mathsf{P}(\underline{\alpha})$ and the functors $\underline{\mathsf{H}}$ and $\mathsf{R}^{*}$ are equivalences, this implies that $\mathsf{P}(\underline{\alpha})$ is epic. We infer that  $\mathsf{P}(\underline{\alpha})$ is regular in $(\T/\X^{\bot})[\mathcal R^{-1}]$ and therefore $\mathsf{P}(\underline{\alpha})$ is invertible since  $(\T/\X^{\bot})[\mathcal R^{-1}]$ is abelian. However since $\mathcal R$ is saturated, this means that $\underline{\alpha}$ lies in $\mathcal R$. Conversely if $\underline{\alpha}$ lies in $\mathcal R$, then $\mathsf{R}(\underline{\alpha})$ is invertible. Then so is $\mathsf{H}(\alpha) = \underline{\mathsf{H}}\mathsf{R}(\underline{\alpha})$. Then by Lemma $4.1$, we infer that $\underline{\alpha}$ is regular. We conclude that $\mathcal R$ consists of the regular maps in $\T/\X^{\bot}$. 

Next we show that $\mathcal R$ admits a calculus of left fractions. Since $\mathcal R$ is saturated and $\T/\X^{\bot}$ has cokernels, one can check easily that the conditions in \cite[19.3.5(b)]{Schubert} (ensuring that a given class of maps in a category with cokernels admits a calculus of left fractions) are satisfied by the class $\mathcal R$. Alternatively, by \cite[Proposition I.3.4]{GZ} we have to check that the functor $\mathsf{R} \colon \T/\X^{\bot} \lxr  (\X\star \X[1])/\X^{\bot}$ preserves cokernels. Since $\underline{\mathsf{H}}$ is an equivalence, it is equivalent to check that the functor $\underline{\mathsf{H}} \circ \mathsf{R} \colon \T/\X^{\bot} \lxr \smod\X$, $(\underline{\mathsf{H}} \circ \mathsf{R})(\unA) = \mathsf{H}(\Cell_{1}(A))$ preserves cokernels. However  if $\unA \lxr \unB \lxr \unC \lxr 0$ is exact in $\T/\X^{\bot}$, then by the construction of cokernels in $\T/\X^{\bot}$, see diagram (4.1) in Lemma $4.2$, it follows directly that $(*): \mathsf{H}(A) \lxr \mathsf{H}(B) \lxr \mathsf{H}(C) \lxr 0$ is exact in $\smod\X$. Since by Remark $4.3$, the map $\gamma_{A} \colon \Cell_{1}(A) \lxr A$ is made invertible by $\mathsf{H}$, it follows that  the induced sequence $\mathsf{H}(\Cell_{1}(A)) \lxr \mathsf{H}(\Cell_{1}(B)) \lxr \mathsf{H}(\Cell_{1}(C)) \lxr 0$ is exact in $\smod\X$ since it is isomorphic to the exact sequence $(*)$. We infer that the functor $\underline{\mathsf{H}} \circ \mathsf{R}$, hence the functor $\mathsf{R}$, preserves cokernels and therefore  $\mathcal R$ admits a calculus of left fractions. Then  $\mathsf{P}$ preserves cokernels by the dual of \cite[I.3.1]{GZ}. 

Finally from $(4.3)$ we have $\mathsf{P} = (\mathsf{R}^{*})^{-1} \circ \mathsf{R}$ and  the following natural isomorphisms, $\forall A\in (\T/\X^{\bot})[\mathcal R^{-1}]$ and $\forall B \in \T/\X^{\bot}$, show that the functor $\mathsf{L} \circ \mathsf{R}^{*} \colon (\T/\X^{\bot})[\mathcal R^{-1}] \lxr \T/\X^{\bot}$  is  a left adjoint of $\mathsf{P}$:
\[
\begin{CD}
{\Hom}\big(\mathsf{L}\mathsf{R}^{*}(A),B) \,\, @> \cong >>  \,\, {\Hom}\big(\mathsf{R}^{*}(A),\mathsf{R}(B)\big) \,\, @> \cong >> \,\,   {\Hom}\big(A,(\mathsf{R}^{*})^{-1}\mathsf{R}(B)\big) \,\, @> \cong >>\end{CD}
\]
\[
\begin{CD}   @> \cong >> \,\, {\Hom}\big(A,\mathsf{P}(B)\big) 
\end{CD}
\qedhere
\] 
\end{proof}
\end{thm} 

\

As a consequence we have the following. 

\begin{cor} The following are equivalent:
\begin{enumerate}
\item The factor category $\T/\X^{\bot}$ is abelian.
\item $\forall A \in \T$: the epimorphism $\underline{f}^{0}_{A} \colon \unX^{0}_{A} \lxr \unA$ in $\T/\X^{\bot}$ is a cokernel. 
\item  $\forall A \in \T$: the regular map $\underline{\gamma}_{A} \colon \underline{\Cell}_{1}(A) \lxr \unA$ in $\T/\X^{\bot}$ is invertible.
\end{enumerate}
\begin{proof} If $\forall A \in \T$, $\underline{f}^{0}_{A}$ is a cokernel, then  from the construction of cokernels in Lemma $4.2$ it follows that $A \in\X\star \X[1]$, hence $\T/\X^{\bot} = (\X\star\X[1])/\X^{\bot}$ and $\T/\X^{\bot}$ is abelian by Lemma $4.4$. Finally  if the map $\underline{\gamma}_{A}$ is invertible, $\forall A \in \T$, then the functor $\mathsf{R} \colon \T/\X^{\bot} \lxr (\X\star\X[1])/\X^{\bot}$ is an equivalence, hence $\T/\X^{\bot}$ is abelian.  
\end{proof}
\end{cor} 

The following is the dual version of Theorem $4.6$.

\begin{thm} Let $\T$ be a triangulated category and $\X$ a covariantly finite rigid subcategory of $\T$. Let $\mathcal R_{\op}$ be the class of regular maps in $\T/{^{\bot}}\X$. Then the class $\mathcal R_{\op}$ admits a calculus of  left and right fractions, the localization category $(\T/{^{\bot}}\X)[\mathcal R^{-1}_{\op}]$ exists, and the canonical functor $\mathsf{P} \colon \T/{^{\bot}}\X \lxr (\T/{^{\bot}}\X)[\mathcal R^{-1}_{\op}]$ is faithful, preserves kernels and cokernels and admits a fully faithful left adjoint. Moreover  there are equivalences: 
\[
\begin{CD} 
(\T/{^{\bot}}\X)[\mathcal R^{-1}_{\op}]\,\,\  @> \approx >> \,\,\ \big(\X[-1]\star\X\big)/{^{\bot}}\X  \,\,\ @> \approx >> \,\,\  (\X\lsmod)^{\op}
\end{CD}
\]     
and $\T/{^{\bot}}\X$ is abelian if and only if for any object $A \in \T$, the monomorphism $\underline{f}^{A}_{0} \colon \unA \lxr \unX^{A}_{0}$ is a kernel in $\T/{^{\bot}}\X$.  
 \end{thm}
 
 Theorem $4.6$ admits an equivalent formulation as follows. Let $\mathsf{H} \colon \T \lxr \A$ be a homological functor, where $\A$ is abelian. An object $P \in \T$ is called $\mathsf{H}$-{\em projective} if $\mathsf{H}(P)$ is projective in $\A$ and for any object $A \in \A$ the map $\T(P,A) \lxr \A(\mathsf{H}(P), \mathsf{H}(A))$ is invertible. We denote by $\mathcal P(\mathsf{H})$ the full subcategory of $\T$ consisting of the $\mathsf{H}$-projective objects. We say that $\T$ {\em has enough} $\mathsf{H}$-{\em projectives}, if for any object $A \in \T$ there exists a triangle $K \lxr P \lxr A \lxr K[1]$ in $\T$, where $P \in \mathcal P(\mathsf{H})$ and the map $\mathsf{H}(P) \lxr \mathsf{H}(A)$ is an epimorphism.  Finally we say that the $\mathsf{H}$-dimension of $\T$ is at most one, $\mathsf{H}$-$\mathsf{dim}\, \T \leq 1$, if for any object $A\in \T$ there exists a triangle  $P_{1} \lxr P_{0} \lxr A \lxr P_{1}[1]$ in $\T$ where the $P_{i} \in \mathcal P(\mathsf{H})$ and the map $\mathsf{H}(P_{0}) \lxr \mathsf{H}(A)$ is an epimorphism.
  
 \begin{thm} Let $\T$ be a triangulated category and $\mathsf{H} \colon \T \lxr \A$ a homological functor, where $\A$ is abelian with enough projectives. We assume that $\T$ has enough $\mathsf{H}$-projectives and $\mathcal P(\mathsf{H})[1] \subseteq \Ker\mathsf{H}$. If $\Proj\A \subseteq \Image \mathsf{H}$,  
 then the class $\mathcal R$ of regular maps in $\T/\Ker\mathsf{H}$ admits a calculus of left and right fractions and there is an equivalence
 \[
\begin{CD}
(\T/\Ker\mathsf{H})[\mathcal R^{-1}] \, @> \approx >> \, \A
\end{CD}
\] 
\begin{proof} We sketch the proof for the convenience of the reader. Using that $\T$ has enough $\mathsf{H}$-projectives and any projective object of $\A$ lies in the image of $\mathsf{H}$, one shows that $\mathcal P(\mathsf{H})$ is contravariantly finite in $\T$ and the functor $\mathsf{H}$ restricts to an equivalence between $\mathcal P(\mathsf{H})$ and $\Proj\A$. Next using that $\A$ has enough projectives and the fact that $\A = \smod\Proj\A$, we have $\A = \smod\mathcal P(\mathsf{H})$ and then the homological functor $\mathsf{H} \colon \T \lxr \A$ is isomorphic to the restriction functor $\T \lxr \smod\mathcal P(\mathsf{H})$, $A \longmapsto \T(-,A)|_{\mathcal P(\mathsf{H})}$. It follows that $\Ker\mathsf{H} = \mathcal P(\mathsf{H})^{\bot}$ and therefore the condition $\mathcal P(\mathsf{H})[1] \subseteq \Ker\mathsf{H}$ implies that $\mathcal P(\mathsf{H})$ is rigid.  The rest follows from Theorem  $4.6$ with $\X = \mathcal P(\mathsf{H})$. 
\end{proof}
\end{thm} 

\begin{cor} Let $\T$ be a triangulated category and $\mathsf{H} \colon \T \lxr \A$ a homological functor, where $\A$ is abelian. We assume that $\mathsf{H}$ is surjective on objects and $\mathcal P(\mathsf{H})[1] \subseteq \Ker\mathsf{H}$. If $\mathsf{H}$-$\mathsf{dim}\, \T \leq 1$, then $\T/\Ker\mathsf{H}$ is abelian, $\mathcal P(\mathsf{H})^{\bot} = \mathcal P(\mathsf{H})[1]$, so $\mathcal P(\mathsf{H})$ is $2$-cluster tilting, and we have an equivalence $\T/\Ker\mathsf{H} \stackrel{\approx}{\lxr} \A$. 
\begin{proof} Since $\mathsf{H}$-$\mathsf{dim}\, \T \leq 1$ it follows that $\T$ has enough $\mathsf{H}$-projectives, and since $\mathsf{H}$ is surjective on objects, it follows that $\A$ has  enough projectives. Then the assertion follows from Theorems $4.6$ and $4.9$ using the fact that $\mathsf{H}$-$\mathsf{dim}\,\T \leq 1$ implies that $\T = \mathcal P(\mathsf{H}) \star \mathcal P(\mathsf{H})[1]$.
\end{proof} 
 \end{cor}

 \begin{exam} Let $\T$ be a triangulated category which admits all small coproducts and let $\mathcal P$ be a set of compact objects of $\T$ such that $\T(\mathcal P,\mathcal P[n]) = 0$, $\forall n \neq 0$; for instance we may take $\T = {\bf D}(\Mod\Lambda)$, the unbounded derived category of an associative ring $\Lambda$, and $\mathcal P = \{T\}$ is a tilting complex.  Then the set $\mathcal P$ induces a $t$-structure $(\T^{\leq 0}, \T^{\geq 0})$ in $\T$, where $\T^{\geq 0} = \{A \in \T \ | \ \T(\mathcal P,A[n]) = 0, \ \forall n \leq -1\}$ and  $\T^{\leq 0} = {^{\bot}}\T^{\geq 1}$, see \cite{BR}.  Let $\mathcal H = \T^{\leq 0} \cap \T^{\geq 0}$ be the heart of the $t$-structure and $\mathsf{H} \colon \T \lxr \mathcal H$ be the induced homological functor. It is not difficult to  that the assumptions of Theorem $4.9$ are satisfied, hence we have an equivalence $(\T/\Ker\mathsf{H})[\R^{-1}] \, \stackrel{\approx}{\lxr} \, \mathcal H$. 
 \end{exam}

 \subsection{$\T/\X^{\bot}$ as a Semiabelian Category} We have seen that the factor category $\T/\X^{\bot}$ has cokernels. So it is natural to ask whether the category $\T/\X^{\bot}$ has kernels.  In this connection we recall that  an additive category $\A$ is called {\em preabelian} if $\A$ has kernels and cokernels.  A preabelian category $\A$ is called  {\em integral}, see \cite{Rump}, 
if the class of epimorphisms is closed under pull-backs and the class of monomorphisms is closed under push-outs. Finally a preabelian category $\A$ is called {\em semiabelian},  if for any map $f \colon A \lxr B$ in $\A$, the naturally induced map $\widetilde{f} \colon \Coker(\mathsf{ker} f) \, \lxr \, \Ker(\coker f)$ is regular.

\begin{prop}  Let $\T$ be a triangulated category and $\X$ a contravariantly finite rigid subcategory of $\T$. If  $\X^{\bot}$ is contravariantly finite in $\T$, then the stable category $\T/\X^{\bot}$ is semi-abelian and integral.
\begin{proof} Since $\X^{\bot}$ is functorially finite, it follows from Lemma $2.1$ that we have an adjoint pair
\[
\big(\Sigma^{1}_{\X^{\bot}},\Omega^{1}_{\X^{\bot}}\big) \  \colon  \xymatrix@C=3.5pc {\T/\X^{\bot} \,\, \ar@<0.75ex>[r]^-{{\mathsf{}}} & \,\, \ar@<0.3ex>[l]^-{{\,\, \mathsf{}}}\T/\X^{\bot}}
\]
Since the map $h^{0}_{A} \colon A \lxr \Omega^{1}_{\X}(A)[1]$ is a left $\X^{\bot}$-approximation of $A$, by construction, we have $\Sigma^{1}_{\X^{\bot}}(A) =  X^{0}_{A}[1]$.  Since $\X$ is rigid we have $X^{0}_{A}[1] \in \X^{\bot}$ and this means that $\Sigma^{1}_{\X^{\bot}}(A) = 0$ in $\T/\X^{\bot}$. We infer that the suspension functor $\Sigma^{1}_{\X^{\bot}}$ is zero. Clearly then $\Omega^{1}_{\X^{\bot}}$ is zero, hence $\T/\X^{\bot}$ is a pretriangulated category with zero suspension and loop functor. This trivially implies that $\T/\X^{\bot}$ is preabelian. Alternatively,  let $\underline{\alpha} \colon \unA \lxr \unB$ be a map in $\T/\X^{\bot}$. Let 
$(T) :\, \Omega^{1}_{\X^{\bot}}(B)  \lxr Y^{0}_{B} \lxr B \lxr \Omega^{1}_{\X^{\bot}}(B)[1]$ be a triangle in $\T$, where the middle map is a right $\X^{\bot}$-approximation of $B$, and form the homotopy pull-back diagram of $(T)$ along the map $\alpha$:
\begin{equation}
  \xymatrix@C=1.5cm{
    \Omega^{1}_{\X^{\bot}}(B)  \ar[r]^{}  \ar@{=}[d]  &  
   K \ar[r]^{\kappa}\ar[d]^{}&
    A \ar[d]^{\alpha} \ar[r]^{}&
    \Omega^{1}_{\X^{\bot}}(B)[1] \ar@{=}[d] \\
    \Omega^{1}_{\X^{\bot}}(B) \ar[r]^{}& Y^{0}_{B} \ar[r]^{}& B \ar[r]^{} & \Omega^{1}_{\X^{\bot}}(B)[1]}
\end{equation} 
Since $\Omega^{1}_{\X^{\bot}}(B)$ lies in $\X^{\bot}$, it is easy to see that the map $\underline{\kappa} \colon \unK \lxr \unA$ is the kernel of $\underline{\alpha}$ in $\T/\X^{\bot}$. Hence $\T/\X^{\bot}$ has kernels and therefore $\T/\X^{\bot}$ is preabelian since $\T/\X^{\bot}$ has cokernels by Lemma $4.2$.  
By Theorem $4.6$ the localization functor $\mathsf{P} \colon  \T/\X^{\bot} \lxr (\T/\X^{\bot})[\mathcal R^{-1}]$ is faithful and preserves  kernels and cokernels. By a result of Rump, see \cite[Proposition 7]{Rump}, any preabelian category which admits a kernel and cokernel preserving faithful functor into an abelian category is integral.  We infer that $\T/\X^{\bot}$ is integral and then $\T/\X^{\bot}$ is semiabelian by \cite[Corollary 1]{Rump}.   
\end{proof}
\end{prop}

Next we show that $\X^{\bot}$ is contravariantly finite in $\T$ if $\X$ is covariantly finite and $\T$ has Serre duality.

\begin{lem} Assume that $\T$ admits a Serre functor $\mathbb S$ and  let $\X$ be a functorially finite rigid subcategory of $\T$. Then the subcategories $\X^{\bot}$ and ${^{\bot}}\X$ are functorially finite in $\T$. Moreover the semiabelian category $\T/\X^{\bot}$ has enough projectives and injectives and: \, $\Proj\T/\X^{\bot} = \X/\X^{\bot} \approx \X$ \, and \,   $\Inj\T/\X^{\bot} = \mathbb S(\X)/\X^{\bot} \approx \mathbb S(\X)$.   
\begin{proof} Since $\X$ is functorially finite, Lemma $2.1$ implies that $\X^{\bot}$ is covariantly finite and ${^{\bot}}\X$ is contravariantly finite. Since we always have $\X^{\bot} = \mathbb{S}\big({^{\bot}}\X\big)$ and $\mathbb{S}^{-1}\big(\X^{\bot}\big) = {^{\bot}}\X$, 
and ${^{\bot}}\X$ is contravariantly finite, it follows that so is $\mathbb{S}\big({^{\bot}}\X\big) = \X^{\bot}$. Then by Proposition $4.12$ the category $\T/\X^{\bot}$ is semiabelian.  Since by Lemma $4.2$, $\T/\X^{\bot}$ has enough projectives, it remains to show that  $\T/\X^{\bot}$ has enough injectives. This follows from Lemma $4.5$  and the fact that, since $\mathbb S({^{\bot}}\X) = \X^{\bot}$, the Serre functor $\mathbb S$ induces an equivalence $\mathbb S \colon \T/{^{\bot}}\X \, \lxr \, \T/\X^{\bot}$. 
\end{proof}
\end{lem}

\begin{rem} Theorem $4.6$ was proved using different methods, by Buan and Marsh in \cite{BM1, BM2} in case $\T$ is a $k$-linear triangulated category, with finite-dimensional Hom-spaces over a field $k$, satisfying Serre duality and $\X = \add T$ is the additive closure of a rigid object $T$; in this case $\X$ is always functorially finite and the category of coherent functors $\smod\X$ is the category of finite-dimensional modules over the endomorphism algebra of $T$.  
\end{rem}

However the following presents a situation where the setting and the methods of \cite{BM1, BM2} are not applicable. 

\begin{exam} Let $\T$ be a triangulated category with all small coproducts and $T$ a rigid object of $\T$. Let $\Add T$ be the full subcategory of $\T$ consisting of the direct factors of small coproducts of copies of $T$. It is well-known that $\Add T$ is contravariantly finite in $\T$.  We assume that $T$ is {\em self-compact}, i.e. the functor $\T(T,-)$ commutes with all small coproducts of copies of $T$.  It is easy to see that there is an equivalence $\smod\Add T \approx \Mod\End_{\T}(T)$, so the functor $\mathsf{H} = \T(T,-)  \colon \T \lxr \Mod\End_{\T}(T)$ induces an equivalence between $\Add T$ and the category of projective modules over the endomorphism ring of $T$. Since self-compactness of $T$ implies that $\Add T$ is rigid, it follows by Theorem $4.6$ that we have an equivalence $(\T/T^{\bot})[\mathcal R^{-1}] \approx \Mod\End_{\T}(T)$, where $T^{\bot} = \{A \in \T \ | \ \T(T,A) = 0\} =  (\Add T)^{\bot}$,  and $\mathcal R$ is the class of regular maps in $\T/T^{\bot}$. For instance let $\T = {\bf D}(\Mod\Lambda)$ for a ring $\Lambda$, and let $T = \Lambda$ be concentrated in degree zero. Then $\{{\bf D}(\Mod\Lambda)/\Lambda^{\bot}\}[\mathcal R^{-1}] \approx \Mod\Lambda$, where $\mathcal R$ is the class of maps $A \lxr B$ in ${\bf D}(\Mod\Lambda)/\Lambda^{\bot}$ such that the induced map $\mathrm{H}^{0}(A) \lxr \mathrm{H}^{0}(B)$ on $0$th cohomology is invertible. 
\end{exam}

\subsection{Fullness of $\mathsf{H} \colon \T \lxr \smod\X$} Recall that if $\X$ is a contravariantly finite rigid subcategory of $\T$, then the homological functor $\mathsf{H} \colon \T \lxr \smod\X$ is surjective on objects: any coherent functor $F \colon \X^{\op} \lxr \ab$ with presentation $(-,X^{1}) \lxr (-,X^{0}) \lxr F \lxr 0$ is isomorphic to $\mathsf{H}(A)$, where $A$ is a cone of $X^{1} \lxr X^{0}$ in $\T$. Clearly $\mathsf{H}$ is faithful iff $\T = \X$ iff $\T = \{0\}$. We close this section by observing that the above results allow us to characterize when the homological functor $\mathsf{H}$ is full, see also \cite{B}.  This characterization is similar in spirit to Theorem $3.3$ of \cite{BK} where $\X$ is assumed to be closed under all shifts instead of being rigid.   First we need the following notation:  if $\A, \B$ are classes of objects of $\T$, then we denote by $\A \oplus \B$ the full subcategory of $\T$ consisting of the direct summands of direct sums $A \oplus B$, where $A \in \A$ and $B \in \B$.

\begin{prop} Let $\T$ be a triangulated category and $\X$ a rigid subcategory of $\T$.  
\begin{enumerate} 
\item If $\X$ is contravariantly finite, then the following statements are equivalent:
\begin{enumerate} 
\item The homological functor $\mathsf{H} \colon \T \lxr \smod\X$ is full.  
\item $\T = (\X \star \X[1]) \oplus \X^{\bot}$. 
\item $\T/\X^{\bot}$ is abelian. 
\end{enumerate} 
\item If $\X$ is covariantly finite, then the following statements are equivalent:
\begin{enumerate} 
\item The cohomological functor $\mathsf{H} \colon \T \lxr \X\lsmod$ is full.  
\item $\T = (\X[-1] \star \X) \oplus {^{\bot}}\X$.
\item $\T/{^{\bot}}\X$ is abelian. 
\end{enumerate} 
\end{enumerate}
\begin{proof} (i) By Theorem $4.6$, the functor $\mathsf{H}$ admits a factorization $\mathsf{H} = \underline{\mathsf{H}} \circ \mathsf{R}^{*} \circ \mathsf{P} \circ \pi$ where the functors $\underline{\mathsf{H}}$ and $\mathsf{R}^{*}$ are equivalences and the functor $\pi$ is full. It follows from this that $\mathsf{H}$ is full if and only if so is $\mathsf{P}$. Since $\mathsf{P}$ is faithful and surjective on objects, we infer that $\mathsf{H}$ is full if and only of $\mathsf{P}$ is an equivalence and this happens if and only if $\T/\X^{\bot}$ is abelian. Hence (a) $\Leftrightarrow$ (c) and clearly (b) implies (c) since then $\T/\X^{\bot} = (\X\star\X[1])/\X^{\bot}$ which is abelian by Lemma $4.4$. Finally assuming (a), we have seen that $\mathsf{P}$ is an equivalence and this implies that the inclusion $\mathsf{L} \colon \T/\X^{\bot} \lxr (\X\star \X[1])/\X^{\bot}$ is an equivalence. Then for any object $A \in \T$ the map $\underline{\mathsf{Cell}_{1}(A)} \lxr \unA$ is invertible in $\T/\X^{\bot}$,  hence  $A \oplus K \cong \mathsf{Cell}_{1}(A) \oplus L$ for some objects $K, L \in \X^{\bot}$. Since $\mathsf{Cell}_{1}(A) \in \X\star\X[1]$ we have that $A$ lies in $(\X \star \X[1]) \oplus \X^{\bot}$ and part (b) follows.  

Part (ii) follows by duality using Theorem $4.8$. 
\end{proof}
\end{prop}

As a direct consequence of  Proposition $4.16$ and (the proof of) Lemma $4.13$ we have the following.

\begin{cor} Let $\T$ be a Hom-finite $k$-linear triangulated category with Serre duality over a field $k$. If $\X$ is a functorially finite rigid subcategory of $\T$, then $\T/\X^{\bot}$ is abelian if and only if $\T/{^{\bot}}\X$ is abelian. 
\end{cor}

Recall  that $\X$ is called $n$-{\em rigid}, where $n \geq 1$, if \, $\T(\X,\X[i]) = 0$, $1\leq i \leq n$. 

\begin{cor} Let $\T$ be a triangulated Krull-Schmidt category and $\X$ a contravariantly finite $n$-rigid subcategory of $\T$, where $n \geq 2$. If the functor $\mathsf{H} \colon \T \lxr \smod\X$ is full, then $\T(\X,\X[-i]) = 0$, $1\leq i \leq n-1$. 
\begin{proof} Let $X \in \X$ and $1\leq i \leq n-1$. Since $\T$ is Krull-Schmidt, by Proposition $4.16$, the object $X[-i]$  admits a direct sum decomposition $X[-i] = A \oplus B$, where $A \in \X\star \X[1]$ has no non-zero direct summands from $\X^{\bot}$, and $B \in \X^{\bot}$, and then clearly  $A = X^{\prime}[-i]$, for some $X^{\prime} \in \X$. Since $A$ lies in $\X\star\X[1]$, there exists a triangle  $X_{1} \lxr X^{\prime}[-i] \lxr X_{2}[1] \lxr X_{1}[1]$, where the $X_{i}$'s lie in $\X$, and therefore a triangle 
$X_{1}[i] \lxr X^{\prime} \lxr X_{2}[i+1] \lxr X_{1}[i+1]$. Since $i \leq n-1$ and $\X$ is $n$-rigid, applying the homological functor $\mathsf{H}$ to this triangle, we have $\mathsf{H}(X^{\prime}) = 0$, i.e. $X^{\prime} = 0$, and therefore $A = 0$. We infer that $X[-i] = B \in \X^{\bot}$, for $1\leq i \leq n-1$. 
\end{proof}
\end{cor} 

We note  that important examples of $n$-rigid subcategories are the $(n+1)$-cluster tilting subcategories in the sense of \cite{KR1, IY}.  Recall that $\X$ is $(n+1)$-{\em cluster tilting subcategory} of $\T$, if  $\X$ is functorially finite in $\T$ and:
\[
\big\{A \in \T \,\, | \,\, \T(\X,A[k]) = 0, \, 1\leq k \leq n\big\} \,\ = \,\ \X \,\  = \,\ \big\{A \in \T \,\, | \,\, \T(A,\X[k]) = 0, \, 1\leq k \leq n\big\}
\]
Note that, by \cite{KZ}, $\X$ is $2$-cluster tilting if $\X$ is contravariantly finite and $\X^{\bot} = \X[1]$. Hence if $\X$ is an $(n+1)$-cluster tilting subcategory of $\T$, where $n\geq 2$, and the functor $\mathsf{H}$ is full, then $\T(\X,\X[-i]) = 0$, for $1\leq i \leq n-1$.

The following, based on an example of Iyama, see \cite{KR1}, shows that in general: $(\alpha)$ the functor $\mathsf{H}$ is not full, and $(\beta)$ the factor category $\T/\X^{\bot}$ is not abelian, even for $3$-cluster tilting subcategories.

 \begin{exam}  Let $k$ be a field with $\mathsf{char}(k) \neq 3$  and let $\omega$ be a primitive third root of unity.  The cyclic group $G = \mathbb Z/3\mathbb Z = \langle g\rangle$ acts on the power series algebra $S = k[[t,x,y,z]]$ by $g \cdot t = \omega t, \, g\cdot x = \omega x, \, g \cdot y = \omega^{2} y, \, g\cdot z = \omega^{2} z$. Then the algebra of invariants $S^{G}$ is an isolated singularity  and the stable triangulated category $\T = \underline{\E}$, where $\E = \mathsf{CM}(S^{G})$ is the Frobenius category of Cohen-Macaulay modules over $S^{G}$, contains $\X = \add\underline{S}$ as a $3$-cluster tilting subcategory, see \cite{KR1}.   Clearly then $\smod\X \approx \smod\uEnd(S)$ is the stable endomorphism algebra of $S$, and $\mathsf{H} \colon \T \lxr  \smod\uEnd(S)$ is given by $\uHom(S,-)$. It is not difficult to check, see \cite{B}, that  no nonzero  direct summand $X$ of $\underline{S}$ satisfies the vanishing condition $\T(X,X[-1]) = 0$ so by Corollary $4.18$ the functor $\mathsf{H}$ is not full.  This reflects the fact that the endomorphism algebra   $\uEnd(S)$ is not Gorenstein, see \cite{KR1, B}, or Theorem $5.2$ below, for more details. By Proposition $4.12$,  the stable category $\T/\X^{\bot}$ is semi-abelian and integral but not abelian. 
\end{exam}

\begin{rem} If $\X$ is contravariantly finite rigid, then by Proposition $4.16$ fullness of $\mathsf{H} \colon \T \lxr \smod\X$ implies that $\T/\X^{\bot}$ is abelian and $\T/\X^{\bot} \approx \smod\X$. Clearly in this case the stable category $\T/\X^{\bot}$ is identified with the stable category $\T/\mathsf{Gh}_{\X}(\T)$ defined by the ideal of $\X$-ghost maps and therefore we have an equivalence $\T/\mathsf{Gh}_{\X}(\T) \approx \smod\X$ induced by $\mathsf{H}$.  Note that in this case the functor $\mathsf{H}$ does not necessarily reflect isomorphisms and the global dimension of $\smod\X$ can be arbitrary.  On the other hand if $\X$ is contravariantly finite, not rigid but stable, i.e. $\X = \X[1]$, then, as follows easily from \cite[Theorem 5.3]{B:3cats}, if the functor $\mathsf{H}$ reflects isomorphisms, then: 
\[\mathsf{H} \ \  \text{induces an equivalence} \ \  \T/\mathsf{Gh}_{\X}(\T) \approx \smod\X \ \ \ \text{if and only if} \ \ \ \smod\X \ \  \text{is hereditary}\]  
This covers for instance \cite[Theorem $3.1$]{KYZ}. 

If $\T = {\bf D}^{b}(\A)$ is the bounded derived category of an abelian category $\A$ with enough projectives, and  $\X$ is the full subcategory of Cartan-Eilenberg projective complexes, see \cite{B:3cats}, then the functor $\mathsf{H}^{*} \colon {\bf D}^{b}(\A) \lxr \mathsf{Gr}^{b}(\A)$ associated to $\X$  is the usual cohomology functor, where $\smod\X = \mathsf{Gr}^{b}(\A)$ is the  subcategory of bounded objects of the graded category $\mathsf{Gr}(\A)$ of $\A$.  Since $\mathsf{H}^{*}$ reflects isomorphisms, $\mathsf{H}^{*}$ induces an equivalence ${\bf D}^{b}(\A)/\mathsf{Gh}_{\X}({\bf D}^{b}(\A) ) \approx \mathsf{Gr}^{b}(\A)$ if and only if $\mathsf{Gr}^{b}(\A)$, or equivalently $\A$, is hereditary.  Note that $\mathsf{Gh}_{\X}({\bf D}^{b}(\A))$ is the ideal of classical ghost maps, i.e. the maps $f$ in ${\bf D}^{b}(\A)$ such that $\mathsf{H}^{*}(f) = 0$.  
  \end{rem}

\section{Abelian Factors: The Gorenstein and the Calabi-Yau Property}
  
Let $\X$ be a contravariantly finite rigid subcategory of $\T$.  In this section we show that when the factor category $\T/\X^{\bot}$ is abelian, or equivalently the associated functor $\mathsf{H} \colon \T \lxr \smod\X$ is full, then the category $\smod\X$ is $1$-Gorenstein. Moreover if $\T$ is $2$-Calabi-Yau, then we show that the stable triangulated category modulo projectives  of the Cohen-Macaulay objects of $\smod\X$ is $3$-Calabi-Yau. Thus we extend  basic results of Keller-Reiten \cite{KR1}, see also Koenig-Zhu \cite{KZ}, concerning $2$-cluster tilting objects to the more general setting of fully rigid objects in the sense of the following definition which is inspired by the results of subsection $4.4$:

\begin{defn} A contravariantly finite rigid subcategory $\X$ of a triangulated category  $\T$ is called {\em fully rigid} if the associated homological functor $\mathsf{H} \colon \T \lxr \smod\X$ is full. 
\end{defn}

By Example $4.19$ there exist contravariantly finite rigid subcategories, even $3$-cluster tilting, which are not fully rigid. On the other hand, it is clear that any $2$-cluster tilting subcategory is fully rigid. Note that in contrast to $2$-cluster tilting, or maximal rigid,  subcategories, which are defined in terms of vanishing properties inside  $\T$, fully rigid subcategories  are defined in terms of properties of the associated homological functor from $\T$  to the category of coherent functors or in terms of the associated factor category $\T/\X^{\bot}$.

  Recall that an abelian category $\A$ is called ($n$-){\em Gorenstein} if $\A$ has enough projectives and enough injectives and its {\em Gorenstein dimension} 
\[\Gor\mathsf{dim}\A : =\max\{\spli\A, \silp\A\}\]
 is finite ($\leq n$), where $\spli\A = \sup\{\pd I \ | \ I \in \Inj\A \}$ and $\silp\A = \{\id P \ | \ P \in \Proj\A\}$.   An Artin algebra $\Lambda$ is called $n$-Gorenstein if $\smod\Lambda$ is $n$-Gorenstein; this is equivalent to $\Gor\mathsf{dim}\Lambda = \{\id{_{\Lambda}}\Lambda,\id\Lambda_{\Lambda}\} \leq n$. 
  
  \begin{thm} Let $\T$ be a Hom-finite  $k$-linear triangulated category $\T$ over a field $k$ and $\X$ a  fully rigid subcategory of $\T$. Assume that either $(\alpha)$ $\X = \add T$, or else $(\beta)$ $\X$ is covariantly finite and $\T$ has Serre duality.  
\begin{enumerate}
\item  The category $\smod\X$ of coherent functors over $\X$ is $1$-Gorenstein.
\item  Either $\gd\smod\X \leq 1$  or $\gd\smod\X = \infty$. 
\item In case $(\beta)\!:$ \, $\smod\X$ is Frobenius if and only if $\X = \mathbb S(\X)$.
 \end{enumerate}
  \begin{proof} First observe that in both cases the subcategory $\X$ is functorially finite in $\T$ and $\smod\X$ is a dualizing $k$-variety in the sense of \cite{AR}. Indeed in case $(\alpha)$ $\smod\X$ is the category $\smod\Lambda$ of finitely generated $\Lambda$-modules over the endomorphism algebra $\Lambda := \End_{\T}(T)$, and in case $(\beta)$ this follows easily from the existence of a Serre functor $\mathbb S \colon \T \lxr \T$ of $\T$. Hence in both cases the duality functor $\mathsf{D} = \Hom_{k}(-,k)$ induces a duality between $\X\lsmod$ and $\smod\X$, in particular $\smod\X$ has enough projectives and enough injectives.  

  \underline{\textsf{Case}} $(\alpha)$:  We show first that $\Inj\smod\X \subseteq \add \mathsf{H}(\X^{\bot}[1])$. Let $F$ be an injective $\Lambda$-module. Since the functor $\mathsf{H}$ is surjective on objects, we have $F \cong \mathsf{H}(A)$ for some object $A\in \T$. Let 
\[
\Omega^{1}_{\X}(A[-1]) \, \lxr \, X^{0}_{A[-1]} \, \lxr \,  A[-1] \, \lxr \, \Omega^{1}_{\X}(A[-1])[1]
\] be a triangle, where the middle map is a right $\X$-approximation of $A[-1]$. By Lemma $2.2$, we have $\Omega^{1}_{\X}(A[-1])[1] \in \X^{\bot}$. Applying the homological functor $\mathsf{H}$ to the triangle  
$
X^{0}_{A[-1]}[1] \, \lxr \,  A \,  \lxr \, \Omega^{1}_{\X}(A[-1])[2] \, \lxr \,  X^{0}_{A[-1]}[2]
$ and using that $\X$ is rigid, we have a monomorphism $\mathsf{H}(A) \lxr \mathsf{H}\big(\Omega^{1}_{\X}(A[-1])[2]\big)$ which splits since $\mathsf{H}(A)$ is injective. Since  $\Omega^{1}_{\X}(A[-1])[2]$ lies in $\X^{\bot}[1]$, it follows that any injective object of $\smod\X$ is a direct summand of an object from $\mathsf{H}(\X^{\bot}[1])$, i.e. $\Inj\smod\X \subseteq \add \mathsf{H}(\X^{\bot}[1])$. 

Next we show that $\pd\mathsf{H}(A) \leq 1$, $\forall A \in \X^{\bot}[1]$. Since $\mathsf{H}$ is full, by Proposition $4.16$  we can write $A = A_{1} \oplus A_{2}$, where $A_{1} \in \X\star \X[1]$ and $A_{2} \in \X^{\bot}$. Then clearly $\mathsf{H}(A) \cong \mathsf{H}(A_{1})$ and $A_{1} \in \X^{\bot}[1]$, or equivalently $A_{1}[-1] \in \X^{\bot}$, as a direct summand of $A$. Since $A_{1} \in \X\star \X[1]$, there exists a triangle $X^{1} \lxr X^{0} \lxr A_{1} \lxr X^{1}[1]$, where the $X^{i}$ lie in $\X$. Applying the functor $\mathsf{H}$ and using that $\mathsf{H}(A_{1}[-1]) = 0$, we have a short exact sequence $0 \lxr \mathsf{H}(X^{1}) \lxr \mathsf{H}(X^{0}) \lxr \mathsf{H}(A_{1}) \lxr 0$ which shows that $\pd\mathsf{H}(A_{1}) = \pd\mathsf{H}(A) \leq 1$. Hence any injective object of $\smod\X$ has projective dimension $\leq 1$ or equivalently $\pd\mathsf{D}({_{\Lambda}}\Lambda) \leq 1$, and therefore  $\id{_{\Lambda}}\Lambda \leq 1$. However for Artin algebras $\Lambda$, this already implies that  $\id\Lambda_{\Lambda} \leq 1$. We recall the argument. Since $\pd\mathsf{D}({_{\Lambda}}\Lambda) \leq 1$, $\mathsf{D}({_{\Lambda}}\Lambda)$ is a partial tilting  $\Lambda$-module and then by a classical result of Bongartz \cite{Bongartz}, $\mathsf{D}({_{\Lambda}}\Lambda)$ is a direct summand of a tilting $\Lambda$-module.  This implies that $\mathsf{D}({_{\Lambda}}\Lambda)$ is a tilting $\Lambda$-module since it has the same number of indecomposable direct summands as the regular module $\Lambda$. Hence there exists a short exact sequence $0 \lxr \Lambda_{\Lambda} \lxr I^{0} \lxr I^{1} \lxr 0$, where the $I^{i}$ lie in $\add\mathsf{D}({_{\Lambda}}\Lambda)$ and therefore $\id\Lambda_{\Lambda} \leq 1$. We infer that $\mathsf{G}$-$\mathsf{dim}\Lambda \leq 1$ and  $\Lambda$  is $1$-Gorenstein.

   \underline{\textsf{Case}} $(\beta)$: By Serre duality we have natural isomorphisms 
\[
\mathsf{H}(\mathbb S(?)) \, = \, \T(-,\mathbb S(?))|_{\X} \, \cong \,  \mathsf{D}\T(\X,?)|_{\X} \, = \, \mathsf{D}\mathsf{H}^{\op}(?)
\]
Hence $\mathsf{D}\mathsf{H}^{\op}(\X) = \mathsf{H}(\mathbb{S}(\X))$ is the full subcategory of injective objects of $\smod\X$. If $\alpha \colon \mathbb S(X_{1}) \lxr \mathbb S(X_{2})$ is $\X$-ghost, where the $X_{i}$ lie in $\X$, then $\alpha$ factorizes through $\X^{\bot}$. Since, by Serre duality,  $\T(\X^{\bot},\mathbb S(\X))$ is isomorphic to $\mathsf{D}\T(\X,\X^{\bot}) = 0$, it follows that the map $\alpha$ is zero. Since $\mathsf{H}$ is full, this means that $\mathsf{H}$ induces an equivalence $\mathbb S(\X) \approx \Inj\smod\X$. Since $\X$ is fully rigid, by Proposition $4.16$ any object $A$ of $\T$ can be written as: 
 \[
  A = A_{1} \oplus A_{2} = A^{1} \oplus A^{2}, \ \  \ A_{1} \in \X\star \X[1], \ \ A_{2} \in \X^{\bot}  
\ \  \text{and} \ \ \ A^{1} \in \X[-1] \star \X, \   A^{2} \in {^{\bot}}\X \eqno (\dag)
 \]
Let $\mathsf{H}(\mathbb S(X))$ be an injective object of $\smod\X$, where $X \in \X$. By $(\dag)$ we can write $\mathbb S(X) = A_{1} \oplus A_{2}$, where $A_{1} \in \X\star \X[1]$ and $A_{2} \in \X^{\bot}$. Then clearly $\mathsf{H}(\mathbb S(X)) \cong \mathsf{H}(A_{1})$ and $A_{1} \in \mathbb S(\X)$ as a direct summand of $A$. Note that  $\mathsf{H}(A_{1}[-1]) = 0$ since $\mathsf{H}(\mathbb S(\X)[-1]) \cong \mathsf{D}\mathsf{H}^{\op}(\X[-1]))=0$. Since $A_{1} \in \X\star \X[1]$, there exists a triangle 
\[
X^{1} \lxr X^{0} \lxr A_{1} \lxr X^{1}[1]
\]
 where the $X^{i}$ lie in $\X$. Applying the functor $\mathsf{H}$  we have a short exact sequence 
\[
0 \lxr \mathsf{H}(X^{1}) \lxr \mathsf{H}(X^{0}) \lxr \mathsf{H}(A_{1}) \lxr 0
\]
 which shows that $\pd\mathsf{H}(A_{1}) = \pd\mathsf{H}(\mathbb S(X)) \leq 1$. We infer that any injective object of $\smod\X$ has projective dimension $\leq 1$, i.e. $\spli\smod\X \leq 1$. 

 Now let $\mathsf{H}(X)$ be a projective object of $\smod\X$, where $X \in \X$.  By $(\dag)$ we can write $\mathbb S^{-1}(X) = A^{1} \oplus A^{2}$, where $A^{1} \in \X[-1]\star\X$ and $A^{2} \in {^{\bot}}\X$. Then $X = \mathbb S(A^{1}) \oplus \mathbb S(A^{2})$. Since $\mathsf{H}(\mathbb S(A^{2})) \cong \mathsf{D}\mathsf{H}^{\op}(A^{2}) = 0$, it follows that $\id\mathsf{H}(X) = \id\mathsf{H}(\mathbb S(A^{1}))$. Since $A^{1} \in \X[-1]\star\X$, there exists a triangle $X_{1}[-1] \lxr A^{1} \lxr X_{0} \lxr X_{1}$, where the $X_{i}$ lie in $\X$, and therefore a triangle 
\[
\mathbb S(X_{1}[-1]) \lxr \mathbb S(A^{1}) \lxr \mathbb S(X_{0}) \lxr \mathbb S(X_{1})
\]
 Observe that 
$\mathsf{H}(\mathbb S(X_{1}[-1])) \cong \mathsf{D}\mathsf{H}^{\op}(X_{1}[-1]) = 0$.  Since $A_{1}$ is a direct summand of $\mathbb S^{-1}(X)$, it follows that  $\mathbb S(A^{1}[1])$ lies in $\mathbb S\mathbb S^{-1}(\X[1]) = \X[1]$, hence $\mathsf{H}(\mathbb S(A^{1}[1])) = 0$.  Hence applying the homological functor $\mathsf{H}$ to the above triangle we have an exact sequence 
\[
0 \lxr \mathsf{H}(\mathbb S(A^{1})) \lxr \mathsf{H}(\mathbb S(X_{0})) \lxr \mathsf{H}(\mathbb S(X_{1})) \lxr 0
\]
which shows that $\id \mathsf{H}(\mathbb S(A^{1})) = \id\mathsf{H}(X) \leq 1$, i.e. $\silp\smod\X \leq 1$. 
 We infer that $\Gor\mathsf{dim}\smod\X \leq 1$ and  $\smod\X$  is $1$-Gorenstein.
 
 Since for a Gorenstein category $\A$ we have either $\gd\A = \infty$ or else $\gd\A = \Gor\dim\A$ if $\gd\A < \infty$, and since $\inj\smod\X = \mathbb S(\X)$, parts (ii) and (iii) are left to the reader (see \cite{KR1}, \cite{KZ} in  case  $\X$ is $2$-cluster tilting). 
\end{proof}
\end{thm}

Let $\A$ be an abelian category. Recall that an acyclic complex $P^{\bullet}$ consisting of projective objects of $\A$ is called {\em totally acyclic} if the induced complex $\A(P^{\bullet},P)$ is acyclic, for any projective object $P$ of $\A$. An object $A$ is called {\em Gorenstein-projective}, or {\em Cohen-Macaulay}, if  $A$ is isomorphic to $\Ker d^{0}_{P^{\bullet}}$ for a totally acyclic complex $(P^{\bullet},d^{\bullet}_{P^{\bullet}})$ of projectives of $\A$. We denote by $\mathsf{CM}(\A)$ the full subcategory of Cohen-Macaulay objects of $\A$. It is an exact Frobenius subcategory of $\A$ with $\Proj\A$ as the full subcategory of projective-injective objects. It follows that $\underline{\mathsf{CM}}(\A)$ is a triangulated category with suspension functor $[1]_{\underline{\mathsf{CM}}(\A)} = \Omega^{-1}$. It is easy to see that for a $1$-Gorenstein abelian category $\A$, $\mathsf{CM}(\A)$ consists of the subobjects of the projective objects \cite{B:gor}. If $\X \subseteq \T$ is contravariantly finite, we use the notation $\mathsf{CM}(\X)$ for the category of Cohen-Macaulay objects of $\smod\X$ and $\underline{\mathsf{CM}}(\X)$ for the corresponding stable triangulated category modulo projectives. 

For the notion of Auslander's transpose duality functor $\mathsf{Tr} \colon \umod\X \lxr \umod\X^{\op}$ we refer to \cite{AR}.   

\begin{lem} Let $\T$ be a Hom-finite triangulated category over a field $k$ which admits a Serre functor $\mathbb S$. Let $\X$ be a functorially finite rigid subcategory of $\T$. Then for any object $A\in \X \star \X[1]$: 
\[
\mathsf{DTr}(\mathsf{H}(A)) \,  \cong \,  \mathsf{H}(\mathbb S(A)[-1])
\]  
\begin{proof} Let $X^{1} \stackrel{\alpha}{\lxr} X^{0} \lxr A \lxr X^{1}[1]$ be a triangle, where the $X^{i}$ lie in $\X$. Then we have a projective presentation $\mathsf{H}(X^{1}) \lxr \mathsf{H}(X^{0}) \lxr \mathsf{H}(A) \lxr 0$ of $\mathsf{H}(A)$ in $\smod\X$  and therefore a projective presentation $\mathsf{H}^{\op}(X^{0}) \lxr \mathsf{H}^{\op}(X^{1}) \lxr \mathsf{Tr}\mathsf{H}(A) \lxr 0$ of $\mathsf{Tr}\mathsf{H}(A) := \Coker(\mathsf{H}^{\op}(\alpha))$ in $\smod\X^{\op}$. Applying $k$-duality we have an exact sequence $0 \lxr \mathsf{DTr}\mathsf{H}(A) \lxr  \mathsf{D}\mathsf{H}^{\op}(X^{1}) \lxr  \mathsf{D}\mathsf{H}^{\op}(X^{0})$ in $\smod\X$. Using Serre duality we have $\mathsf{D}\mathsf{H}^{\op}(X) \cong \mathsf{D}\T(X,-)|_{\X} \cong \T(-,\mathbb S(X))|_{\X} = \mathsf{H}(\mathbb S(X))$, $\forall X \in \X$.  Therefore the above exact sequence is isomorphic to the exact sequence  $0 \lxr \mathsf{DTr}\mathsf{H}(A) \lxr \mathsf{H}(\mathbb S(X^{1})) \lxr \mathsf{H}(\mathbb S(X^{0}))$. Since the Serre functor $\mathbb S$ is triangulated, we have a triangle $\mathbb S(X^{1}) \lxr \mathbb S(X^{0}) \lxr \mathbb S(A) \lxr \mathbb S(X^{1})[1]$ in $\T$.  Using that $\mathsf{H}(\mathbb S(X^{0})[-1]) \cong \mathsf{D}\mathsf{H}^{\op}(X^{0}[-1]) = \mathsf{D}\T(X^{0}[-1],\X) = 0$, applying $\mathsf{H}$ to the last triangle we have an exact sequence $0 \lxr \mathsf{H}(\mathbb S(A)[-1]) \lxr \mathsf{H}(\mathbb S(X^{1})) \lxr \mathsf{H}(\mathbb S(X^{0}))$. Hence we have an isomorphism $\mathsf{DTr}\mathsf{H}(A) \cong \mathsf{H}(\mathbb S(A)[-1])$.
\end{proof}
\end{lem}

\begin{thm} Let $\T$ be a $2$-Calabi-Yau triangulated category over a field $k$, and $\X$ a functorially finite fully rigid subcategory of $\T$. Then the stable triangulated category $\underline{\mathsf{CM}}(\X)$ is $3$-Calabi-Yau.
\begin{proof} We divide the proof into three steps:
\begin{enumerate} 
\item[{\em Step 1}:]  For any object $A\in \X\star\X[1]$ and any Cohen-Macaulay object $\mathsf{H}(B)$ there is a natural isomorphism:
\begin{equation}
\begin{CD}
\mathsf{D}\uHom\big(\mathsf{H}(A),\mathsf{H}(B)\big) \ @> \cong >> \  \uHom\big(\Omega^{1}\mathsf{H}(B),\mathsf{H}(A[1])\big)
\end{CD}
\end{equation}

{\em Proof.} As in the proof of Theorem $5.2$, $\smod\X$ is a dualizing $k$-variety in the sense of \cite{AR} and therefore it has Auslander-Reiten duality (we refer to \cite{AR} for all the relevant notions):   
\[
\begin{CD}
\mathsf{D}\uHom\big(\mathsf{H}(A),\mathsf{H}(B)\big) \  @> \cong >> \ \ {\Ext}^{1}\big(\mathsf{H}(B),\mathsf{DTr}\mathsf{H}(A)\big) 
\end{CD}
\eqno (\dag\dag)
\]
Since $\T$ is $2$-Calabi-Yau, by Lemma $5.3$ we have an isomorphism $\mathsf{DTr}\mathsf{H}(A) \cong \mathsf{H}(A[1])$. Hence Auslander-Reiten formula $(\dag\dag)$ gives the desired isomorphism $(5.1)$:
\[
\begin{CD}
\mathsf{D}\uHom\big(\mathsf{H}(A),\mathsf{H}(B)\big) \  @> \cong >> \  {\Ext}^{1}\big(\mathsf{H}(B),\mathsf{H}(A[1])\big) \  @> \cong >>  \  \uHom\big(\Omega^{1}\mathsf{H}(B),\mathsf{H}(A[1])\big)
\end{CD}
\]
where the last isomorphism holds since $\mathsf{H}(B)$ is Cohen-Macaulay. 

\smallskip
 
\item[{\em Step 2}:] For any Cohen-Macaulay object $\mathsf{H}(A)$ of $\smod\X$ we have an isomorphism in $\umod\X$: 
\begin{equation}
\begin{CD}
\Omega^{2}\mathsf{H}(A[1]) @> \cong >>  \mathsf{H}(A) \ \ \ \ \ \text{\em or equivalently} \ \ \ \ \  \Omega^{1}\mathsf{H}(A[1]) @> \cong >>  \Omega^{-1}\mathsf{H}(A)
\end{CD}
\end{equation}

{\em Proof.}  By Proposition $4.16$ we can write $A[1] = M \oplus N$, where $M\in \X\star \X[1]$ and $N \in \X^{\bot}$.  
Let $X^{1} \lxr X^{0} \lxr M \lxr X^{1}[1]$ be a triangle in $\T$ where the $X^{i}$ lie in $\X$.  Then we have the triangle 
\[
X^{1} \oplus N[-1] \ \lxr \ X^{0} \  \lxr \ A[1] \ \lxr \ X^{1}[1]\oplus N
\]
Applying the homological functor $\mathsf{H}$ we have an exact sequence:
\begin{equation}
\begin{multlined}
\mathsf{H}(X^{0}[-1]) \, \lxr \, \mathsf{H}(A) \,  \lxr \, \mathsf{H}(X^{1} \oplus N[-1]) \, \lxr \, \mathsf{H}(X^{0}) \, \lxr \, \\ \, \lxr \, \mathsf{H}(A[1]) \, \lxr \, \mathsf{H}(X^{1}[1] \oplus N)
\end{multlined}
\end{equation}
Since $\mathsf{H}(A)$ is Cohen-Macaulay, it follows that it is a subobject of a projective, i.e. there is a monomorphism $\mathsf{H}(A) \lxr \mathsf{H}(X)$ for some $X \in \X$. Since $\X$ is fully rigid, this monomorphism is induced by a map $A \lxr X$ which when composed with the map $X^{0}[-1] \lxr A$ is zero since $\X$ is rigid. Hence the composition $\mathsf{H}(X^{0}[-1]) \lxr \mathsf{H}(A) \lxr \mathsf{H}(X)$ is zero and therefore   the map $\mathsf{H}(X^{0}[-1]) \lxr \mathsf{H}(A)$ is zero.  Using that the objects $X^{1}[1]$ and $N$ lie in $\X^{\bot}$, we infer that $(5.3)$ becomes an exact sequence:
\begin{equation}
0 \ \lxr \  \mathsf{H}(A)  \ \lxr \ \mathsf{H}(X^{1} \oplus N[-1]) \  \lxr \ \mathsf{H}(X^{0}) \ \lxr \ \mathsf{H}(A[1])  \ \lxr \ 0
\end{equation}
Hence to prove assertion $(5.2)$ it suffices to show that $\mathsf{H}(N[-1])$ is projective. By Proposition $4.16$ we may write $N[-1] = K \oplus L$, where $K \in \X\star\X[1]$ and $L \in \X^{\bot}$ so that $\mathsf{H}(N[-1]) \cong \mathsf{H}(K)$. Let $T^{1} \lxr T^{0} \lxr K \lxr T^{1}[1]$ be a triangle where the $T^{i}$ lie in $\X$. By the $2$-Calabi-Yau property we have $\T(K[1],T^{1}[2]) \cong \mathsf{D}\T(T^{1},K[1])$ and the last space is zero since $K[1]$ lies in $\X^{\bot}$ as a direct summand of $N$. Hence the map $K \lxr T^{1}[1]$ is zero and then $K$ lies in $\X$ as a direct summand of $T^{0}$. We infer that $\mathsf{H}(N[-1]) \cong \mathsf{H}(K)$ is projective in $\smod\X$ and then $(5.4)$ implies that $\Omega^{2}\mathsf{H}(A[1]) \cong \mathsf{H}(A)$ in the stable category $\umod\X$.  Since $\mathsf{H}(A)$ and $\Omega^{1}\mathsf{H}(A[1])$ are Cohen-Macaulay (the latter as a subobject of a projective object), using that $\Omega^{-1}|_{\underline{\mathsf{CM}}(\X)}$ is an equivalence, we have $\Omega^{1}\mathsf{H}(A[1]) \cong \Omega^{-1}\mathsf{H}(A)$ in $\underline{\mathsf{CM}}(\X)$.

\smallskip

\item[{\em Step 3}:] The triangulated category  $\underline{\mathsf{CM}}(\X)$ is $3$-Calabi-Yau. 

\smallskip 

{\em Proof.} Let $\mathsf{H}(A)$ and $\mathsf{H}(B)$ be Cohen-Macaulay objects; clearly we may assume that both objects $A$ and $B$ lie in $\X\star\X[1]$. Since $\smod\X$ is Gorenstein, it follows from \cite{B:cm} that the inclusion $\underline{\mathsf{CM}}(\X) \lxr \umod\X$ admits the functor $\Omega^{-1}\Omega^{1} \colon \umod\X \lxr \underline{\mathsf{CM}}(\X)$ as a right adjoint. Using this and equations $(5.1)$ and $(5.2)$ we have  the following natural isomorphism in $\umod\X$ which shows that $\underline{\mathsf{CM}}(\X)$ is $3$-Calabi-Yau: 
\[
\begin{CD}
\mathsf{D}\uHom\big(\mathsf{H}(A),\mathsf{H}(B)\big)  @> \cong >> \uHom\big(\Omega^{1}\mathsf{H}(B),\mathsf{H}(A[1])\big) @> \cong >>  \uHom\big(\Omega^{1}\mathsf{H}(B),\Omega^{-1}\Omega^{1}\mathsf{H}(A[1])\big)  
\end{CD}
\]
\[
\begin{CD}
@> \cong >> \uHom\big(\Omega^{1}\mathsf{H}(B),\Omega^{-2}\mathsf{H}(A)\big) @> \cong >>  \uHom\big(\mathsf{H}(B),\Omega^{-3}\mathsf{H}(A)\big) 
\end{CD} \qedhere 
\]
\end{enumerate}
\end{proof}  
\end{thm}

\,

\begin{rem} Since $2$-cluster tilting subcategories are fully rigid, Theorem $5.2$ extends a basic result of Keller-Reiten \cite[Proposition $2.1$]{KR1}, see also \cite[Theorem $4.3$]{KZ}, namely that the category of coherent functors over a $2$-cluster tilted subcategory in a ($2$-Calabi-Yau) $k$-linear triangulated category $\T$, with finite-dimensional Hom spaces, is $1$-Gorenstein.   On the other hand, Theorem $5.4$ extends another basic result of Keller-Reiten \cite[Theorem $3.3$]{KR1} which plays an important role in cluster tilting theory, namely  that the stable category of  Cohen-Macaulay modules over  a $2$-cluster tilting subcategory in a $2$-Calabi-Yau triangulated category is $3$-Calabi-Yau. Note that the proof of Theorems $5.2$ and $5.4$ is different from the proof of the $2$-cluster tiling analogues in \cite{KR1} or \cite{KZ}.  For higher cluster tilting theoretic analogues of the above results we refer to \cite{B}.   
\end{rem}

\section{Hereditary Categories and Tilting} 

Let $\T$ be a triangulated category and $\X$ a fully rigid subcategory of $\T$. In this section we study the connections between $2$-cluster tilting subcategories of $\T$ and tilting subcategories of $\smod\X$.

  To proceed further we need the following result from \cite{B}, where $\mathsf{Gh}^{[2]}_{\X}(A,B)$ denotes the subgroup of $\T(A,B)$ consisting of all maps which can be written as compositions of an $\X$-ghost $A \lxr C$ and an $\X[1]$-ghost $C \lxr B$.

\begin{lem} \cite{B} If $\X$ is a contravariantly finite subcategory of $\T$, then $\forall A, B \in \T$ there exists an exact sequence
\begin{multline}
0 \, \lxr \, \mathsf{Gh}_{\X}(A,B) \,  \lxr \, \T(A,B) \, \lxr \, \Hom[\mathsf{H}(A),\mathsf{H}(B)] \, \lxr \, \\ 
\, \lxr \,   \mathsf{Gh}_{\X}(\Omega^{1}_{\X}(A),B)  \, \lxr \, \mathsf{Gh}^{[2]}_{\X}(A,B[1]) \, \lxr \, 0
\end{multline}
\end{lem} 

The following result gives sufficient conditions for a contravariantly finite rigid subcategory $\X$ to be fully rigid or $\smod\X$ to be hereditary. In the sequel we denote by $\T_{\X^{\bot}}$ the full subcategory of $\T$ consisting of all objects with no non-zero direct summands from $\X^{\bot}$. 

\begin{prop} Let $\X$ be a contravariantly finite rigid subcategory of $\T$. 
\begin{enumerate}
\item Let $\X$ be fully rigid.  If $A \in \T$, then $\pd \mathsf{H}(A) \leq 1$ if and only if for any direct summand $A_{1}$ of $A$ lying in $(\X\star\X[1]) \cap \T_{\X^{\bot}}$, any map $A_{1}[-1] \lxr X$ in $\T$, where $X \in \X$, is $\X$-ghost. 

In particular:  $\smod\X$ is hereditary if and only if $\mathsf{Gh}_{\X}(A[-1],\X) = \T(A[-1],\X)$, $\forall A \in (\X\star\X[1]) \cap \T_{\X^{\bot}}$. 
\item If $\smod\X$ is hereditary  and $\mathsf{Gh}_{\X}(\Omega^{i}_{\X}(A),\X) = 0$, $i = 1,2$,  $ \forall A \in \T_{\X^{\bot}}$, then $\X$ is fully rigid. 
\end{enumerate} 
\begin{proof} (i) $(\Longleftarrow)$ Let $A\in \T$. Then by Proposition $4.16$ we may write $A = A_{1} \oplus A_{2}$, where $A_{1} \in \X\star\X[1]$ has no non-zero direct summands from $\X^{\bot}$ and $A_{2} \in \X^{\bot}$. Clearly $\mathsf{H}(A) = \mathsf{H}(A_{1})$. Let $(\dag): \, X^{1} \lxr X^{0} \lxr A_{1} \lxr X^{1}[1]$ be a triangle in $\T$, where the $X^{i}$ lie in $\X$. By hypothesis the map $A_{1}[-1]\lxr X_{1}$ is $\X$-ghost. Hence applying $\mathsf{H}$ to $(\dag)$, we have an exact sequence $0\lxr \mathsf{H}(X^{1}) \lxr \mathsf{H}(X^{0}) \lxr \mathsf{H}(A_{1}) \lxr 0$ which shows that $\pd\mathsf{H}(A)  = \pd\mathsf{H}(A_{1}) \leq 1$.

$(\Longrightarrow)$ Conversely let $\pd\mathsf{H}(A) \leq 1$ and let $A_{1}$ be a direct summand of $A$ with no non-zero direct summands from $\X^{\bot}$. Since $\X$ is fully rigid, this implies that $A_{1}$ lies in $\X\star \X[1]$ and clearly $\pd \mathsf{H}(A_{1}) \leq 1$. Let $\alpha \colon A_{1}[-1] \lxr X$ be a map
 and consider the triangle $(\dag)$ as above. Since the composition $X^{0}[-1] \lxr A_{1}[-1] \lxr X$ is zero, the map $\alpha$ factorizes through the map $A_{1}[-1] \lxr X^{1}$. Hence it suffices to show that the last map is $\X$-ghost. Since $\pd\mathsf{H}(A_{1}) \leq 1$,  $\Image\big(\mathsf{H}(A_{1}[-1]) \lxr \mathsf{H}(X^{1})\big)$ is projective; hence it is isomorphic to an object $\mathsf{H}(X^{*})$ for some $X^{*} \in \X$. Then the projection $\widetilde{\varepsilon} \colon \mathsf{H}(A_{1}[-1]) \lxr \mathsf{H}(X^{*})$ splits, so there exists a map $\widetilde{\mu} \colon \mathsf{H}(X^{*})  \lxr \mathsf{H}(A_{1}[-1])$ such that $\widetilde{\mu} \circ \widetilde{\varepsilon} = 1_{\mathsf{H}(X^{*})}$.  Since $\X$ is fully rigid, we have  $\widetilde{\mu} = \mathsf{H}(\mu)$ and $\widetilde{\varepsilon} = \mathsf{H}(\varepsilon)$ for some maps $\mu \colon X^{*} \lxr A[-1]$ and $\varepsilon \colon A_{1}[-1] \lxr X^{*}$. This clearly implies that $\mu \circ \varepsilon = 1_{X^{*}}$ and therefore $X^{*}$ is a direct summand of $A_{1}[-1]$, or equivalently $X^{*}[1]$ is a direct summand of $A_{1}$. Since $X^{*}[1] \in \X[1] \subseteq \X^{\bot}$, by hypothesis,  $X^{*} = 0$ and therefore the map  $\mathsf{H}(A_{1}[-1]) \lxr \mathsf{H}(X^{1})$ is zero, i.e. the map $A_{1}[-1] \lxr X^{1}$ is $\X$-ghost.

(ii)  Let $\widetilde{\alpha} \colon \mathsf{H}(A) \lxr \mathsf{H}(B)$ be a map in $\smod\X$. Let $\gamma_{A} \colon \Cell_{1}(A) \lxr A$ and $\gamma_{B} \colon \Cell_{1}(B) \lxr B$ be right $\X\star \X[1]$-approximations of $A$ and $B$. Since, by Remark $4.3$, the maps $\mathsf{H}(\gamma_{A})$ and $\mathsf{H}(\gamma_{B})$ are invertible, we have a map $\mathsf{H}(\gamma_{A}) \circ \widetilde{\alpha} \circ \mathsf{H}(\gamma_{B})^{-1} \colon \mathsf{H}(\Cell_{1}(A)) \lxr \mathsf{H}(\Cell_{1}(B))$. Since the functor $\mathsf{H}|_{\X\star\X[1]}$ is full, there is a map $\beta \colon \Cell_{1}(A) \lxr \Cell_{1}(B)$ such that $\mathsf{H}(\gamma_{A}) \circ \widetilde{\alpha} \circ \mathsf{H}(\gamma_{B})^{-1} = \mathsf{H}(\beta)$ and we have a commutative square 
\[
\xymatrix@C=1.5cm{
  \mathsf{H}(\Cell_{1}(A)) \ar[d]_{\mathsf{H}(\gamma_{A}) \,\, }^{\cong} \ar@{-->}[r]^{\mathsf{H}(\beta)}  & \mathsf{H}(\Cell_{1}(B)) \ar[d]^{\, \mathsf{H}(\gamma_{B})}_{\cong}  \\
  \mathsf{H}(A) \ar[r]^{\widetilde{\alpha}} &  \mathsf{H}(B)  }   
\]
  If $\mathsf{H}(\gamma_{A})^{-1} = \mathsf{H}(\rho)$, for some map $\rho \colon A \lxr \Cell_{1}(A)$, then the map $\alpha := \rho \circ  \beta \circ \gamma_{B} \colon A \lxr B$ is such that $\mathsf{H}(\alpha) = \widetilde{\alpha}$. Hence to show that $\mathsf{H}$ is full, it suffices to show that any (invertible) map $\mathsf{H}(A) \lxr \mathsf{H}(B)$, where $B \in \X\star\X[1]$ lies in the image of $\mathsf{H}$. Consider the triangle 
\[
\Omega^{1}_{\X}(A) \ \stackrel{g^{0}_{A}}{\lxr} \ X^{0}_{A} \ \stackrel{g^{0}_{A}}{\lxr} \ A \ \lxr \ \Omega^{1}_{\X}(A)[1]
\]
 and let $X^{1} \stackrel{\lambda}{\lxr} X^{0} \stackrel{\kappa}{\lxr} B \lxr X^{1}[1]$ be a triangle in $\T$, where the $X^{i}$ lie in $\X$. Then we have a commutative diagram
  \begin{equation}
\begin{CD}
\mathsf{H}(\Omega^{1}_{\X}(A)) @> \mathsf{H}(g^{0}_{A}) >> \mathsf{H}(X^{0}_{A}) @> \mathsf{H}(f^{0}_{A}) >> \mathsf{H}(A) @>  >> 0 \\
 @V{\exists \, \widetilde{\delta}}VV  @VV{\mathsf{H}(\sigma)}V   @VV{\widetilde{\alpha}}V & \ &    \\ 
 \mathsf{H}(X^{1})  @> \mathsf{H}(\lambda) >> \mathsf{H}(X^{0}) @> \mathsf{H}(\kappa) >> \mathsf{H}(B) @>  >> 0 
\end{CD}
\end{equation}  
Indeed the existence of the commutative square on the right of $(6.1)$ is clear since the object $\mathsf{H}(X^{0}_{A})$ is projective. Then we have an induced map $\Image\mathsf{H}(g^{0}_{A}) \lxr \Image\mathsf{H}(\lambda)$ which lifts to a map  $\Image\mathsf{H}(g^{0}_{A}) \lxr \mathsf{H}(X^{1})$ since $\Image\mathsf{H}(g^{0}_{A})$ is projective because $\smod\X$ is hereditary. Then the composition  $\widetilde{\delta} := \mathsf{H}(\Omega^{1}_{\X}(A)) \lxr\Image\mathsf{H}(g^{0}_{A}) \lxr \mathsf{H}(X^{1})$ makes the diagram $(6.1)$ commutative. Since $\mathsf{Gh}_{\X}(\Omega^{2}_{\X}(A),\X) = 0$, by Lemma $6.1$ it follows that there exists a map $\delta \colon \Omega^{1}_{\X}(A) \lxr X^{1}$ such that $\mathsf{H}(\delta) = \widetilde{\delta}$. Then the difference $\delta\circ \lambda - g^{0}_{A} \circ \sigma \colon \Omega^{1}_{\X}(A) \lxr X^{0}$ is $\X$-ghost. Since $\mathsf{Gh}_{\X}(\Omega^{1}_{\X}(A),\X) = 0$, we have $\delta\circ \lambda - g^{0}_{A} \circ \sigma = 0$. It follows that we have a morphism of triangles 
\[
\begin{CD}
\Omega^{1}_{\X}(A) @> g^{0}_{A} >> X^{0}_{A} @> f^{0}_{A} >> A @>  >> \Omega^{1}_{\X}(A)[1] \\
 @V{\delta}VV  @VV{\sigma}V   @V{\exists \, \alpha}VV @VV{\delta[1]}V    \\ 
 X^{1}  @> \lambda >> X^{0} @> \kappa >> \ \ \  B \ \ \ @>  >> X^{1}[1]
\end{CD}
\]
and then clearly $\mathsf{H}(\alpha) = \widetilde{\alpha}$. Hence $\mathsf{H}$ is full and therefore $\X$ is fully rigid. 
\end{proof}
\end{prop}

\begin{prop} Let $\T$ be a triangulated category and $\X$ a fully rigid subcategory  of $\T$. Let $A$ be in $\T_{\X^{\bot}}$ and assume that $\pd\mathsf{H}(A) \leq 1$. Then $\forall B \in\T$, there is an isomorphism:
\[
\begin{CD}
\mathsf{Gh}_{\X}(A,B[1]) \,  @> \cong >> \, {\Ext}^{1}\big(\mathsf{H}(A),\mathsf{H}(B) \big)
\end{CD}
\]
\begin{proof} Since $\X$ is fully rigid, as before, $A$ lies in $\X\star\X[1]$ and by Proposition $6.2$, there is a triangle $X^{1} \stackrel{\alpha}{\lxr} X^{0} \lxr A \lxr X^{1}[1]$ in $\T$, where the $X^{i}$ lie in $\X$ and the map $A[-1] \lxr X$ is $\X$-ghost. Then we have a projective resolution $0 \lxr \mathsf{H}(X^{1}) \lxr \mathsf{H}(X^{0}) \lxr \mathsf{H}(A) \lxr 0$ of $\mathsf{H}(A)$  which gives rise to an exact sequence 
\begin{multline}
0 \lxr \Hom\big(\mathsf{H}(A),\mathsf{H}(B)\big) \lxr \Hom\big(\mathsf{H}(X^{0}),\mathsf{H}(B)\big) \lxr \\ \lxr  \Hom\big(\mathsf{H}(X^{1}),\mathsf{H}(B)\big)  \lxr {\Ext}^{1}\big(\mathsf{H}(A),\mathsf{H}(B)\big) \lxr 0
\end{multline}
Since the middle map is isomorphic to the map $\T(\alpha,B) \colon \T(X^{0},B) \lxr \T(X^{1},B)$, it follows that $\Ext^{1}(\mathsf{H}(A),\mathsf{H}(B))$ is isomorphic to $\Coker\T(\alpha,B)$ which, as easily seen, consists of all maps $A[-1] \lxr B$ factorizing through $X^{1}$. Since $\X$ is rigid, any map $A[-1] \lxr B$ factorizing  through an object from $\X$ factorizes through $X^{1}$. It follows that $\Ext^{1}(\mathsf{H}(A),\mathsf{H}(B))$ is isomorphic to the subgroup $G$ of $\T(A[-1],B)$ consisting of all maps $A[-1] \lxr B$ factorizing through an object from $\X$. Since $A$ lies in $\X\star\X[1]$, clearly $G$ is isomorphic to $\mathsf{Gh}_{\X}(A,B[1])$. 
\end{proof} 
\end{prop} 

\begin{prop} Let $\T$ be a triangulated category with a Serre functor $\mathbb S$, and let $\X$  be a covariantly finite fully rigid subcategory  of $\T$.  If $\Y \subseteq \T_{\X^{\bot}}$  satisfies $\mathbb S(\Y) = \Y[2]$, and $\pd\mathsf{H}(Y) < \infty$, $\forall Y \in \Y$, then:
\[\T(\Y,\Y[1]) = 0 \ \ \ \ \ \ \ \textit{if and only if} \ \ \ \ \ \ \ {\Ext}^{1}(\mathsf{H}(\Y),\mathsf{H}(\Y)) = 0\]
In particular an object $A \in \Y$ is rigid in $\T$ if and only if the object $\mathsf{H}(A)$ is rigid in $\smod\X$.
\begin{proof} By Theorem $5.2$ the category $\smod\X$ is $1$-Gorenstein. It follows that any object of $\smod\X$ with finite projective dimension has projective dimension at most one.

``$\Longrightarrow$" If $\T(A,B[1]) = 0$, $\forall A, B \in \Y$, then by Proposition $6.3$ we have $\Ext^{1}(\mathsf{H}(A),\mathsf{H}(B)) = 0$. 

``$\Longleftarrow$" Conversely assume that $\Ext^{1}(\mathsf{H}(A),\mathsf{H}(B)) = 0$, $\forall A,B \in \Y$.   Since  $A, B \in \T_{\X^{\bot}}$ and $\X$ is fully rigid, it follows that $A$ and $B$ lie in $\X\star \X[1]$. This implies that the map $\T(A,B) \lxr \Hom(\mathsf{H}(A),\mathsf{H}(B))$ is surjective. Since $\pd\mathsf{H}(A) < \infty$,  we have $\pd\mathsf{H}(A) \leq 1$ and then $\mathsf{Gh}_{\X}(A,B[1]) = 0$ by Proposition $6.3$.  Therefore $\T(A,B[1]) = \Hom(\mathsf{H}(A), \mathsf{H}(B[1]))$, $\forall A, B \in \Y$. In particular, since $\Y[1] = \mathbb S(\Y)[-1]$,  we have: $\T(A,\mathbb S(B^{\prime})[-1]) = \Hom(\mathsf{H}(A), \mathsf{H}(\mathbb S(B^{\prime})[-1]))$, $\forall B^{\prime} \in \Y$. Hence to show that $\T(A,B[1]) = 0$, it suffices to show that,  $\forall B^{\prime} \in \Y$:
 \[
 \Hom(\mathsf{H}(A), \mathsf{H}(\mathbb S(B^{\prime})[-1])) = 0
 \]
 By Lemma $5.3$ we have $\mathsf{DTr}(\mathsf{H}(B^{\prime})) = \mathsf{H}(\mathbb S(B^{\prime})[-1])$.  
On the other hand, by using Auslander-Reiten duality in $\smod\X$, we have $\overline{\Hom}(\mathsf{H}(A), \mathsf{DTr}\mathsf{H}(B^{\prime})) \cong \mathsf{D}\Ext^{1}(\mathsf{H}(B^{\prime}),\mathsf{H}(A)) = 0$, where $\overline{\Hom}(\mathsf{H}(A),\mathsf{DTr}\mathsf{H}(B^{\prime}))$ is the factor of the space of all maps $\mathsf{H}(A) \lxr \mathsf{H}(B^{\prime}[1])$ by the subspace of maps factorizing through an injective object. Hence any map  $\gamma \colon \mathsf{H}(A) \lxr \mathsf{DTr}\mathsf{H}(B^{\prime})$ factorizes through an injective object. Since $\smod\X$ is a dualizing $k$-variety and $\pd\mathsf{H}(A) \leq 1$, by \cite[Proposition IV.$1.16$]{ARS}, it follows that there are no non-zero maps from an injective object of $\smod\X$ to $\mathsf{DTr}\mathsf{H}(B^{\prime})$. Hence $\Hom(\mathsf{H}(A), \mathsf{DTr}\mathsf{H}(B^{\prime})) = \Hom(\mathsf{H}(A), \mathsf{H}(\mathbb S(B^{\prime})[-1])) = 0$, $\forall B^{\prime}
\in \Y$. Therefore  $\T(A,B[1]) = 0$, $\forall A, B \in \Y$, as required.    
\end{proof}
\end{prop}

A contravariantly finite subcategory $\U$ of an abelian category $\A$ is called a {\em tilting subcategory} if $\Fac(\U) = \{A \in \A \ | \ \Ext^{1}(\U,A) = 0\}$, where $\Fac(\U)$ is the full subcategory of $\A$ consisting of all factors of objects  of $\U$. It is easy to check, see \cite{B:tilting}, that if $\A$ has enough projectives, then a contravariantly finite subcategory $\U$ of $\A$ is  tilting if and only if: $(\alpha)$ $\Ext^{1}(\U,\U) = 0$, $(\beta)$ $\pd U \leq 1$, $\forall U \in \U$, and $(\gamma)$ any projective object $P$ of $\A$  admits an exact sequence $0 \lxr P \lxr U^{0} \lxr U^{1} \lxr 0$, where the $U^{i}$ lie in $\U$. For instance if $\Lambda$ is an Artin algebra, then $T \in \smod\Lambda$ is a tilting module if and only if $\add T$ is a tilting subcategory of $\smod\Lambda$.

\begin{lem} Let $\X$ be a fully rigid subcategory of $\T$ and $\Y \subseteq \T_{\X^{\bot}}$ a contravariantly finite subcategory of $\T$. If $\mathsf{H}(\Y)$ is a tilting subcategory of $\smod\X$, then any  $A\in \T_{\X^{\bot}}$ such that $\Ext^{1}(\mathsf{H}(\Y),\mathsf{H}(A)) = 0 = \T(A,\Y[1])$, lies in $\Y$.    
\begin{proof} Since $\Ext^{1}(\mathsf{H}(\Y),\mathsf{H}(A)) = 0$, it follows that $\mathsf{H}(A)$ is a factor of an object from $\mathsf{H}(\Y)$, i.e. there is an exact sequence $0 \lxr F \lxr \mathsf{H}(Y^{0}) \lxr \mathsf{H}(A) \lxr 0$, where $Y^{0} \in \Y$, and then clearly  $\Ext^{1}(\mathsf{H}(\Y),F) = 0$. For the same reason $F$ is a factor of an object $\mathsf{H}(Y^{1})$, where $Y^{1} \in \Y$.  Hence using that $\mathsf{H}$ is full, we have an exact sequence 
\[
\begin{CD}
\mathsf{H}(Y^{1})  @> \mathsf{H}(f^{1}) >>  \mathsf{H}(Y^{0}) @> \mathsf{H}(f^{0}) >> \mathsf{H}(A) @>  >> 0 
\end{CD}
\]
Since $Y^{1}$ lies in $\Y \subseteq \T_{\X^{\bot}}$, it follows that $Y^{1}$ lies in $\X\star \X[1]$. Hence there exists a triangle  $X^{1} \stackrel{\alpha}{\lxr} X^{0} \stackrel{\beta}{\lxr} Y^{1} \stackrel{\gamma}{\lxr} X^{1}[1]$ in $\T$, where the $X^{i}$ lie in $\X$.  Embedding that map $(f^{1}, \gamma) \colon Y^{1} \lxr Y^{0} \oplus X^{1}[1]$ in a triangle 
\[
\begin{CD}
Y^{1} \, @> >> \, Y^{0} \oplus X^{1}[1] \, @>  >> \, \widetilde{A} \, @> >> \,  Y^{1}[1]
\end{CD}
 \eqno (*)
 \]
it is easy to see that the map $\widetilde{A} \lxr Y^{1}[1]$ is $\X$-ghost since it factorizes through $\beta[1]$.  Hence applying the functor $\mathsf{H}$ to the triangle $(*)$,  we have an exact sequence  $\mathsf{H}(Y^{1}) \stackrel{}{\lxr} \mathsf{H}(Y^{0}) \lxr \mathsf{H}(\widetilde{A}) \lxr 0$. It follows that $\mathsf{H}(A)$ is isomorphic to $\mathsf{H}(\widetilde{A})$. Since $\X$ is fully rigid, this means that $A \oplus  M \cong \widetilde{A} \oplus N$, where $M, N \in \X^{\bot}$. Since $A$ lies in $\T_{\X^{\bot}}$, it follows that $A$ is a direct summand of $\widetilde{A}$: $\widetilde{A} \cong A \oplus L$, where clearly $L \in \X^{\bot}$. Then the triangle $(*)$ can be written as: 
\[
\begin{CD}
Y^{1} \, @> >> \, Y^{0} \oplus X^{1}[1] \, @> >>  \,  A \oplus L \, @> >> \, Y^{1}[1]
\end{CD}
\] 
Now since  $\T(A,\Y[1]) = 0$, the composition of the split  inclusion $A \lxr A \oplus L$ with the map $A \oplus L \lxr Y^{1}[1]$ is zero and therefore $A \lxr A \oplus L$ factorizes through $Y^{0} \oplus X^{1}[1]$. This clearly implies that $A$ is a direct summand of $Y^{0} \oplus X^{1}[1]$, and since $A \in \T_{\X^{\bot}}$, it follows that $A$ is a direct summand of $Y^{0}$,  hence $A \in \Y$. 
\end{proof}
\end{lem}

\begin{nota} Let $\T$ be a Krull-Schmidt category and $\X$ contravariantly finite fully rigid subcategory of $\T$. If  $\U$ is  a subcategory of $\smod\X$, we denote by $\mathsf{H}^{-1}(\U)$ the full subcategory $\Y = \{Y \in \T_{\X^{\bot}} \  | \ \mathsf{H}(Y) \in \U\}$. 

 Note that $\Y$ is unique. Indeed if $\Z \subseteq \T_{\X^{\bot}}$ is such that $\mathsf{H}(\Y) = \mathsf{H}(\Z)$, then for any object  $Y\in \Y$,  we have an isomorphism $\mathsf{H}(Y) \cong \mathsf{H}(Z)$ for some $Z \in \Z$.   Since  $\smod\X \approx \T/\X^{\bot}$, it follows that $Y \oplus M \cong Z\oplus N$ for some objects $M,N \in \X^{\bot}$.  Since $\T$  is Krull-Schmidt and  $\Y$ has no non-zero direct summand from $\X^{\bot}$ it follows that $Y$ is a direct summand of $Z$, hence $Y \in \Z$ and therefore $\Y \subseteq \Z$. Dually we have $\Z \subseteq \Y$.  
\end{nota}

Now we can prove the main result of this section which gives a connection between $2$-cluster tilting subcategories of $\T$ and tilting subcategories of $\smod\X$, for a contravariantly finite subcategory $\X$ of $\T$.  

\begin{thm} Let $\T$ be a Hom-finite triangulated category over a field $k$ which admits a Serre functor $\mathbb S$, and let $\X$ be a contravariantly finite subcategory of $\T$. If $\gd\smod\X < \infty$, then the following are equivalent:
\begin{enumerate}
\item $\X$ is $2$-cluster tilting subcategory of $\T$.
\item  $\X$ is fully rigid, $\mathbb S(\X) = \X[2]$, and the map $\Phi \colon \Y \ \longmapsto \ \mathsf{H}(\Y)$ gives a bijective correspondence between:
\begin{enumerate}
\item[$\mathsf{(I)}$] $2$-cluster tilting subcategories $\Y$ of $\T$ contained in $\T_{\X^{\bot}}$. 
\item[$\mathsf{(II)}$] Tilting subcategories $\U$ of $\smod\X$ such that:
\begin{enumerate}
\item[$\mathrm{(a)}$] $\mathsf{H}^{-1}(\U)$ is contravariantly finite in $\T$, and 
\item[$\mathrm{(b)}$] $\mathbb S\mathsf{H}^{-1}(\U) = \mathsf{H}^{-1}(\U)[2]$. 
\end{enumerate}
\end{enumerate}
\end{enumerate}  
Under the above correspondence we have: $\U \approx \Y/\mathsf{Gh}_{\X}(\Y)$.   Moreover if $\Y = \add T$, so $\U =  \add \mathsf{H}(T)$, then there is an algebra isomorphism: 
\[
{\End}_{\T}(T)/\mathsf{Gh}_{\X}(T) \ \cong \ {\End}_{\smod\X}(\mathsf{H}(T))
\] where $\mathsf{Gh}_{\X}(T)$ is the ideal of $\End_{\T}(T)$ consisting of all maps $T \lxr T$ factorizing through an object from $\X[1]$. 
\begin{proof} (ii) $\Longrightarrow$ (i) Assume that $\X$ is fully rigid, $\mathbb S(\X) = \X[2]$, and the map $\Phi$ in the statement is a bijection. Then the tilting subcategory $\U = \proj\smod\X$, satisfies the assumptions of $\mathsf{(II)}$ and we have $\mathsf{H}^{-1}(\proj\smod\X) = \X$.  Therefore there exists a $2$-cluster tilting subcategory $\Y$ of $\T$ contained in $\T_{\X^{\bot}}$ such that $\mathsf{H}(\Y) = \mathsf{H}(\X)$. Since $\Y$, as discussed above, is unique, we have $\Y = \X$, so $\X$ is $2$-cluster tilting.

(i) $\Longrightarrow$ (ii) First note that any $2$-cluster tilting subcategory $\M$ of $\T$ satisfies: $\mathbb S(\M) = \M[2]$. Indeed we clearly have $\M^{\bot} = \M[1]$ and ${^{\bot}}\M = \M[-1]$, and $\T(\mathbb S^{-1}(\M[2]),\M[1]) = \T(\M[2],\mathbb S(\M[1])) =\mathsf{D}\T(\M[1],\M[2]) = \mathsf{D}\T(\M,\M[1]) = 0$, hence $\mathbb S^{-1}(\M[2]) \in {^{\bot}}\M[1] = \M$. It follows that $\M[2] \subseteq \mathbb S(\M)$. Dually we have $\mathbb S(\M) \subseteq \M[2]$ and therefore $\mathbb S(\M) = \M[2]$. Now we fix a $2$-cluster tilting subcategory $\X$ of $\T$. Note that since $\gd\smod\X < \infty$, it follows by Theorem $5.2$ that $\smod\X$ is hereditary. 

$\mathsf{(I)} \lxr \mathsf{(II)}$ Let $\Y$ be a $2$-cluster tilting subcategory of $\T$ contained in $\T_{\X^{\bot}}$. Then $\Y$ is rigid, hence by Proposition $6.3$ we have $\Ext^{1}(\U,\U) = 0$, and $\T = \Y \star \Y[1]$. In particular for any object $X[1]$, where $X \in \X$, there exists a triangle $Y^{1} \lxr X[1] \lxr Y^{0}[1] \lxr Y^{1}[1]$, or equivalently a triangle $X\lxr Y^{0} \lxr Y^{1} \lxr X[1]$, in $\T$ where the $Y^{i}$ lie in $\Y$. Since the $Y^{i}$ lie in $\T_{\X^{\bot}}$ and $\pd\mathsf{H}(Y^{i}) \leq 1$, by Proposition $6.2$ it follows that the map $Y^{1}[-1]\lxr X$ is $\X$-ghost. Therefore we have an exact sequence $0 \lxr \mathsf{H}(X) \lxr \mathsf{H}(Y^{0}) \lxr \mathsf{H}(Y^{1}) \lxr 0$ in $\smod\X$, where the $\mathsf{H}(Y^{i})$ lie in $\U$. Since clearly $\U$ is contravariantly finite in $\smod\X$ and since $\pd U\leq 1$, $\forall U \in \U$, we infer that $\U = \mathsf{H}(\Y)$ is a tilting subcategory of $\smod\X$ and $\mathsf{H}^{-1}(\U) = \Y$. By the above argument the subcategory  $\mathsf{H}^{-1}(\U)$ is contravariantly finite in $\T$  and satisfies satisfies $\mathbb S\mathsf{H}^{-1}(\U) = \mathsf{H}^{-1}(\U)[2]$.

$\mathsf{(II)} \lxr \mathsf{(I)}$ Let $\U$ be a tilting subcategory of $\smod\X$ such that  $\Y = \mathsf{H}^{-1}(\U)$ is contravariantly finite in $\T$ and $\mathbb S(\Y) = \Y[2]$. We have to show that $\Y$ is $2$-cluster tilting. Since $\U = \mathsf{H}(\Y)$ is tilting,  we have $\Ext^{1}(\mathsf{H}(\Y),\mathsf{H}(\Y)) = 0$ and $\pd\mathsf{H}(Y) \leq 1$, $\forall Y \in \Y$. Then by Proposition $6.4$, it follows that $\T(Y_{1},Y_{2}[1]) = 0$. Hence $\Y$ is rigid. Let $A$ be an indecomposable object in $\T$ such that $\T(\Y,A[1]) = 0$. If $A$ lies in $\X^{\bot}$, then since $\X$ is $2$-cluster tilting, we have $\X^{\bot}  = \X[1]$ and therefore $A[1] = X^{\prime}[2]$, for some $X^{\prime} \in \X$. Since $\mathbb S(\X) = \X[2]$ we further have $A[1] = \mathbb S(X^{\prime\prime})$, for some $X^{\prime\prime} \in \X$. Then $0 = \T(\Y,A[1]) = \T(\Y,X^{\prime}[2])  = \T(\Y,\mathbb S(X^{\prime\prime}))\cong \mathsf{D}\T(X^{\prime\prime},\Y)$. Since $\U = \mathsf{H}(\Y)$ is tilting, the projective object $\mathsf{H}(X^{\prime\prime})$ is a subobject of an object $\mathsf{H}(Y)$, where $Y \in \Y$. Hence $\T(X^{\prime\prime},\Y) \cong \Hom(\mathsf{H}(X^{\prime\prime}),\mathsf{H}(\Y)) \neq 0$ and this is a contradiction. We infer that $A \notin \X^{\bot}$ and therefore $A \in \T_{\X^{\bot}}$. Since $\T(\Y,A[1]) = 0$, we have $\Ext^{1}(\mathsf{H}(\Y),\mathsf{H}(A)) = 0$. On the other hand  by Serre duality we have $0  = \mathsf{D}\T(\Y,A[1]) \cong \T(A[1],\mathbb S(\Y))$. Since $\mathbb S(\Y) = \Y[2]$, we have $ 0 = \T(A[1],\mathbb \Y[2]) = \T(A,\Y[1])$. Then Lemma $6.5$ applies and gives that  $A \in \Y$. We infer that $\Y$ is a $2$-cluster tilting subcategory of $\T$. Clearly the map $\Phi$ is a bijection.
\end{proof}
\end{thm} 

\begin{rem} (i) The map $\mathsf{(II)} \longmapsto \mathsf{(I)}$ in Theorem $6.6$ does not uses that $\gd\smod\X < \infty$. It follows that any tilting subcategory $\U$ of $\smod\X$ satisfying the conditions in $\mathsf{(II)}$ of Theorem $6.6$, for instance if $\T$ is $2$-Calabi-Yau and $\U = \add M$ for some $M \in \smod\X$,  lifts to a $2$-cluster tilting subcategory $\Y = \mathsf{H}^{-1}(\U)$ of $\T$. 

(ii) In connection with part $\mathsf{(I)}$ of Theorem $6.6$, it is easy to see that a $2$-cluster tilting subcategory $\Y$ of $\T$ is contained in $\X^{\bot}$ if and only if $\Y = \X[1]$.  
\end{rem}

We point out some consequences of the above result. 

\begin{cor} Let $\T$ be a Hom-finite triangulated category with Serre duality over a field $k$. Let $\X$ be a $2$-cluster tilting subcategory of $\T$ and assume that $\gd\smod\X < \infty$. If $\Y$ and $\Z$ are  $2$-cluster tilting subcategories of $\T$ contained in $\T_{\X^{\bot}}$, then there are triangle equivalences:
\[
\begin{CD}
{\bf D}^{b}(\smod\mathsf{H}(\Y))   @< \approx <<  {\bf D}^{b}(\smod\X)    @> \approx >> {\bf D}^{b}(\smod\mathsf{H}(\Z))
\end{CD}
\]
\begin{proof} By Theorem $6.6$, $\mathsf{H}(\Y)$ and $\mathsf{H}(\Z)$ are tilting subcategories of $\smod\X$. Then the assertion follows from a well-known result of Happel \cite{Happel}.  
\end{proof}
\end{cor} 

In \cite[Theorem I.$1.8$]{BIRS} it is shown that if $\T$ is an algebraic $2$-Calabi-Yau triangulated category with a $2$-cluster tilting object $T$, then any rigid object is a direct summand of a $2$-cluster tilting object. The next result proves the same assertion for, not necessarily algebraic, $2$-Calabi-Yau triangulated categories having a $2$-cluster tilting object whose endomorphism ring has finite global dimension. 

\begin{cor} Let $\T$ be a $2$-Calabi-Yau triangulated category over a field $k$, and $T \in \T$ a $2$-cluster tilting object such that $\gd\End_{\T}(T) < \infty$. Then any rigid object of $\T$ is a direct summand of a $2$-cluster tilting object of $\T$. 
\begin{proof} Let $S$ be a rigid object.  By Proposition $6.3$, the object $\mathsf{H}(S)$ of $\smod\X$ is rigid in $\smod\End_{\T}(T)$.  Since $\gd\End_{\T}(T) < \infty$,  it follows that from Theorem $5.2$ that $\End_{\T}(T)$ is hereditary and therefore $\pd\mathsf{H}(S) \leq 1$. By Bongartz's Lemma, see \cite{Bongartz}, there exists a $\End_{\T}(T)$-module $F$ such that $\mathsf{H}(S) \oplus F$ is a tilting $\End_{\T}(T)$-module.   Let $K$ be an object in $\T_{\T^{\bot}}$ such that $\mathsf{H}(K) \cong F$. Then $\mathsf{H}(S\oplus K)$ is a tilting module. Since $\T$ is $2$-Calabi-Yau and $\add(S\oplus K)$ is clearly contravariantly finite in $\T$, by Theorem $6.6$ the object $S \oplus K$ is $2$-cluster tilting. 
\end{proof}
\end{cor}

\begin{rem} The above results generalize and improve several results in the literature. More precisely  Proposition $6.3$ and Theorem $6.6$ generalize results of J{\o}rgensen-Holm, see \cite[Proposition $1.5$, Theorem $3.4$, Proposition $3.6$]{JH}, where related results were proved in case $\X$ is $2$-cluster tilting and any object of $\smod\X$ has finite length. Proposition $6.4$ generalizes a result of Iyama-Yoshino \cite[Corollary $6.5$]{IY}, where an analogous result was proved in case $\T$ is $2$-Calabi-Yau and $\X$ is $2$-cluster tilting.  The lifting of any tilting subcategory $\U$ of $\smod\X$ to a $2$-cluster tilting subcategory of $\T$ in Theorem $6.6$ was proved in \cite[Theorem $1$]{Smith} in case $\T$ is the cluster category of a finite-dimensional hereditary algebra and $\U  = \add M$,  and in \cite[Theorem $3.3$]{FuLiu} in case $\T$ is $2$-Calabi-Yau and $\U = \add M$. Finally we refer to \cite[Proposition $2.26$]{IT} for related results in other contexts. 
\end{rem}

\section{The Free Abelian Category of a Rigid Subcategory}  

Let as before $\T$ be a triangulated category and $\X$ a functorially finite rigid subcategory of $\T$. In this section we are interested in the question of when the stable pretriangulated category $\T/\X$ is abelian. 

By the previous results, associated to $\X$ are  the following subcategories of the pretriangulated category $\T/\X$:
\[
\mathsf{A}(\X) := \frac{\X[-1]\star \X}{\X} \,\,\,\,\ \subseteq \,\,\,\,\ \frac{\T}{\X} \,\,\,\,\ \supseteq \,\,\,\,\  \frac{\X^{\bot}[-1]}{\X}  := \mathsf{T}(\X) 
 \] 
 Note that:
 \begin{itemize}
 \item $\mathsf{A}(\X)$ is abelian since it is equivalent to $\smod\X$,  and
 \item $\mathsf{T}(\X)$ is triangulated provided that $\T$ admits a Serre functor $\mathbb S$ and $\mathbb S(\X) = \X[2]$, see  Theorem $3.15$. 
 
 In this case we have: $\mathsf{T}(\X) = \X^{\bot}[-1]/\X = {^{\bot}}\X[1]/\X$. 
 \end{itemize}
 
 \,

Since $\X$ is rigid, the functor $\T \lxr \smod\X$, $A \longmapsto \T(-,A[1])|_{\X}$ kills the objects of $\X$, so it induces a  functor  
\[
\begin{CD}
\mathbb{R} \ \colon \ \T/\X \ @> >> \ \smod\X, \ \ \ \ \ \mathbb{R}(\underline{A}) = \T(-,A[1])|_{\X}
\end{CD}
\] 

\, 

\begin{lem} The functor $\mathbb R$ admits a fully faithful left adjoint  \ $\mathbb L \colon \smod\X \lxr \T/\X$,  and \[\Ker\mathbb R = \mathsf{T}(\X)\]
\begin{proof} We define an additive functor $\mathbb{L} \colon \smod\X \lxr \T/\X$ as the composition
\[
\begin{CD}
\smod\X  \ @> \approx >>  \  (\X\star\X[1])/\X[1] \  @> \approx >> \ (\X[-1]\star\X)/\X \ @> >>  \ \T/\X
\end{CD}
\]   
where the first functor is the quasi-inverse of the equivalence $(\X\star\X[1])/\X[1] \lxr \smod\X$, see Lemma $4.4$, the second one is given by sending $A$ to $A[-1]$ and the third one is the canonical inclusion. It is easy to see that if $F$ is in $\smod\X$ and $\T(-,X^{1}) \lxr \T(-,X^{0}) \stackrel{\varepsilon}{\lxr} F \lxr 0$ is a projective presentation of $F$ induced by a map $f \colon X^{1} \lxr X^{0}$ between objects of $\X$, then $\mathbb{L}(F) = \underline{C(f)[-1]}$, where $C(f)$ is defined by the triangle $X^{1} \stackrel{f}{\lxr} X^{0} \stackrel{g}{\lxr} C(f) \lxr X^{1}[1]$ in $\T$. Since $\X$ is rigid, we have $\T(-,C(f))|_{\X} = F$ and we may define a map  
\[
\begin{CD}
\varphi_{F,\unA} \,\ \colon \,\,  \Hom\big(\underline{C(f)[-1]},\unA\big) \, @> >> \, \Hom\big( F,\mathbb \T(-,A[1])|_{\X}), \,\,\,\ \ \ \varphi_{F,\unA}(\underline{\alpha}) = \T(-,\alpha[1])|_{\X}
\end{CD}
\] 
If $\varphi_{F,\unA}(\underline{\alpha}) = \T(-,\alpha[1])|_{\X} = 0$, then $\T(X^{0},\alpha[1]) = 0$ and therefore $\T(X^{0},\alpha[1])(g) = g \circ \alpha[1] = 0$. Hence $g[-1] \circ \alpha = 0$ and this implies that $\alpha \colon C(f)[-1] \lxr A$ factorizes through $X^{1}$. It follows that $\underline{\alpha} = 0$, hence $\phi_{F,\unA}$ is injective. Let $\widetilde{\alpha} \colon   F \lxr \mathbb \T(-,A[1])|_{\X}$ be a map in $\smod\X$. By Yoneda's Lemma the composition $\varepsilon \circ \widetilde{\alpha}$ induces a map $\gamma \colon X^{0} \lxr A[1]$ and then clearly $f \circ \gamma = 0$. Hence $\gamma$ factorizes through $g$ and this gives a map $\alpha \colon C(f) \lxr A[1]$ such that $g \circ \alpha = \gamma$. Plainly $\varphi_{F,\unA}(\underline{\alpha[-1]}) = \T(-,\alpha)|_{\X} = \widetilde{\alpha}$ and this shows that 
$\varphi_{F,\unA}$ is surjective, hence an isomorphism which is easily seen to be natural.
 We infer that $(\mathbb R,\mathbb L)$ is an adjoint pair.   By construction $\mathbb L$ is fully faithful and  $\Ker\mathbb R = \big\{\unA \in \T/\X \,\, | \,\, \mathbb R(\unA) = \T(-,A[1])|_{\X} = 0\big\} = (\X^{\bot}[-1])/\X = \mathsf{T}(\X)$. 
\end{proof} 
\end{lem}

\begin{lem} If the factor category $\T/\X$ is abelian, then  $\mathsf{T}(\X)$ is a colocalizing abelian subcategory of $\T/\X$,  the functor $\mathbb R \colon \T/\X \lxr \smod\X$ is exact and induces an exact sequence of abelian categories
\begin{equation}
\begin{CD}
0 \, @> >> \,  \mathsf{T}(\X) \, @> >> \,  \T/\X \, @> >> \,  \mathsf{A}(\X) \, @> >> \,  0
\end{CD}
\end{equation}
If $\T$ admits a Serre functor $\mathbb S$ and $\mathbb S(\X) = \X[2]$, then the triangulated category $\mathsf{T}(\X)$ is semisimple abelian. Moreover  $\X$ is $2$-cluster tilting if and only if the functor $\mathbb R \colon \T/\X \lxr \mathsf{A}(\X)$ is an equivalence. 
\begin{proof} Let $\unf \colon \unB \lxr \unC$ be an epimorphism in $\T/\X$. Let $B \stackrel{f}{\lxr} C \stackrel{g}{\lxr} A \stackrel{h}{\lxr} B[1]$ be a triangle in $\T$. Since  $\unf \circ \ung = 0$, then map $\ung$ is zero, i.e. $g$ factorizes through an object form $\X$. Then the map $C[1] \lxr A[1]$ factorizes through  an object from $\X[1]$ and therefore, since $\X$ is rigid, from the exact sequence  $\T(-,B[1])|_{\X} \lxr \T(-,C[1])|_{\X} \lxr \T(-,A[1])|_{\X}$ we infer that the map $\T(-,B[1])|_{\X} \lxr \T(-,C[1])|_{\X}$ is an epimorphism in $\smod\X$. In other words, the map $\mathbb R(f) \colon \mathbb R(\unB) \lxr \mathbb R(\unC)$ is an epimorphism, i.e. $\mathbb R$ is right exact. Then $\mathbb R$ is exact since it is a right adjoint of $\mathbb L$. It follows that $\Ker\mathbb R$ is a colocalizing subcategory of $\T/\X$ and we have  the exact sequence of abelian categories $(7.1)$. If $\mathbb S(\X) = \X[2]$ for a Serre functor $\mathbb S$ of $\T$,  then  $\mathsf{T}(\X)$ is triangulated by Proposition $3.13$. Since an abelian triangulated category is semisimple, we infer that $\mathsf{T}(\X)$ is semisimple abelian.  Now $\mathbb R$ is an equivalence if and only if $\Ker\mathbb R = 0$ if and only if  $\mathsf{T}(\X) = 0$ if and only if $\X^{\bot}[-1] = \X$. This is equivalent to  $\X = \big\{A \in \T \, | \, \T(\X,A[1]) = 0 \big\}$, i.e. $\X$ is $2$-cluster tilting. 
\end{proof}
\end{lem} 

Now we can prove the following characterization of $2$-cluster tilting subcategories which complements, and is inspired by, related results of Koenig-Zhu, see \cite[Theorems 5.1, 5.2]{KZ}. 

\begin{thm} Let $\T$ be a connected triangulated category with a Serre functor $\mathbb S$ and  $\X$ a non-zero  functorially finite rigid subcategory of $\T$.  Then the following are equivalent:
\begin{enumerate} 
\item $\X$ is $2$-cluster tilting.
\item $\X$ is a maximal extension closed subcategory of $\T$ such that $\mathbb S(\X) = \X[2]$. 
\item $\T/\X$ is abelian \, and \, $\mathbb S(\X) = \X[2]$. 
\end{enumerate}
\begin{proof} (i) $\Longleftrightarrow$ (ii) $\Longrightarrow$ (iii) If $\X$ is $2$-cluster tilting, then $\T = \X \star \X[1]$, $\X^{\bot} = \X[1]$ and ${^{\bot}}\X = \X[-1]$. It follows from Lemma $4.4$ that $\T/\X$ is abelian, and as observed in the proof of Theorem $6.6$, we have $\mathbb S(\X) = \X[2]$. Hence (i) implies (iii) and  the  equivalence between (i) and (ii) follows from remark $3.16$. 

(iii) $\Longrightarrow$ (i) Since $\mathbb S(\X) = \X[2]$ it follows that $\X^{\bot}[-1] = {^{\bot}}\X[1]$ and we have to show that $\mathsf{T}(\X) = 0$, i.e. $\U := \X^{\bot}[-1] = \X$. Let $A$ be an indecomposable object in $\U$ and $A \notin \X$.  By \cite{VdBR}  there exists an Auslander-Reiten triangle $\mathbb S A[-1] \stackrel{g}{\lxr}   B \stackrel{f}{\lxr} A \stackrel{h}{\lxr} \mathbb S A$ in $\T$. In particular $h \neq 0$ and any map $\alpha \colon C \lxr A$ which is not split epic, factorizes through $f$.  Since $A \notin \X$, no map $X \lxr A$ where $X \in \X$, is split epic and therefore $h$ is $\X$-ghost, i.e. $\T(\X,h) = 0$.  It follows that $f^{0}_{A}$ factorizes through $f$ and there exists a morphism of triangles
\begin{equation}
\begin{CD}
\Omega^{1}_{\X}(A) @> g^{0}_{A} >> X^{0}_{A} @> f^{0}_{A} >> A @> h^{0}_{A} >> \Omega^{1}_{\X}(A)[1] & \ & \ & \ & \ & \ & \ & \ & \ & \ & \  (T^{\prime})\\
 @V{\delta}VV  @VV{\sigma}V   \Big\| & \ &  @VV{\delta[1]}V   \\ 
 \mathbb S(A)[-1]  @> g >> B @> f >> \ \ \ A \ \ \ @> h >> \mathbb S(A) & \ & \ & \ & \ & \ & \ & \ & \ & \ & \  (T)
\end{CD}
\end{equation} 
which by subsection 2.3 induces the following left triangle in the pretriangulated category $\T/\X$: 
\[
\begin{CD}
\Omega^{1}_{\X}(A) \, @> \underline{\delta} >> \, \mathbb S(A)[-1] \, @> \underline{g} >> \, B \, @> \underline{f} >> \, A 
\end{CD}
\eqno (T^{\prime\prime})
\] 
Since $\mathbb S(\X) = \X[2]$, it follows from Proposition $3.13$(ii) that $\mathbb S(\U) = \U[2]$ and therefore 
$\mathbb S(A)[-2]\in \X^{\bot}[-1] = \U$ and $\mathbb S^{-1}(A)[2] \in \U$. We shall show in four steps that the assumption $A \notin \X$ leads to a contradiction. Note that, by Lemma $7.2$, $\mathsf{T}(\X)$ is a Serre subcategory of $\T/\X$, so it is closed under subobjects and quotients, and also it is semisimple abelian, so a non-zero object in $\mathsf{T}(\X)$ is indecomposable if and only if it is  simple.
\begin{enumerate}
\item[{\em Step 1}:] The middle term $B$ of the Auslander-Reiten triangle $(T)$ lies in $\X$. 

  {\em Proof:} Let  $\unf \colon \unB \stackrel{\underline{l}}{\lxr} \unC \stackrel{\underline{k}}{\lxr} \unA$ be an epic-monic  factorization of $\unf$ in $\T/\X$. Since $\unA$ is a simple object  in $\U/\X$,  it follows that either $(\alpha)$ $\unC = \unA$,  or $(\beta)$ $\unC = 0$, i.e. $C \in \X$. 
  
 $(\alpha)$ If $\unC = \unA$, then $\unf \colon \unB \lxr \unA$ is an epimorphism. Let $\underline{\kappa} \colon \unK \lxr \unB$ be the kernel of $\unf$ in $\T/\X$. Since $g \circ f = 0$, there exists a unique map $\underline{\lambda} \colon \underline{\mathbb S(A)[-1]} \lxr \unK$ such that $\underline{\lambda} \circ \underline{\kappa} = \ung$ and since, from the triangle $(T^{\prime\prime})$, $\ung$ is a weak kernel of $\unf$, there exists a map $\underline{\mu} \colon \unK \lxr \underline{\mathbb S A[-1]}$ such that $\underline{\mu} \circ \underline{g} = \underline{\kappa}$. It follows from this that $1_{\underline{K}} =  \underline{\mu} \circ \underline{\lambda}$, hence we have a direct sum decomposition $\underline{\mathbb S A[-1]} \cong \unK \oplus \underline{L}$ in $\T/\X$. Since $\mathbb S A[-1]$ is indecomposable in $\T$, so in $\underline{\mathbb S A[-1]}$ in $\T/\X$, hence either $\unK = 0$, i.e. $K \in \X$, or else the map $\underline{\lambda} : \underline{\mathbb S A[-1]} \lxr \unK$ is invertible. In the latter case,  $\ung$ is monic and this implies that the map $\delta \colon \Omega^{1}_{\X}(A) \lxr \mathbb S A[-1]$ factorizes through $\X$. Since $\T(\X,\mathbb S A[-1]) = 0$, it follows that $\delta = 0$ and therefore from (7.2) we have $h = 0$ and this is impossible; hence  $K \in \X$. Then $\unf$ monic and therefore it is invertible since $\T/\X$ is abelian. Then also $\Omega^{1}_{\X}(\unf)$ is invertible and therefore from the triangle $(T^{\prime\prime})$, we deduce that $\underline{\delta} = 0$. Then as above $h = 0$ and this is not the case.

$(\beta)$ We infer that $C := X^{\prime} \in \X$; then since $h$ is $\X$-ghost, the composition $k \circ h$ is zero and therefore there exists a map $m \colon X^{\prime} \lxr B$ such that $m \circ f = k$. This implies that the map $1_{B} - l \circ m \colon B \lxr B$ admits a factorization   $1_{B} - l \circ m = r \circ g$ and then $1_{\unB} = \underline{r} \circ \underline{g}$, i.e. $\underline{g}$ is split epic. Then the triangle $(T^{\prime\prime})$ gives a direct sum decomposition $\underline{\mathbb S A[-1]} \cong \unB \oplus \Omega^{1}_{\X}(A)$ in $\T/\X$. Since $\mathbb S A[-1]$ is indecomposable, as in $(\alpha)$ we have $\unB = 0$ or else $\Omega^{1}_{\X}(\unA) = 0$, i.e. either $B \in \X$ or else $\Omega^{1}_{\X}(A) \in \X$. If $\Omega^{1}_{\X}(A)\in \X$, then the map $A \lxr \Omega^{1}_{\X}(A)[1]$ lies in $\T(A,\X[1])$ which is zero since $A \in \U$. Then $A$ lies in $\X$ as a direct summand of $X^{0}_{A}$ and this is not the case. We infer that $B := X \in \X$. 

\smallskip

\item[{\em Step 2}:] Setting $g^{*} = \underline{g[-1]}$, $f^{*} = \underline{f[-1]}$, and  $h^{*} = \underline{h[-1]}$,  we have an exact sequence in $\T/\X$:
\begin{equation}
\begin{CD}
0 \, @> >> \,  \underline{\mathbb S A[-2]} \, @> g^{*} >> \,  \underline{X[-1]} \, @> f^{*} >>  \,  \underline{A[-1]} \, @> h^{*} >> \,   \underline{\mathbb S A[-1]} \, @> >> \,  0 
\end{CD}
\end{equation}

{\em Proof:} From Step 1 it follows that the triangles $(T^{\prime})$ and $(T)$ in (7.2) are isomorphic, in particular $\mathbb S A[-1] = \Omega^{1}_{\X}(A) \in \U$. Since always $\Omega^{1}_{\X}(A)[1]$ lies $\X^{\bot}$, we infer that $\T(\X,\mathbb S A) = 0$.  Then by Serre duality we have $0 = \T(\X,\mathbb S A) = \mathsf{D}\T(A,\X)$, hence  $\T(A,\X) = 0$. Since $\T(\mathbb S A[-2],\X) = \T(\mathbb S A, \X[2]) = \T(\mathbb S A, \mathbb S X) = \T(A,\X) = 0 = \T(A,X[1])$, the maps $g^{*}$ and $h^{*}$ are non zero in $\T/\X$.  Since $\underline{\mathbb S A[-2]}$ and $\underline{\mathbb S A[-1]}$ are simple objects in $\mathsf{T}(\X)$, it follows that the map $g^{*}$ is  monic and the map $h^{*}$ is epic. Let $\underline{\beta} \colon \unB \lxr \underline{X[-1]}$ be a map such that $\underline{\beta} \circ f^{*} = 0$. Then the map $\beta \circ f[-1]$ factorizes through an object $X^{\prime} \in \X$, i.e.  
$\beta \circ f[-1] = \kappa \circ \lambda$, where $\kappa \colon B \lxr X^{\prime}$ and $\lambda \colon X^{\prime} \lxr A[-1]$. Since $\T(\X,A[1]) = 0$, the composition $\lambda \circ h[-1]$ is zero, hence there is a map $\rho \colon X^{\prime} \lxr X[-1]$ such that $\rho \circ f[-1] = \lambda$.  It follows that $\beta \circ f[-1] = \kappa \circ \rho \circ f[-1]$, hence $\beta - \kappa \circ \rho = \sigma \circ g[-1]$ for some map $\sigma \colon B \lxr A$.  Then in $\T/\X$ we have $\underline{\beta} = \underline{\sigma} \circ g^{*}$. This shows that $g^{*} = \mathsf{ker} f^{*}$. Similarly let $\underline{\alpha} \colon \underline{A[-1]} \lxr \unC$ be a map such that $f^{*} \circ  \underline{\alpha} = 0$, hence $f[-1] \circ \alpha$ factorizes through an object $X^{\prime\prime} \in \X$: $f[-1] \circ  \alpha = \mu \circ \nu$, where $\mu \colon X[-1] \lxr X^{\prime\prime}$ and $\nu \colon X^{\prime\prime} \lxr C$ . Since $\X$ is rigid we have $\mu = 0$ and therefore $f[-1] \circ  \alpha = 0$. Hence $\alpha$ factorizes through $h[-1]$ and then $\underline{\alpha}$ factorizes through $h^{*}$. This shows that $\mathsf{coker}f^{*} = h^{*}$. 

\smallskip

\item[{\em Step 3}:] We have $B := X = 0$ and the Auslander-Reiten triangle starting at $A$ is of the form:
\begin{equation}
\begin{CD}
A[-1] \, @>  >> \, 0 \, @> >> \, A \, @> 1_{A} >>  \, A
\end{CD}
\end{equation}

{\em Proof:} Assume that $X \neq 0$.  Let $\underline{E} = \Image(f^{*})$ be the image of $f^{*}$ in $\T/\X$, and let
\[
\begin{CD}
f^{*} = \underline{e} \circ \underline{m} \,  \colon \,  \underline{X[-1]} \, @> \underline{e} >> \, \underline{E} \, @> \underline{m} >> \,  \underline{A[-1]}
\end{CD}
\]
be its canonical factorization. Since $g^{*} \circ \underline{e} = 0$, the composition $g[-1] \circ e \colon \mathbb S A[-2] \lxr E$ factorizes through an object from  $\X$. Since $\T(\mathbb S A[-2],\X) = 0$, it follows that $f[-1] \circ e = 0$, hence there exists a map $k \colon A[-1] \lxr E$ such that $f[-1] \circ k = e$. Then $\underline{f[-1]} \circ \underline{k} = \underline{e}$ and  $\underline{m} \circ \underline{k} = 1_{\underline{E}}$ since $\underline{e}$ is an epimorphism. Since $\underline{A[-1]}$ is indecomposable in $\T/\X$, its endomorphism algebra is local, so the endomorphism $\underline{k} \circ\underline{m}$ is invertible or the endomorphism $1_{\underline{A[-1]}} - \underline{k} \circ\underline{m}$ is invertible. In the first case $\underline{m}$ is invertible and then $\underline{f[-1]}$ is epic, hence $\underline{\mathbb S A[-1]} = 0$, i.e. $\mathbb S A[-1] \in \X$. Then $\mathbb S A[-1] = 0$, i.e. $A = 0$, since $\mathbb S A[-1]\in \X^{\bot}$, and this is impossible. In the second case $\underline{f[-1]} = 0$ and then  $\underline{g[-1]}$ is invertible. Since $\T(\X[-1],\mathbb S A[-2]) = \T(\X,\mathbb S A[-1]) = 0$ this implies that $A= 0$. This contradiction shows that $X = 0$, so the map $h \colon A \lxr \mathbb S(A)$ is invertible. Clearly then $(T)$ is isomorphic to (7.4).

\item[{\em Step 4}:] Since (7.4) is an Auslander-Reiten triangle in $\T$, it follows that for any indecomposable object $T$ in $\T$, not isomorphic to $A$, we have $\T(A,T) = 0 = \T(T,A)$. Then connectedness of $\T$ implies that $\T = \add A$ and therefore $\X = 0$ since $A$ lies in $\X^{\bot}[-1]$. This contradiction shows that $A$ lies in $\X$.  \qedhere
\end{enumerate}   
\end{proof}
\end{thm}

Theorem $7.3$ and Corollary $3.5$ imply the following result which describes the extreme cases in the connections between the pretriangulated category $\T/\X$, its abelian subcategory $\mathsf{A}(\X) = (\X[-1]\star\X)/\X$, and its triangulated subcategory $\mathsf{T}(\X) = \X^{\bot}[-1]/\X$; there are two degenerate cases pointing into opposite directions:

\begin{cor} Let $\T$ be a connected $2$-Calabi-Yau triangulated category over a field and let $\X$ be a  functorially finite rigid subcategory of $\T$.

\smallskip

\begin{minipage}{8cm}
$\mathsf{(I)}$ If $\X \neq 0$, then the following are equivalent:
\begin{enumerate}
\item $\X$ is $2$-cluster tilting.
\item $\mathsf{T}(\X) = 0$.
\item $\mathsf{A}(\X) = \T/\X$.
\item $\T/\X$ is abelian. 
\end{enumerate}
\end{minipage}
\begin{minipage}{5cm}
$\mathsf{(II)}$ The following are equivalent:
\begin{enumerate}
\item $\X = 0$.
\item $\mathsf{T}(\X) = \T/\X$. 
\item $\mathsf{A}(\X) = 0$.
\item $\T/\X$ is triangulated. 
\end{enumerate} 
\end{minipage}
\end{cor}

\begin{rem} Assume that $\T$ is connected and $2$-Calabi-Yau, and let $0\neq \X$ be functorially finite rigid. Then by Corollary $7.4$ it follows that $\X$ is $2$-cluster tilting if and only if there is an equivalence $\smod\T/\X \approx \smod \,(\smod\X)$. Taking into account the equivalence $(\X\lsmod)^{\op} \approx \smod\X$ induced by the duality with respect to the base field, it follows that $\X$ is $2$-cluster tilting if and only if there is an equivalence:
\[
\begin{CD}
\smod\T/\X \,\ @> \approx >> \,\ \smod \,(\X\lsmod)^{\op} 
\end{CD}
\] 
Hence $\X$ is $2$-cluster tilting if and only if $\smod\T/\X$ is the {\em free abelian} or {\em Auslander}, {\em  category} associated to $\X$, in the sense of \cite{B:freyd}. This means in particular that $\T/\X$ has enough projectives and injectives and: 
\[
\gd\smod\T/\X \,\ \leq \,\ 2 \,\  \leq \,\ \ddom\smod\T/\X
\]
where $\ddom\smod\T/\X$ denotes the {\em dominant dimension} of $\T/\X$, i.e. the largest $n \in \mathbb N \cup \{\infty\}$ such that any projective object has an injective resolution whose first $n$ terms are projective, see \cite{B:freyd} for more details.
Further the functor $\mathsf{F}_{\X} \colon \X \lxr \smod\T/\X$, \ $\mathsf{F}_{\X}(X) = (-,\underline{X[-1]})$, realizes $\X$ as the full subcategory of projective-injective objects of $\smod\T/\X$, and is universal for functors out of $\X$ to abelian categories in the following sense: for any functor $G \colon \X \lxr \A$ to an abelian category $\A$, there exists a unique {\em exact} functor $G^{*} \colon \smod\T/\X \lxr \A$ such that  $G^{*} \circ \mathsf{F}_{\X} = G$. It follows that $\smod\X$ is of finite representation type if and only if so is $\T/\X$ and in this case, $\smod\T/\X$ is the category of finitely generated modules over the Auslander algebra of $\smod\X$.
\end{rem}

\,

{\bf Acknowledgement.}  The author thanks the referee for the useful comments and remarks.

\

\end{document}